\documentclass[leqno,twoside,11pt]{amsart}

\usepackage{latexsym,esint,url, amssymb, yfonts, enumitem}
\usepackage{graphicx,color}
\usepackage{floatrow}

\theoremstyle{plain}

\newtheorem{theorem}{Theorem}[section]
\newtheorem{corollary}[theorem]{Corollary}
\newtheorem{lemma}[theorem]{Lemma}

\theoremstyle{definition}

\newtheorem{remark}[theorem]{Remark}
\newcommand{\R}{\mathbb{R}} 

\numberwithin{equation}{section}
\numberwithin{figure}{section}

\begin{document}

\title[Geodesics and isometric immersions in kirigami] {Geodesics and isometric immersions in kirigami} 
\author{Qing Han, Marta Lewicka and L. Mahadevan}
\address{Qing Han: Department of Mathematics, University of Notre Dame, Notre Dame, IN 46556}
\address{Marta Lewicka: University of Pittsburgh, Department of
  Mathematics, 139 University Place, Pittsburgh, PA 15260}
\address{L. Mahadevan: Harvard University, School of Engineering and Applied Sciences, and Department of Physics,
29 Oxford Street, Cambridge, MA 02138}
\email{Qing.Han.7@nd.edu, lewicka@pitt.edu, lmahadev@g.harvard.edu} 

%\subjclass[2000]{}

%\keywords{}
%\date{April 10, 2021}

\begin{abstract} 
Kirigami is the art of cutting paper to make it articulated and
deployable, allowing for it to be shaped into complex two and
three-dimensional geometries. The mechanical response of a kirigami
sheet when it is pulled at its ends is enabled and limited by the
presence of cuts that serve to guide the possible
non-planar deformations. Inspired by the geometry of this art
form, we ask two questions: (i) What is the shortest path
between points at which forces are applied? (ii) What is the nature of
the ultimate shape of the sheet when it is strongly stretched?

Mathematically, these questions are related to the nature and form of
geodesics in the Euclidean plane with linear obstructions (cuts), and
the nature and form of isometric immersions of the sheet with cuts
when it can be folded on itself.  We provide a constructive proof
that the geodesic connecting any two points in the plane is piecewise
polygonal.  We then prove that the family of polygonal geodesics can
be simultaneously rectified into a straight line by flat-folding the sheet so that its
configuration is a (non-unique) piecewise affine planar isometric
immersion.  
\end{abstract}

\maketitle
%\tableofcontents

\section{Introduction}\label{intro}

\begin{figure}
\floatbox[{\capbeside\thisfloatsetup{capbesideposition={right,center},capbesidewidth=10cm}}]{figure}[\FBwidth]
{\caption{A circular sheet of paper with a cut in it becomes soft,
  because the cut allows the sheet to buckle out of the plane when
  pulled by two equal and opposite forces at its boundary. (i) The red
  lines are the geodesics connecting the points of force application.
  (ii) As the sheet deforms, it creates conical structures that allow
  the sheet to deform further, thus causing the geodesics to
  straighten out. (iii) The ultimate shape of the sheet in the
  strongly deformed limit causes each polygonal geodesic to straighten
  out and yields a flat-folded sheet that is piecewise affine
  isometric to the plane,  accompanied by set of sharp
  folds. [We thank G. Chaudhary for the photographs.] }\label{Fig00}}
{\includegraphics[scale=0.5]{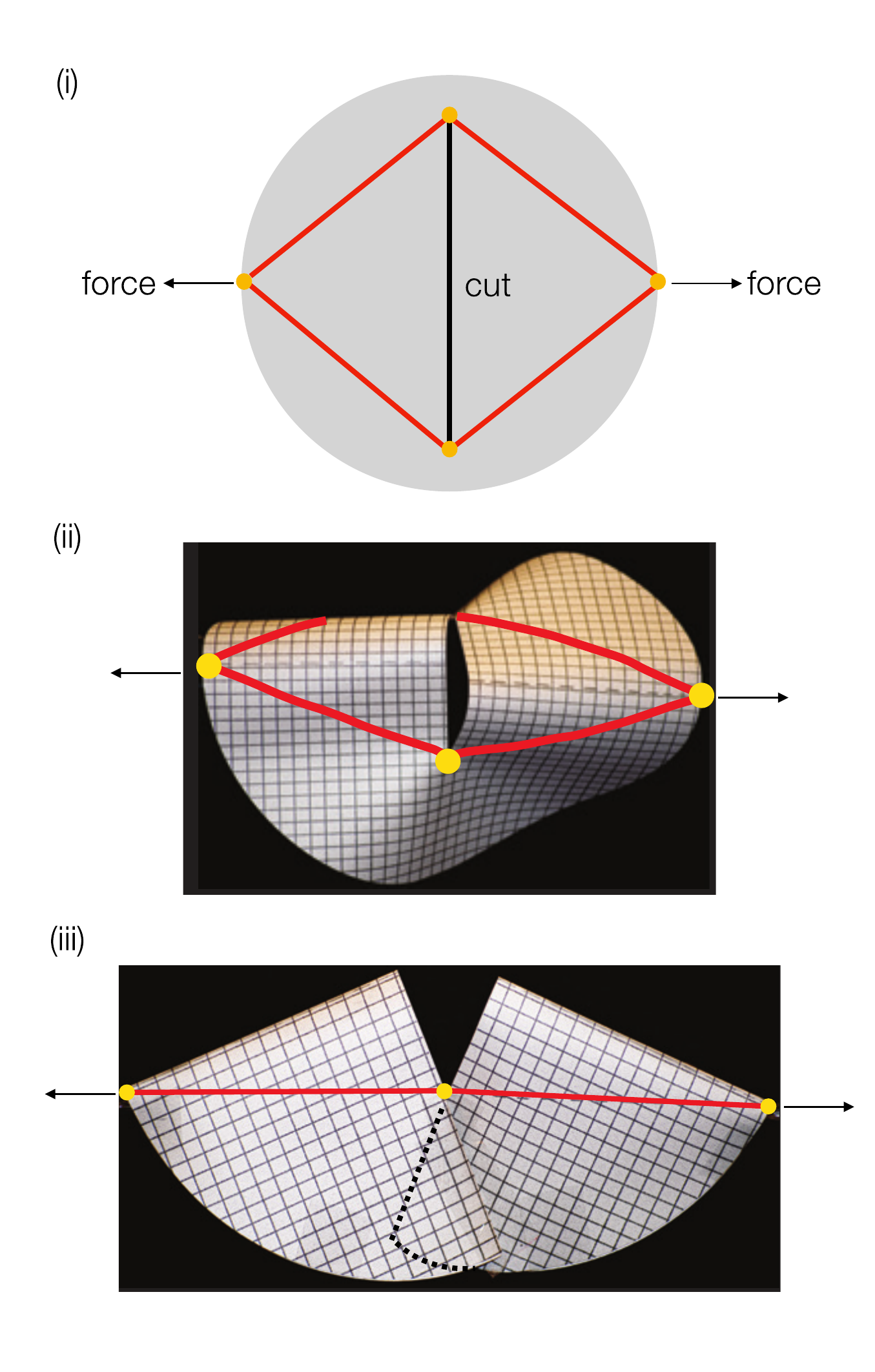}}
\end{figure}

A thin rectangular sheet of paper pulled at its corners is almost
impossible to stretch. Introducing a cut in its interior
changes its topology, and thence changes its physical response. The
corners can now be pulled apart as the sheet  bends out of the
plane, see Figure \ref{Fig00}. The physical reason for this is that the geometric
scale-separation associated  with a sheet of thickness $h$ and size
$L$ (where $h \ll L$), makes it energetically expensive to stretch and easy
to bend, since the elastic potential energy of the sheet per unit area can be written as: 
$$U= Eh ({\rm stretching ~ strain})^2 + Eh^3 ({\rm
  curvature})^2,$$ 
where the {\em stretching strain} and {\em curvature}
characterize the modes of deformation of the sheet, and $E$ is the
elastic modulus of the material. Thus, as $h/L
\rightarrow 0$, for given boundary conditions it is energetically
cheaper to deform by bending (curving) rather than stretching, as can
be observed readily with any thin sheet of any material. This
observation and its generalizations are behind the Sino-Japanese art
of kirigami (kiri = cut, gami = paper). Recently, this ability to make
cuts in a sheet of paper that allow it to be articulated and deployed
into complex two and three-dimensional patterns has become the
inspiration for a new class of mechanical metamaterials
\cite{Kirigami-review, Metamaterial-review}. The
geometrical and topological properties of the slender sheet-like
structures, irrespective of their material constituents, can then be exploited to create
functional structures on scales ranging from the nanometric
\cite{Mceuen} to centimetric and beyond \cite{Bertoldi, M1, M2}.  

Of the various mathematical and physical questions that arise from
this ability to control the configurational degrees of freedom of the
sheet using the geometry and topology of the cuts,  perhaps the
simplest is the following: if a sheet with random cuts was pulled at
two points on the boundary, what is the nature of paths of stress
transmission through the sheet? In the absence of cuts, the lines of
force transmission are straight lines connecting the points,
i.e. geodesics, but this needs to be revisited in the presence of
obstructing cuts. One might ask about the nature of the paths of force
transmission, i.e. the geodesics in this situation. The results of
qualitative experiments with a sheet of paper that has a single cut
along the perpendicular bisector to the line joining the
points of force application, are shown in Figure \ref{Fig00}. For small
forces, the sheet deforms into two conical regions that allow the edge
of the cut to curve out of the plane, and when the forces are large
enough, the ends of the cut become approximately collinear with the
line joining the points of forcing.  Observations of sheets with
multiple cuts are suggestive of a generalization, namely that
cuts cause the sheet to buckle out of the plane until a straight
geodesic in $\mathbb{R}^3$ connects the points of force
application. Furthermore, as the sheet thickness becomes vanishingly
small, allowing the sheet to form sharp creases with a large
curvature, the sheet can fold on itself and become flat again, as seen
in Figure \ref{Fig00} (iii).  

\smallskip

These observations suggest two conjectures:
\begin{itemize} 
\item[(i)] geodesics in a planar sheet with cuts are piecewise linear, i.e. they are polygonals; 
\item[(ii)] on pulling at two points in a sheet with cuts, these polygonal geodesics straighten out
by allowing the sheet to deform in the third dimension, which when
flat-folded causes the geodesic to be rectified, leading to a
configuration that is a piecewise affine isometric immersion.% of the plane. 
\end{itemize}
Here, we prove the above two statements.

\medskip

We point out that a combination of physical and
numerical experiments can be used to characterize the geometric
mechanics of kirigamized sheets as a function of the number, size, and
orientation of cuts. This will be the topic of our forthcoming work
\cite{Kirigami-jointpaper}, which in particular shows that by varying
the geodesic lengths, one can shape the deployment trajectory of a sheet as a composition of
developable units: flats, cylinders and cones, as well as
control its compliance across orders of magnitude.

\medskip

\noindent {\bf Acknowledgement.}  Marta Lewicka was partially supported by the
NSF grant DMS-2006439, L Mahadevan was partially supported by NSF
grants BioMatter DMR 1922321 and MRSEC DMR 2011754 and EFRI 1830901.

\section{The set-up and the main results of this paper}

Let $\Omega\subset \R^2$ be a convex, bounded, planar domain and let
$L$ be the union of finitely many closed segments 
contained in $\Omega$. We study the geodesic distance and the
structure of geodesics in $\Omega\setminus L$. 

% These  immersions are constructed as continuous, piecewise affine maps --- hence
% TODO KIRIGAMI.

\medskip

Specifically, we work under the following setup:
\begin{equation}\tag{S}\label{S}
\left[~\mbox{\begin{minipage}{15.3cm}
In an open, bounded, convex set $\Omega\subset \R^2$, given is a
graph $G$, consisting of $\bar n\geq 2$ vertices
$V=\{a_i\}_{i=1}^{\bar n}$ and $n\geq 1$ edges $E=\{l_i\}_{i=1}^n$, represented by:
$$ l_i = \big\{(1-t)a_{j}+ta_{k};~ t\in [0,1]\big\} \mbox{ for some
$a_j \neq a_k\in V$.} $$
We denote $L=\bigcup_{i=1}^n l_i$ and call $L$ the set of {\em
  cuts}. Without loss of generality, we further assume that $G$ is a
planar graph, i.e. for all $i\neq j$ the intersection $l_i\cap l_j$
is either empty or consists of a single point that is a common vertex of $l_i$ and $l_j$.
\end{minipage}}\right.
\end{equation}

\medskip

\noindent The collection of cuts in $L$ is thus finite but completely arbitrary, i.e. the
cuts may have any length and orientation, and are allowed to intersect
each other. 

\medskip

We next define: 
\begin{equation}\tag{G}\label{G}
\left[~\mbox{\begin{minipage}{15.3cm}
Given two points $p\neq q$ that belong to the same connected component
of $\bar \Omega\setminus L$, we set:
$$\mathrm{dist}(p,q) = \inf\big\{length(\tau);~ \tau:
  [0,1]\to\bar\Omega\setminus L \mbox{ piecewise $\mathcal{C}^1$ with $\tau(0)=p$, $\tau(1)=q$} \big\}.$$
Further, any piecewise $\mathcal{C}^1$ curve
$\sigma:[0,1]\to\bar\Omega$ with $\sigma(0) = p$, $\sigma(1)=q$ is called
a {\em geodesic from $p$ to $q$ in $\Omega\setminus L$}, provided that:
\begin{itemize}
\item[(i)]  $length(\sigma) = \mathrm{dist}(p,q),$
\item[(ii)] $\sigma$ is the uniform limit as $k\to\infty$, of a sequence of
  piecewise $\mathcal{C}^1$ curves  $\{\tau_k:[0,1]\to\bar\Omega\setminus
  L\}_{k=1}^\infty$, each satisfying $\tau_k(0)=p$, $\tau_k(1)=q$.
\end{itemize}
\end{minipage}}\right.
\end{equation}

\medskip

\begin{figure}[htbp]
\centering
\includegraphics[scale=0.6]{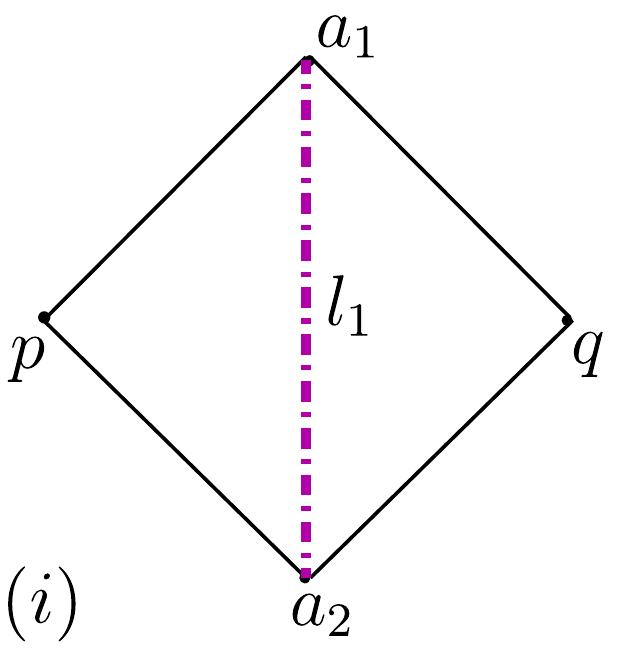}\qquad 
\quad \includegraphics[scale=0.6]{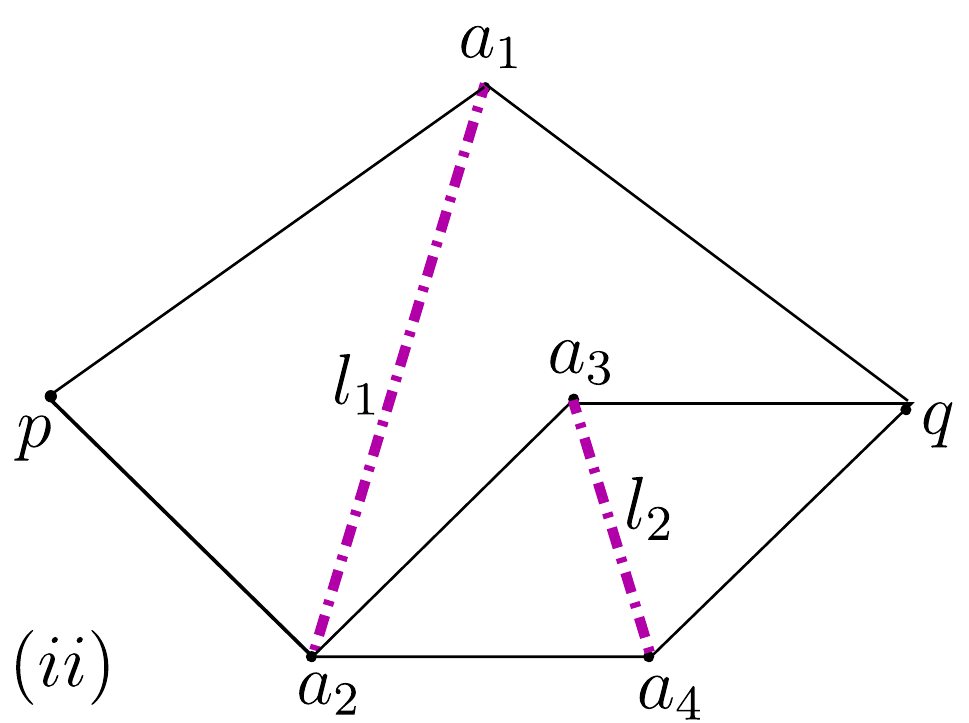}
\qquad \quad\includegraphics[scale=0.6]{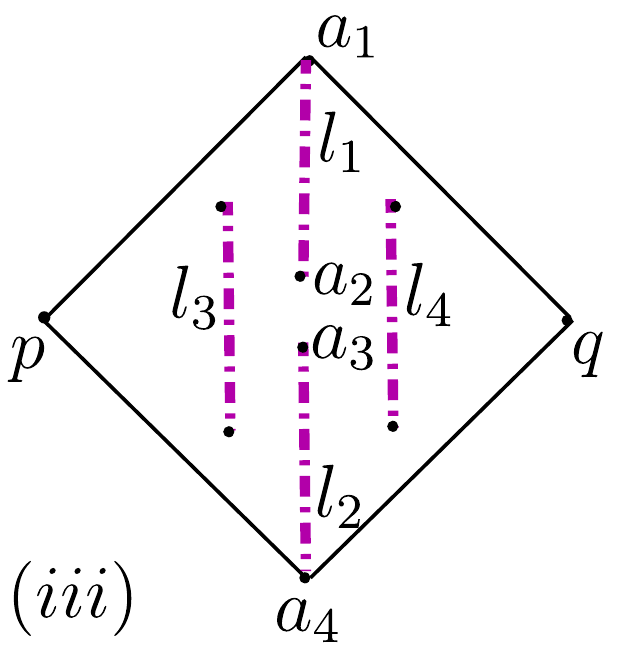}
\caption{{Three configurations of $G, p, q$ with pairwise nonintersecting cuts: (i) $p=(0,-1), q=(1,0)$,
   $a_1=(0,1), a_2=(0,-1)$ yield two geodesics:
    $\sigma_1=\overline{pa_1q}, \sigma_2=\overline{pa_2q}$; 
(ii) $p=(0,0), q=(4,0)$, $a_1= (2,(2^{3/2}-1)^{1/2}),  a_2=(1,-1),
a_3=(3,-1), a_4=(4,0)$ yield three geodesics:
    $\sigma_1=\overline{pa_1q}, \sigma_2=\overline{pa_2a_3q}, \sigma_3=\overline{pa_2a_4q}$; 
(iii) $p=(-1,0), q=(1,0)$, $a_1= (0,1),  a_4=(0,-1)$ with
$a_2, a_3=(0,\pm\epsilon)$, $a_5, a_6, a_7, a_8 = (\pm\epsilon,
\pm(1-{\epsilon}^{1/2}))$ for a sufficiently small $\epsilon>0$ yield two geodesics:
    $\sigma_1=\overline{pa_1q}, \sigma_2=\overline{pa_4q}$.}}
\label{Fig1}
\end{figure}

\begin{figure}[htbp]
\centering
\includegraphics[scale=0.6]{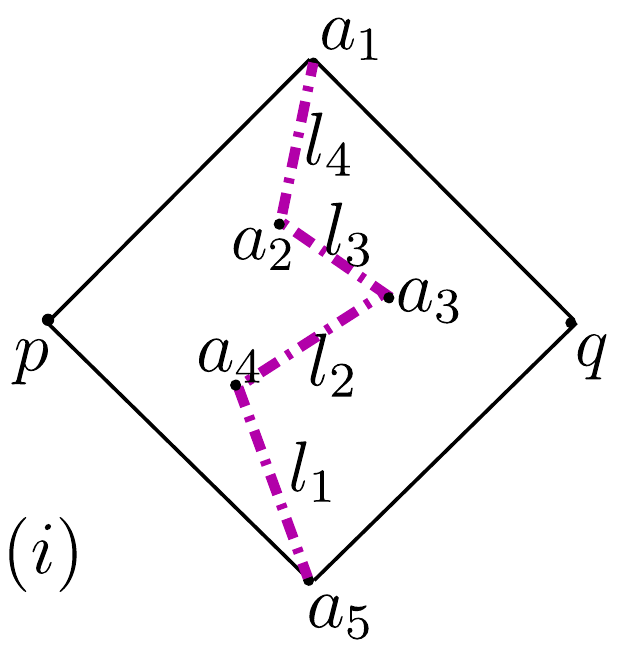}\qquad 
\quad \includegraphics[scale=0.6]{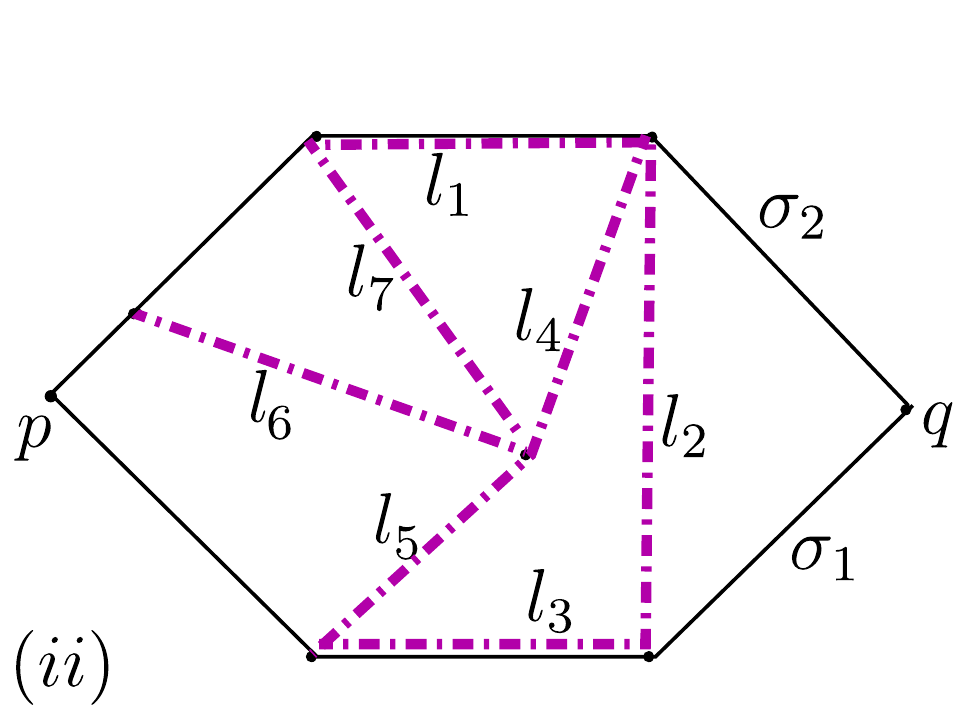} \qquad \quad 
\includegraphics[scale=0.6]{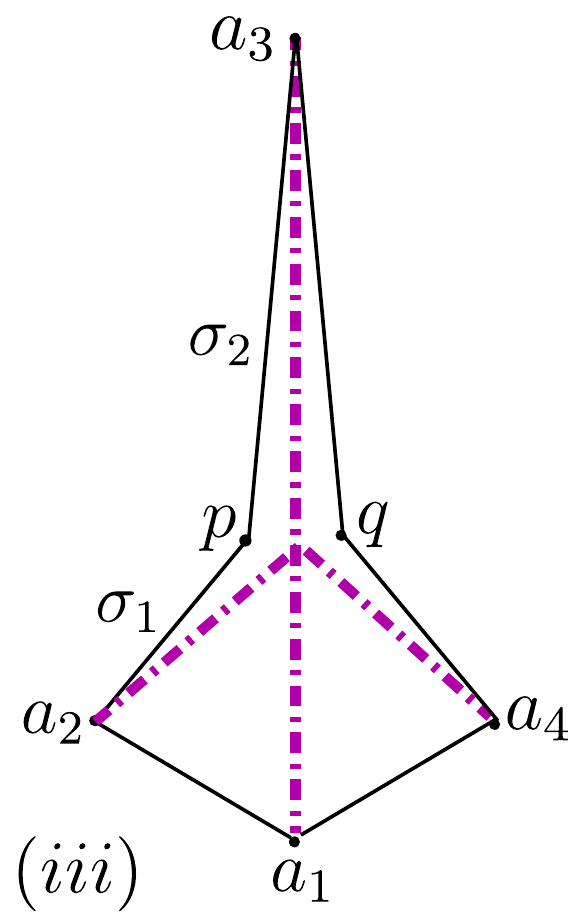}
\caption{{Three configurations of $G, p, q$ with intersecting cuts:
    (i) $n=4$, $\bar n=5$ result in two geodesics:
    $\sigma_1=\overline{pa_1q}, \sigma_2=\overline{pa_5q}$, this
    configuration is minimal in the sense introduced in Section \ref{sec21}; 
(ii) $n=7$, $\bar n=6$ and two geodesics: $\sigma_1$,
$\sigma_2$; (iii) $n=4$, $\bar n=5$ with $length(\overline{a_1a_2p}) =
length(\overline{pa_3})$, this is also a minimal configuration
resulting in two geodesics $\sigma_1, \sigma_2$.}}
\label{Fig1.5}
\end{figure}

\noindent The above definition abuses the notion of a geodesic
slightly, because it  allows $\sigma$ to be  not entirely contained in
$\Omega\setminus L$ (although we call it a geodesic in 
$\Omega\setminus L$), that is,  we allow cuts to be parts of $\sigma$.
Figures \ref{Fig1} and \ref{Fig1.5} show a few examples of $G, p, q$
and the resulting geodesics.
% as in the definition (\ref{G}). 
Here and below, by $\overline{p a_{i_1}
  a_{i_2}\ldots a_{i_k}q}$ we denote the polygonal joining $p$ and $q$ through the consecutive points $a_{i_1},
a_{i_2},\ldots, a_{i_k}$. % in the indicated order. 

\medskip

Our first result is as follows:

%\smallskip

\begin{theorem}\label{thm1}
Assume (\ref{S}) and let $p,q\in \bar\Omega\setminus L$ belong to one
connected component of $\bar\Omega\setminus L$. Then there exists
at least one geodesic from $p$ to $q$, as defined in (\ref{G}). Each
such geodesic $\sigma$ satisfies:
\begin{itemize}
\item[(i)] $\sigma$ is a finite polygonal joining 
$p$ and $q$, with all its other vertices distinct and chosen from $V$,
\item[(ii)] for each $i=1\ldots n$, if $\sigma\cap l_i\neq\emptyset$ then either $l_i\subset
  \sigma$ or $\sigma\cap l_i\subset \{a_{j}, a_{k}\}$, where $l_i=\overline{a_ja_k}$.
\end{itemize}
\end{theorem}

\medskip

\begin{figure}[htbp]
\centering
\includegraphics[scale=0.28]{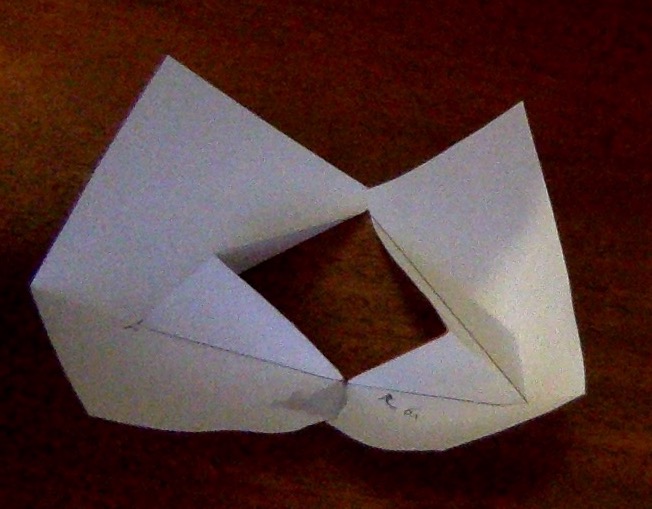}\qquad 
\quad \includegraphics[scale=0.338]{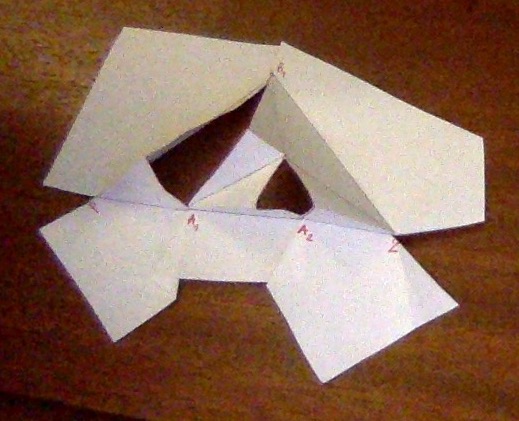}
\caption{{Examples of isometric immersions with properties as in Theorem
    \ref{thm2}, for configurations $V, L, p,q$ as in Figure \ref{Fig1} (i) and (ii).}}
\label{Fig2}
\end{figure}

Our second result and the main contribution of this paper, is motivated by the general considerations in
Section \ref{intro}. We prove the existence of an isometric immersion of $\Omega\setminus L$ into
$\R^3$, which bijectively maps each geodesic between two chosen
boundary points $p,q$ onto one segment in $\R^3$ of appropriate
length (see Figure \ref{Fig2}). More precisely, we have:

\begin{theorem}\label{thm2}
Assume (\ref{S}) and let $p,q\in \partial\Omega$. Then, there exists
a continuous and piecewise affine map $u:\bar\Omega\setminus L\to\R^3$
with the following properties:
\begin{itemize}
\item[(i)]  $u$ is an isometry, i.e.: $(\nabla u)^T\nabla u =Id_2$
  almost everywhere in $\Omega\setminus L$,
\item[(ii)] the image $u(\sigma)$ of every geodesic $\sigma$ from $p$
  to $q$ in $\Omega\setminus L$, coincides with the segment $\overline{u(p) u(q)}$.
\end{itemize}
In particular, $|u(p)-u(q)|=length(\sigma)$ for each geodesic
$\sigma$ (as defined in (\ref{G})).
\end{theorem}

\medskip

We prove Theorem \ref{thm1} in section \ref{sec1} and Theorem
\ref{thm2} in sections \ref{sec21}-\ref{sec24}. Our proofs are
constructive and describe: a specific algorithm to find the polygonal geodesics in Theorem
\ref{thm1}, and a folding procedure that yields the isometric immersion $u$ in Theorem \ref{thm2}.
Even when all cuts in $L$ are non-intersecting (i.e. the edges of the
underlying graph $G$ are pairwise disjoint), the construction of $u$
is far from obvious. The general case requires a further refinement
of the previous arguments, because of the completely
arbitrary planar geometry of each connected component of $G$. 

\medskip

The algorithm that yields the isometry $u$ in Theorem \ref{thm2} consists of:

\smallskip

\begin{itemize}
\item[(i)] identifying and sealing the portions of {\em inessential cuts}, which do not affect $\mbox{dist}(p,q)$;
\item[(ii)] ordering the geodesics and ordering the remaining cuts,
  that now form a new planar graph $G$ consisting of {\em trees} (i.e. $G$ is a forest);
\item[(iii)] constructing $u$ on each region between two consecutive
   trees and two consecutive geodesics;
\item[(iv)] constructing $u$ on regions within each tree;
\item[(v)] constructing $u$ on the {\em exterior region} that is not enclosed by any two geodesics.
\end{itemize}

\smallskip

\noindent The points (i) and (ii) above are introduced in sections \ref{sec21}
and \ref{sec22}, respectively. The main arguments towards (iii) in the simplified setting are presented in section
\ref{sec23}. The general case is resolved  in section \ref{sec23.2}, which carries the heaviest
technical load of this paper. Section \ref{sec24} completes the proofs
and presents an example explaining the necessity of $p,q$ being located on
the boundary of $\Omega$ in Theorem \ref{thm2}.

\section{Proof of Theorem \ref{thm1}}\label{sec1}

Given $p,q\in\bar\Omega\setminus L$ and a piecewise $\mathcal{C}^1$ curve
$\tau:[0,1]\to\bar\Omega\setminus L$ with $\tau(0)=p$, $\tau(1)=q$, we
first demonstrate a general procedure to produce a finite polygonal
$\sigma$ which joins $p$ and $q$, whose other vertices are (not
necessarily distinct) points in $V$, which satisfies condition (ii) in (\ref{G}), and such that:
$$length(\sigma)\leq length(\tau).$$ 
Applying this procedure to curves $\tau$ with $length(\tau)\leq
\mbox{dist}(p,q)+1$ yields a family of polygonals with the listed
properties, each of them having number of edges bounded by:
$$\frac{\mbox{dist}(p,q)+1}{\min_{a_j\neq a_k}|a_j-a_k|}.$$
Hence, all geodesics from $p$ to $q$ in
$\Omega\setminus L$ are precisely the length-minimizing polygonals
among such (finitely many) polygonals. We further show that any
length-minimizing polygonal satisfying condition (ii) of (\ref{G})
cannot pass through the same vertex in $V$ multiple times. Theorem
\ref{thm1} is then a direct consequence of these statements.

\medskip

Without loss of generality, the path $\tau$ has no self-intersections. We
construct $\sigma$ by successive replacements of portions of
$\tau$ by segments, as follows:

\smallskip

\begin{equation*}
\left[\quad\mbox{\begin{minipage}{15.5cm}
{\bf 1.} It is easy to show that for all $t>0$ sufficiently small there holds:
$\overline{p\tau(t)}\subset\bar\Omega\setminus L$. Set $\tau_1=\tau$ and define:
\begin{equation*}
t_1 = \sup\big\{t\in (0,1);~
\overline{p\tau_1(s)}\subset\bar\Omega\setminus L\;\mbox{ for all } \;s\in (0,t)\big\},\quad q_1=\tau_1(t_1).
\end{equation*}
There further holds: $t_1\in (0,1]$ and $\overline{pq_1}$ is a geodesic from $p$ to $q_1$. If $q_1=q$ then we set
$\sigma=\overline{pq}$ and stop the process.
Otherwise, by construction, the segment $\overline{pq_1}$ must contain some of the
vertices in $V$. Call $p_1$ the closest one of these points to $q_1$ and note that $p_1\neq q_1$. Consider the
concatenation of the segment $\overline{p_1q_1}$ and 
the curve ${\tau_1}_{\mid [t_1,1]}$. After re-parametrisation, it yields a
piecewise $\mathcal{C}^1$ curve $\tau_2:[0,1]\to\bar\Omega$,
with the property that $\tau_2((0,1])\subset \bar\Omega\setminus L$ and also:
$$|p-p_1| + length(\tau_2) \leq length(\tau).$$

\medskip

{\bf 2.} We inductively define a finite sequence of endpoints
$\{p_i\}_{i=2}^k\subset V$ and a sequence of
piecewise $\mathcal{C}^1$ curves $\{\tau_i:[0,1]\to\bar\Omega\}_{i=3}^{k+1}$,
by applying the procedure in Step 1 to curve $\tau_{i}$ and points 
$\tau_{i}(t_{i})=p_i$ and $q$, until $q_{k+1}=q$ so that $\overline{p_kq}$ is a geodesic from $p_k$
to $q$. Along the way, we get: $\tau_i(0)=p_i\neq p_{i-1}$,
$\tau_i(1)=q$, $\tau_i((0,1])\subset \bar\Omega\setminus L$,  and:
$$ \mbox{the sequence }\; \Big\{|p-p_1|+\sum_{j=2}^{j=i}|p_j-p_{j-1}|
+ length(\tau_{i+1})\Big\}_{i=1}^k \; \mbox{ is non-increasing.}$$
Also, the subset of $\bar\Omega$ enclosed by the concatenation of
$\overline{pp_1\ldots p_iq_i}$ with the portion of the curve $\tau$ between $p$
and $q_i$, contains no cuts in its interior. Consequently, each
polygonal $\overline{pp_1\ldots p_iq_i}$ is a uniform limit of
$\mathcal{C}^1$ curves contained in $\bar\Omega\setminus L$.

\bigskip

{\bf 3.} We finally define:
$\sigma = \overline{pp_1\ldots p_kq}.$
\end{minipage}}\right.
\end{equation*}

\smallskip

\begin{figure}[htbp]
\centering
\includegraphics[scale=0.6]{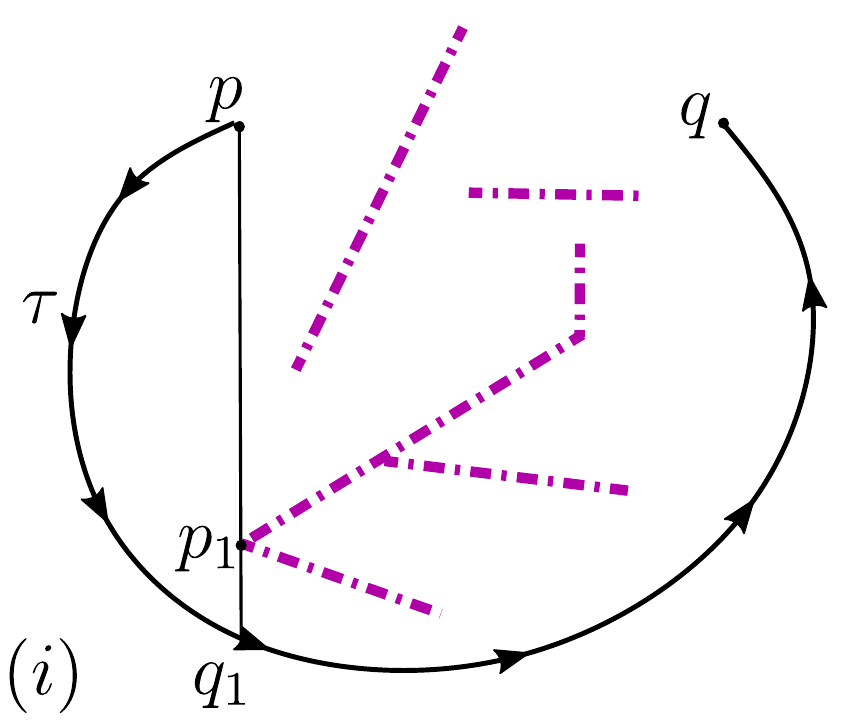}\qquad 
\quad \includegraphics[scale=0.6]{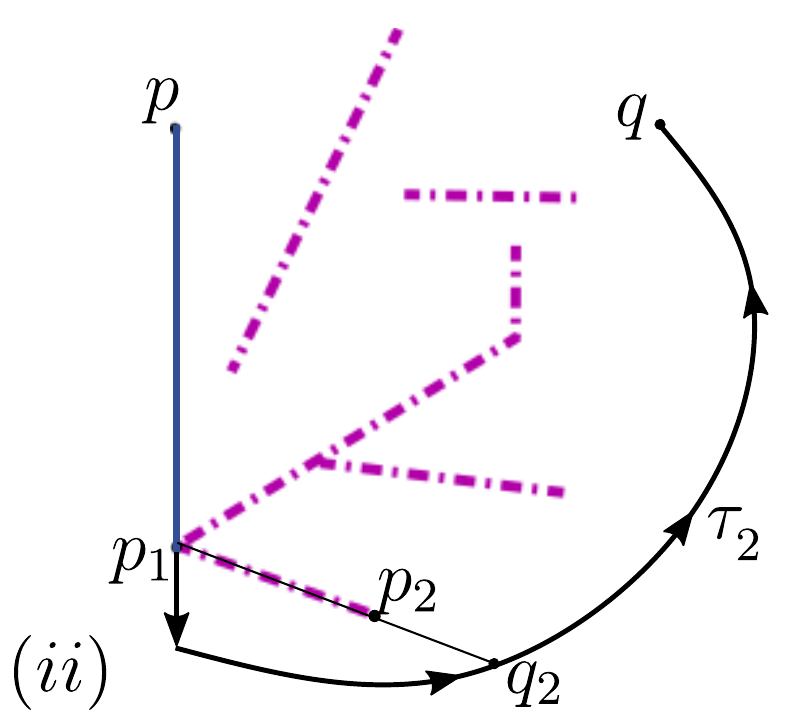} \vspace{0.8cm} \\
 \includegraphics[scale=0.6]{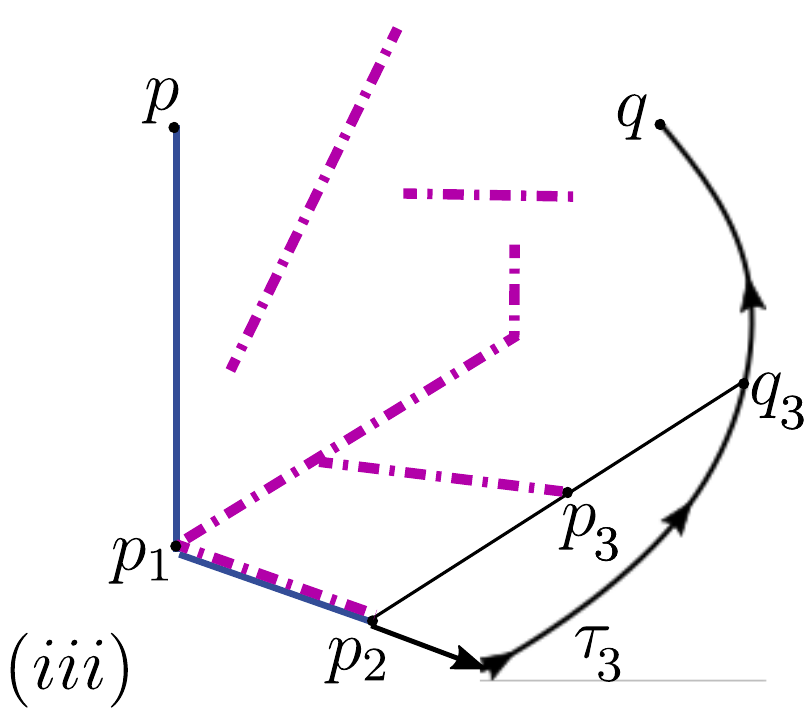} \qquad 
\quad \includegraphics[scale=0.6]{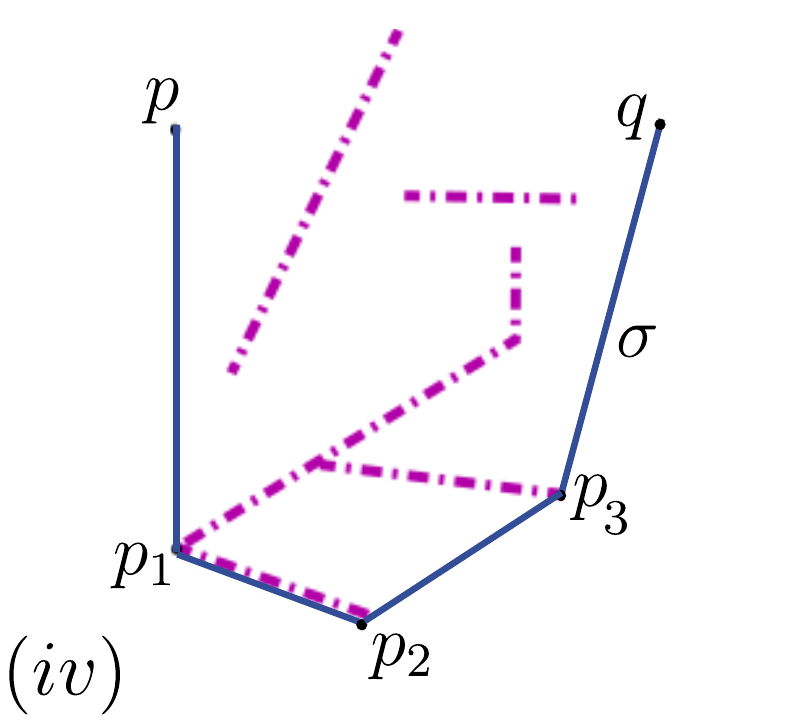}
\caption{{The path-shortening algorithm in the proof of Theorem \ref{thm1}.}}
\label{Fig3}
\end{figure}

\noindent The above process indeed terminates in a finite number $k$
of steps, because the length of each polygonal $\overline{pp_1\ldots
  p_i}$ is bounded by $length(\tau)$, and at each step this length
increases by at least: $\min_{a_j\neq a_k} |a_j-a_k|>0$.  
See Figure \ref{Fig3} for an example of $L, p, q,
\tau$ and the resulting polygonal $\sigma$.

\bigskip

The following observation concludes the proof of Theorem \ref{thm1}:

\begin{lemma}\label{lem_passonce}
Let $\sigma$ be a geodesic from $p$ to $q$ in $\Omega\setminus L$, as
in (\ref{G}). Then for every $a_i\in V$, there 
holds $a_i=\sigma(t)$ for at most one $t\in (0,1)$.
\end{lemma}
\begin{proof}
We argue by contradiction and assume that a geodesic polygonal
$\sigma$ passes through some vertex $a_i\in V$ at
least twice. Without loss of generality, we take $a_i$ to be the first
vertex in $\sigma$ (counting from $p$) with this property. Consider
the portion of $\sigma$ containing the first and second occurrences of
$a_i$, namely: $\overline{a_{i_0}a_ia_{i_1}a_{i_2}\ldots
  a_{i_s}a_ia_{i_{s+1}}}$, and consider the approximating curve $\tau$
as in definition (\ref{G}). From the approximate length-minimizing
property of $\tau$, it
follows that both angles $\angle(a_{i_0} a_i a_{i_1})$ and $\angle(a_{i_s} a_i a_{i_{s+1}})$
must be at least $\pi$. Consequently, they are both equal to
$\pi$. Another application of the same minimality condition yields that at
least one of the segments $\overline{a_{i_0}a_i}$ and
$\overline{a_ia_{i_{s+1}}}$ must be a cut in $L$. This contradicts
with $p,q\not\in L$ and $a_i$ being the first multiple vertex of $\sigma$.
\end{proof}

\section{Proof of Theorem \ref{thm2}. Step 1: sealing the inessential cuts}\label{sec21}

Assume (\ref{S}) and let $p,q$ be two distinct points belonging to one
connected component of $\bar\Omega\setminus L$. We describe a procedure which ``seals''
portions of cuts in $L$ without decreasing the geodesic distance
between $p$ and $q$ in $\Omega\setminus L$.
First, consider $i=1\ldots n$ and $j, k=1\ldots \bar n$ so that $l_{i} = \overline{a_ja_k}$.
Given $t\in \big[0, length(l_i)\big]$, define the altered endpoint of the cut $l_{i}$:
$$a_j(t)= (1-t)a_j + ta_k.$$
Let $L(t)$ be  the new set of cuts in which $l_{i}$ has been replaced by $l_{i}(t) = 
\overline{a_j(t)a_k}$, while all other cuts are left unchanged (this
construction alters the $V$ of the underlying graph $G$ as well). 
We have the following observation:

\begin{lemma}\label{lem1}
With the above notation, the geodesic distance between $p$ and $q$ in
$\Omega\setminus L(t)$:
$$t\mapsto \mathrm{dist}_t(p,q)=\inf\big\{length(\tau);~
\tau:[0,1]\to\bar\Omega\setminus L(t) \mbox{ piecewise $\mathcal{C}^1$,
  with } \tau(0)=p, \tau(1)=q\big\}$$
is left-continuous as a function of $t\in [0,1]$, and right-continuous in $t\in
(0,1]$. It is also right-continuous at $t=0$ when $a_j$ is not the
end-point of any other cut in $L$ besides $l_i$.
\end{lemma}
\begin{proof}
{\em Step 1.} To prove the asserted left-continuity, take a sequence
$\{t_m\in (0,1)\}_{m=1}^\infty$ that is strictly increasing to some
$t_0>0$. It is clear that $\mbox{dist}_{t_0}(p,q)\leq
\liminf_{m\to\infty}\mbox{dist}_{t_m}(p,q)$, because $\bar\Omega\setminus
L(t_0)\subset \bar\Omega\setminus L(t_m)$ so that
$\mbox{dist}_{t_0}(p,q)\leq \mbox{dist}_{t_m}(p,q)$ for all $m$.

For the reverse bound, fix $\epsilon>0$ and let $\tau:[0,1]\to\bar\Omega\setminus L(t_0)$
be piecewise $\mathcal{C}^1$ with $\tau(0)=p$, $\tau(1)=q$, and such
that $length(\tau)\leq \mathrm{dist}_{t_0}(p,q)+\epsilon$. 
We observe that if $\tau$ intersects $L(t_m)$, it must do so within
$l_i(t_m)\setminus l_i(t_0)$. Since for sufficiently large $m$ there holds: 
$length(l_i(t_m)) - length(l_i(t_0))<\epsilon,$
it follows that there exists $\tau_\epsilon:[0,1]\to\bar\Omega\setminus L(t_m)$ which is a local
modification of $\tau$, increasing its length by at most $2\epsilon$. Here, we are
taking advantage of the fact that $a_j(t_m)$ is not the endpoint of
any other cut besides $l_i(t_m)$ in $L(t_m)$. Consequently, we get:
$$\mbox{dist}_{t_m}(p,q)\leq length(\tau_\epsilon) \leq \mbox{dist}_{t_0}(p,q)+3\epsilon.$$
Since $\epsilon>0$ is arbitrary, this implies: 
$\limsup_{m\to\infty}\mbox{dist}_{t_m}(p,q)\leq \mathrm{dist}_{t_0}(p,q)$.

\medskip

{\em Step 2.} To show right-continuity of the function $\mbox{dist}_t(p,q)$ at $t_0\in (0,1)$, let $\{t_m\in
(0,1)\}_{m=1}^\infty$ that is strictly decreasing to $t_0$.  As in Step 1, we get:
$\limsup_{m\to\infty}\mbox{dist}_{t_m}(p,q)\leq
\mathrm{dist}_{t_0}(p,q)$. In virtue of Theorem \ref{thm1}, for each $m$ there holds:
$$\mbox{dist}_{t_m}(p,q)=
length\big(\overline{pa_{i_{1,m}}(t_m)a_{i_{2,m}}(t_m)\ldots a_{i_{k(m),m}}(t_m)q}\big),$$
where for $s\neq i$ we set $a_s(t) = a_s$. Since the number of finite sequences 
of distinct indices chosen among $\{1\ldots n\}$ equals $\sum_{k=1}^n
k!$ and it is finite, it follows that at least
one of such sequences $(i_1, i_2\ldots i_k)$ represents the order of the vertices in a geodesic
polygonal as above, for infinitely many $t_m$-s. Passing to a
subsequence if necessary, we may thus write:
$$\mbox{dist}_{t_m}(p,q)= length\big(\overline{pa_{i_{1}}(t_m)a_{i_{2}}(t_m)\ldots
  a_{i_{k}}(t_m)q}\big)\quad \mbox{ for all $m$.}$$
We emphasize that at most one of the vertices changes
as $m\to\infty$ and all others remain fixed.  
Further, as $m\to\infty$, the geodesics $\overline{pa_{i_{1}}(t_m)a_{i_{2}}(t_m)\ldots
  a_{i_{k}}(t_m)q}$ converge to the polygonal $\sigma =
\overline{pa_{i_{1}}\ldots a_{i_{k}}q}$ that satisfies condition (ii)
of (\ref{G}). Consequently:
$$\lim_{m\to\infty} \mbox{dist}_{t_m}(p,q)=length(\sigma)\geq \mbox{dist}_{t_0}(p,q).$$
This concludes the proof of the lemma. The same argument is valid at
$t_0=0$ under the indicated condition on $a_j$.
\end{proof}

\begin{figure}[htbp]
\vspace{1cm}
\hspace{-11cm}\includegraphics[scale=0.6]{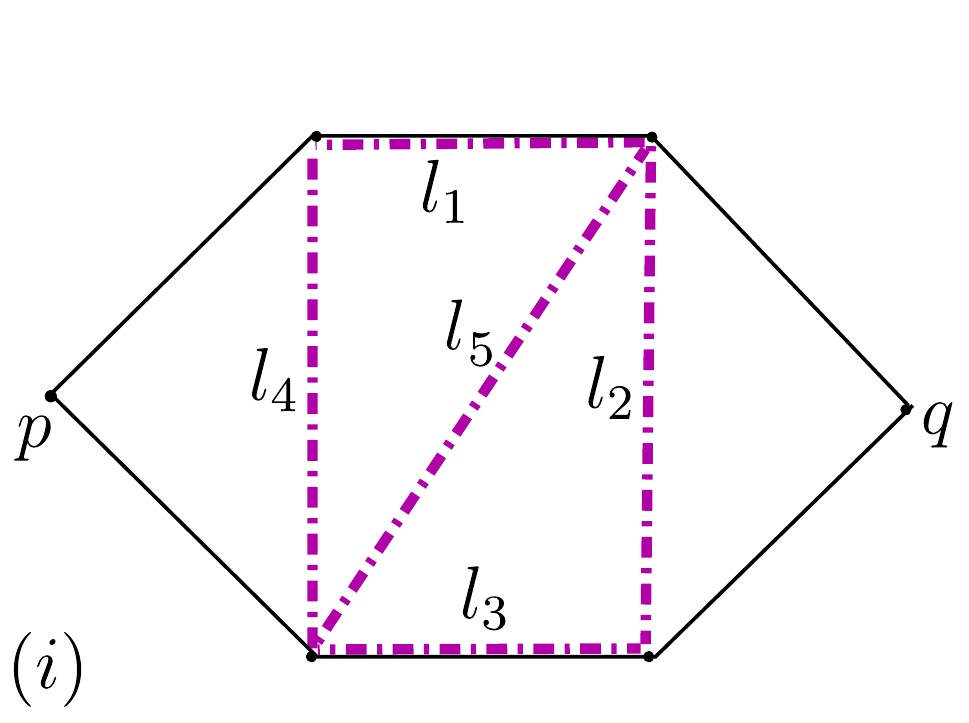}\\

\vspace{-6cm} \hspace{6.2cm}
\includegraphics[scale=0.5]{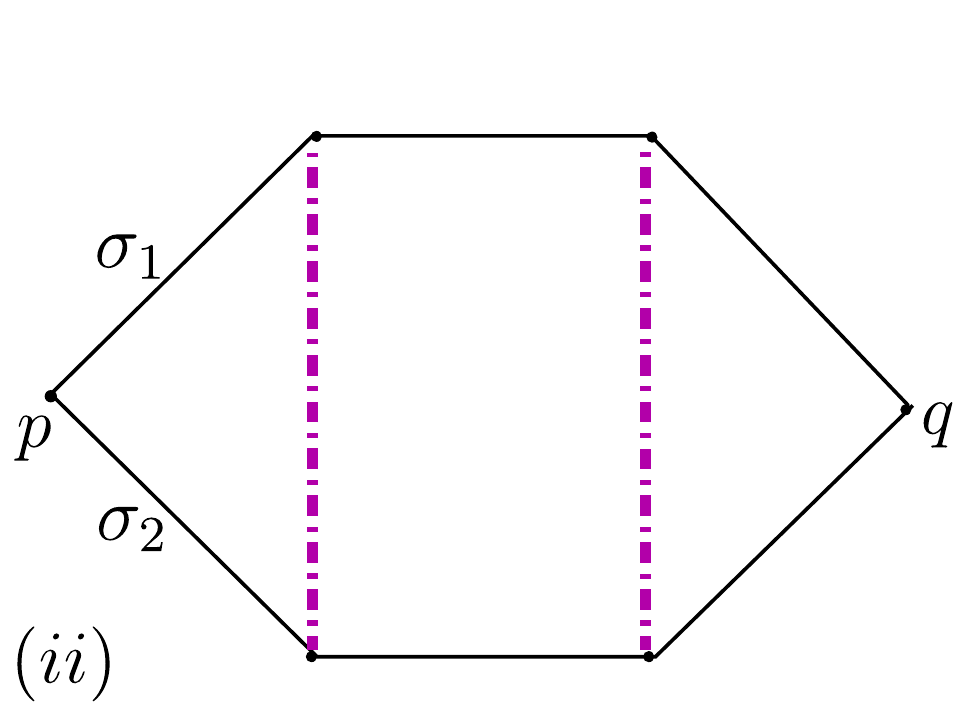} \quad \includegraphics[scale=0.5]{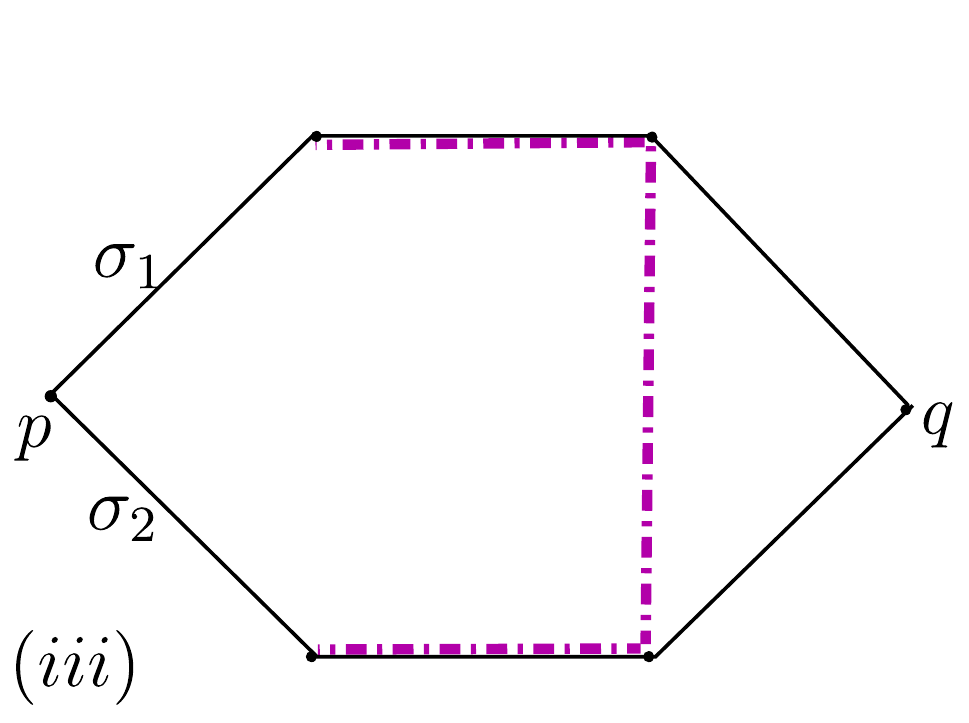} \\

\hspace{6.2cm}
\includegraphics[scale=0.5]{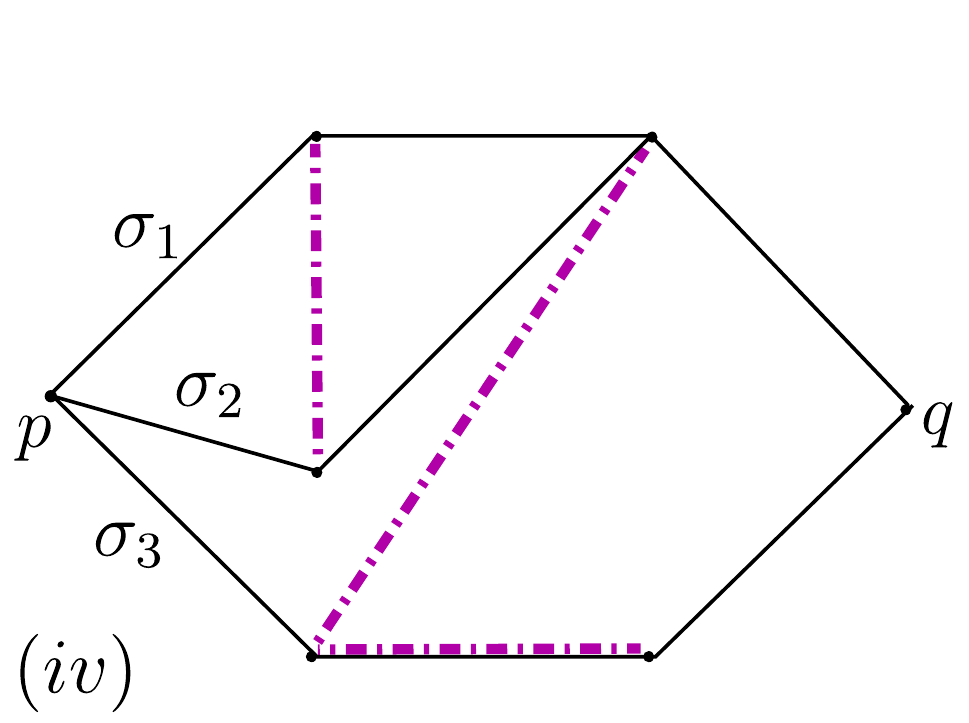} \quad \includegraphics[scale=0.5]{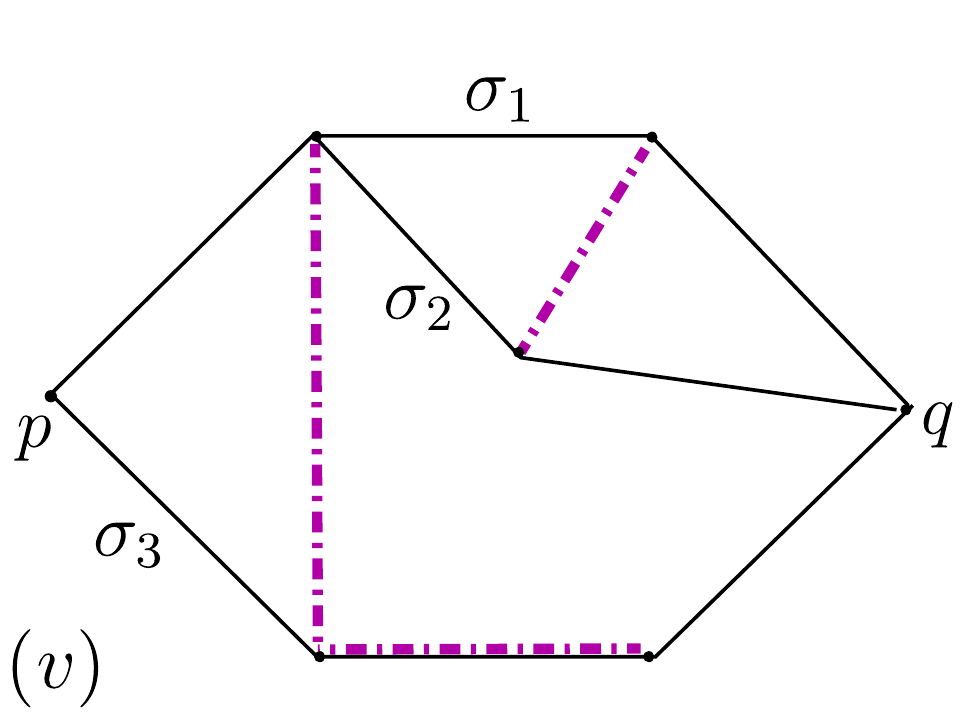}
\caption{{Different minimal configurations resulting from the original
  graph $G$ in (i), obtained by the sealing procedure upon changing
  the order of edges in $E$ and vertices in $V$: (ii) and (iii) yield
  two geodesics, while (iv) and (v) yield three geodesics. }}
\label{fig_minimal}
\end{figure}

\medskip

We now define an inductive procedure in which lengths of all cuts
are decreased as much as possible. Any resulting configuration $\bar
G, \bar V, \bar L$ (see Figure \ref{fig_minimal} for examples) will
be called {\em minimal}.

\begin{equation*}
\left[\quad\mbox{\begin{minipage}{15.5cm}
{\bf 1.} Fix $i=1$, write $l_1=\overline{a_ja_k}$ and define: 
$$t_1=\sup\big\{t\in [0, length(l_1)];~ \mbox{dist}_t(p,q)=\mbox{dist}_0(p,q)\big\},$$ 
where $\mathrm{dist}_t(p,q)$ is as in Lemma \ref{lem1}.
Replace the endpoint $a_j$ by $a_j(t_1)$, and replace the cut $l_1$ by
the segment $l_1(t_1) = \overline{a_j(t)a_k}$. If $a_j(t_1) = a_k$ then
we remove $l_1$ altogether. Consider the problem of finding an isometric immersion $u_1$ with properties
(i), (ii) in Theorem \ref{thm2}, for the same points $p,q$ but with $L$
replaced by $L_1=L(t_1)$. Then $u={u_1}_{\mid \bar\Omega\setminus L}$
is a continuous, piecewise affine map fulfilling Theorem \ref{thm2}.

\medskip

{\bf 2.} Write now $l_1=\overline{a_ka_j}$ and let $t_2$ be defined as
above, where we decrease the length of the already modified cut $l_1$
starting from the so far unaltered vertex $a_k$, up to
$a_k(t_2)$. Replace $l_1$ by $l_2(t_2)= \overline{a_k(t_2)a_j}$ or
remove it all together in case $a_k(t_2)=a_j$. Call the new set of cuts $L_2$.

\medskip

{\bf 3.} Having constructed $L_{2i}$ for some $1\leq i<n$, consider
the next cut $l_{i+1}=\overline{a_ja_k}$ and define: 
$$t_{2i+1} = \sup\Big\{t\in\big[0,length(l_{i+1})\big];~ \mbox{dist}_t(p,q) = \mbox{dist}_0(p,q)\Big\},$$
where $\mbox{dist}_t(p,q)$ is taken
with respect to the previously obtained set of cuts $L_{2i}$. Replace
the endpoint $a_{j}$ by $a_{j}(t_{2i+1})$ and replace the cut $l_{i+1}\subset L_{2i}$ by $l_{i+1}(t_{2i+1})$. This
defines the new collection of cuts $L_{2i+1}\subset L_{2i}$.

\medskip

{\bf 4.} In the same manner, by possibly modifying the endpoint $a_k$
of the already considered cut $l_{i+1}$, we construct the new set of cuts $L_{2i+2}\subset L_{2i+1}$.

\medskip

{\bf 5.} We finally set: 
$$\bar L = L_{2n}.$$ 
As in Step 1 of the algorithm, this ultimate collection $\bar L\subset
L$ of cuts in $\Omega$ has the property that the
validity of Theorem \ref{thm2} for the configuration $p,q,\bar L$ implies its validity
for the original configuration $p,q,L$.
\end{minipage}}\right.
\end{equation*}

\bigskip

\noindent Informally speaking, the above procedure starts by moving
the first endpoint vertex of $l_1$ toward its second vertex, 
whereas we start ``sealing'' the portion of the cut $l_1$ left
behind.  The length of the geodesics connecting $p$ and $q$ may drop
initially, in which case we leave the configuration unchanged. Otherwise,
the geodesic distance is continuously nonincreasing, in view of Lemma
\ref{lem1} (it may initially remain constant). We 
stop the sealing process when the aforementioned distance becomes
strictly less than the original one, and label the new position point
as the new vertex endpoint of $l_1$. In the next step, we move the
remaining endpoint along $l_1$ toward the (new) first endpoint and
repeat the process, thus possibly sealing the cut $l_1$ further. The
procedure is carried out for each $l_i$ in the given order $i=1, 2, \ldots, n$. 
We now claim that the distance between $p$ and $q$ cannot be further
decreased, upon repeating
the same process for the newly created configuration.  

\begin{lemma}\label{lem2}
With respect to the cuts in $\bar L = \bigcup_{i=1}^n \bar l_i$, for any
$i=1\ldots n$, any of the endpoint vertices of $\bar l_i$, and any
$t>0$ there holds:
$$\mathrm{dist}_t(p,q) <\mathrm{dist}_0(p,q).$$
\end{lemma}
\begin{proof}
Denote $d=\mbox{dist}_t(p,q)$ as above. At the $(2i-1)$-th step of 
construction of $\bar L$, we have:
$$d\leq \mbox{dist}_{t_{2i-1}+t}(p,q),$$
because $d$ corresponds to the geodesic distance between $p$ and $q$
in the complement of the cut set $\bar L$ with $\bar l_{i}$ further decreased, while $\mbox{dist}_{t_{2i-1}+t}(p,q)$
corresponds to the geodesic distance in the subset of the
aforementioned complement, obtained by enlarging all cuts $\{\bar
l_j\}_{j>i}$ to their original lengths in $L$. On
the other hand, directly by construction of $\bar L$ we get:
$$\mbox{dist}_{t_{2i-1}+t}(p,q)< \mbox{dist}_{t_{2i-1}}(p,q) = \mbox{dist}_0(p,q).$$
This ends the proof of the lemma.
\end{proof}

\begin{corollary}\label{cor3}
The set of cuts $\bar L$ constructed above coincides with the set of
edges $\bar E$ of the modified graph $\bar G$, with the new set of
vertices $\bar V$, which have the following properties:
\begin{itemize}
\item[(i)] $\bar G$ has no loops, and consequently it is a forest,
  consisting of finitely many trees,
\item[(ii)] each vertex in $\bar V$ that is an endpoint of only one
  edge in $\bar E$ (i.e. a leaf of the forest $\bar G$), is a vertex of some geodesic $\sigma$ from $p$ to
  $q$ in $\Omega\setminus \bar L$.
\end{itemize}
\end{corollary}
\begin{proof}
{\em Step 1.} To prove (i), we show that $\R^2\setminus \bar L$ must
be connected.  Indeed, in the opposite case, the boundary of the
connected component $R_1$ of $\R^2\setminus \bar L$ containing $p$ and
$q$, must contain a cut $\overline{a_ja_k}$ that is also a part of the
boundary of some other connected component $R_2$ of $\R^2\setminus
\bar L$. By the minimality property of $\bar L$ in Lemma \ref{lem2},
it follows that the sealing procedure with respect to the indicated
cut $\overline{a_ja_k}$ and its endpoint $a_{j}$ results in the
decrease of $\mbox{dist}_t(p,q)$ for any $t>0$. Consequently, there exists
a piecewise $\mathcal{C}^1$ curve $\tau:[0,1]\to\bar\Omega\setminus
\bar{L}(t)$ with $\tau(0)=p$, $\tau(1)=q$ and
$length(\tau)<\mbox{dist}_0(p,q)$, where this last distance is taken
in $\Omega\setminus \bar L$. The curve $\tau$ must both enter and exit
$R_2$ through the segment
$\overline{a_{j}a_{j}(t)}$. This means that $\tau$ may be further
shortened by replacing its portion contained in the aforementioned
interior region by an appropriate straight segment. The resulting
curve is $\bar\tau: [0,1]\to\Omega\setminus\bar L$, with:
$$length(\bar\tau)<length(\tau)<\mbox{dist}_0(p,q),$$
which is a contradiction.

\smallskip

{\em Step 2.} Let $a_i\in\bar V$ be as requested in (ii).  Consider
the modified endpoints $a_i(1/m)$ and the cut collections $\bar
L(1/m)$ as described in the sealing algorithm. As in the proof of
Lemma \ref{lem1}, there must exist a finite sequence $(i_1,\ldots,
i_k)$ such that:
$$\sigma_m = \overline{a_{i_1}(1/m) a_{i_2}(1/m)\ldots a_{i_k}(1/m) }$$
is a geodesic from $p$ to $q$ in $\Omega\setminus \bar L(1/m)$ for
infinitely many $m$-s. By the maximality assertion in Lemma
\ref{lem2}, there must be: $i\in\{i_1,\ldots, i_k\}$. But then Lemma \ref{lem1} yields:
$$\mbox{dist}_0(p,q)=\lim_{m\to\infty}\mbox{dist}_{1/m}(p,q) =
\lim_{m\to\infty} length(\sigma_m) = length(\overline{pa_{i_1}\ldots a_{i_k}q}).$$
Consequently, $\overline{pa_{i_1}\ldots a_{i_k}q}$ is a geodesic from $p$ to $q$ in $\Omega\setminus
  \bar L$, whose existence is claimed in (ii).
\end{proof}

\begin{figure}[htbp]
\hspace{-7.5cm}
\includegraphics[scale=0.6]{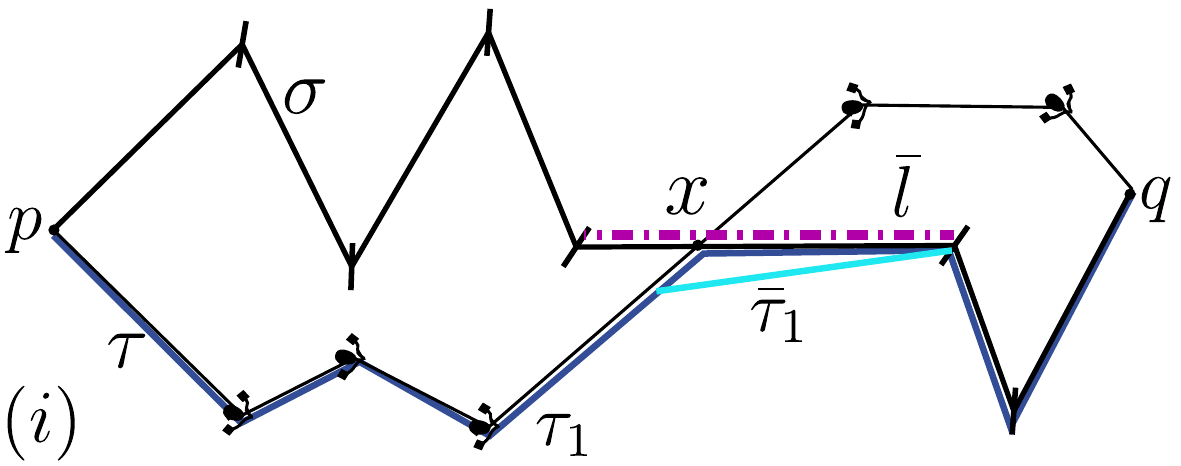}\\
\vspace{-4cm}
\hspace{8cm}  \includegraphics[scale=0.6]{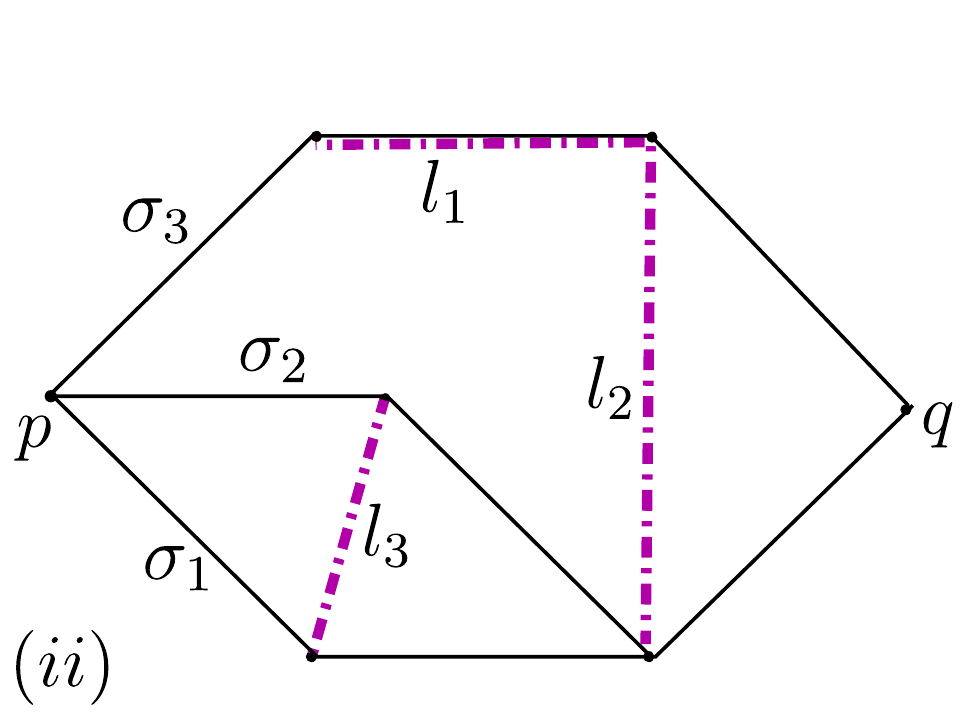} 
\caption{{Concatenating and shortening of the geodesic in 
    Remark \ref{rem1}. In (i), the turning vertices of the base polygonals
 $\sigma$ and $\tau$ are indicated by, respectively, dashes and
 mid-markers. The concatenated shortened polygonal $\tau_1$ is in blue; it can
 be approximated by a polygonal $\bar\tau_1$ with values in
 $\Omega\setminus \bar L$, by means of a segment (in light blue) that
 avoids $\bar l$. In (ii) the displayed configuration of cuts is
 minimal, yet cuts are not separated. There are three geodesic
 polygonals $\{\sigma_i\}_{i=1}^3$ from $p$ to $q$.}}
\label{Fig4}
\end{figure}

\begin{remark} \label{rem1}
When the minimal configuration $\bar L$ consists of disjoint segments
  $\{l_i\}_{i=1}^n$, then each geodesic $\sigma$ from $p$ to $q$ in $\Omega\setminus
  \bar L$ does not contain any cuts. Indeed, assume by contradiction that $\bar l\subset\sigma$, for
some cut $\bar l\subset\bar L$. Denote $\bar{\bar L} = \bar L
\setminus \bar l$, then by the maximality condition in Lemma
\ref{lem2}, there exists a piecewise $\mathcal{C}^1$ curve $\tau$ from $p$ to $q$ in
$\Omega\setminus \bar{\bar L}$ satisfying
$length(\tau)<length(\sigma)$. Hence, $\tau$ must intersect $\bar l$
at only one point which we call $x$. Consider two piecewise
$\mathcal{C}^1$ curves: the curve $\tau_1$ obtained by concatenating the portion
of $\tau$ from $p$ to $x$, with the portion of $\sigma$ from $x$ to
$q$, and the curve $\tau_2$ obtained by concatenating the portion
of $\sigma$ from $p$ to $x$, with the portion of $\tau$ from $x$ to
$q$. One of these curves, say $\tau_1$, must satisfy:
$$length(\tau_1)<length(\sigma).$$ 
But then one can approximate $\tau_1$
by another piecewise $\mathcal{C}^1$ curve $\bar \tau_1:[0,1]\to
\Omega\setminus \bar L$ (see Figure \ref{Fig4} (i)), to the effect that
$length(\bar \tau_1)<length (\sigma)$, which contradicts $\sigma$
being a geodesic.  
Note that in the general case of $\bar L$ supported on the minimal
graph $\bar G$ with vertex degrees possibly exceeding $1$, the above
property is no more true (see Figure \ref{Fig4} (ii)).
\end{remark}

\section{Proof of Theorem \ref{thm2}. Step 2: ordering the geodesics}\label{sec22}

Assume (\ref{S}) and let $p,q$ be two distinct points belonging to one
connected component of $\bar\Omega\setminus L$. In the previous section we showed that, without loss of
generality, the set of cuts $L=\bigcup_{i=1}^n l_i$ satisfies assertions in Corollary \ref{cor3}.
From now on, we work assuming these additional properties and denote
by $(G, V, L)$  a minimal configuration (instead of the notation $(\bar F, \bar V,
\bar L)$ used in section \ref{sec21}).

% of the graph $(G,V, E)$.

\medskip

The (finite, nonempty) set of all geodesics from $p$ to $q$ in
$\Omega\setminus L$ has a partial order relation $\preceq$ in:
\begin{equation}\tag{O}\label{O}
\left[~\mbox{\begin{minipage}{15.3cm}
Given two geodesics from $p$ to $q$ in $\Omega\setminus L$:
$$\sigma_1=\overline{pa_{i_1} a_{i_2}\ldots a_{i_k}q},\qquad \sigma_2=\overline{pa_{j_1} a_{j_2}\ldots a_{j_s}q},$$
we write $\sigma_1\preceq\sigma_2$ provided that the concatenated polygonal: 
$$\sigma=\sigma_1\ast(\sigma_2)^{-1} = \overline{pa_{i_1} a_{i_2}\ldots a_{i_k}qa_{j_s}
  a_{j_{s-1}}\ldots a_{j_1}p}$$ 
is the boundary of (finitely many) open bounded
connected regions in $\R^2$, and moreover $\sigma$ is oriented counterclockwise with
respect to all of these regions.
\end{minipage}}\right.
\end{equation}

\smallskip

\begin{lemma}\label{lem4}
In the above setting, we have:
\begin{itemize}
\item[(i)]  there exist the unique geodesic $\sigma_{min}$ and the
  unique geodesic $\sigma_{max}$ such
  that $\sigma_{min}\preceq\sigma\preceq \sigma_{max}$ for all geodesics $\sigma$ from $p$
  to $q$ in $\Omega\setminus L$; we call $\sigma_{min}$ the least and $\sigma_{max}$ the greatest geodesic,
\item[(ii)] there exists a chain of geodesics
  $\sigma_1\preceq\sigma_2\ldots\preceq\sigma_N$, such that $\sigma_1
  = \sigma_{min}$, $\sigma_N=\sigma_{max}$ and that the
  consecutive geodesics cover each other, i.e. for all $i=1\ldots N-1$ there holds:
$\sigma_i\neq \sigma_{i+1}$ and if $\sigma_i\preceq\sigma\preceq
\sigma_{i+1}$ for some other geodesic $\sigma$, then $\sigma=\sigma_i$ or $\sigma=\sigma_{i+1}$.
\end{itemize}
\end{lemma}
\begin{proof}
{\em Step 1.} For the least geodesic statement in (i), it suffices to
show that if $\sigma_1, \sigma_2$ are two minimal
elements for the partial order $\preceq$, then necessarily
$\sigma_1=\sigma_2$. To this end, we will construct a geodesic
$\sigma$ with $\sigma\preceq\sigma_1$ and $\sigma\preceq\sigma_2$. 
The statement for the greatest geodesic follows by a symmetric argument.
 
\smallskip

We write: $\sigma_1=\overline{pa_{i_1} a_{i_2}\ldots a_{i_k}q}$ and
$\sigma_2=\overline{pa_{j_1} a_{j_2}\ldots a_{j_s}q}$ .
Observe first that $\sigma_1$ and $\sigma_2$ cannot have a common
point $x\not\in\{p,q\}$ that is not a vertex in $V$, unless they have a common edge
$\overline{a_{i_m}a_{i_{m+1}}}=\overline{a_{j_l}a_{j_{l+1}}}$, then
containing $x$. Indeed, in the aforementioned situation, both polygonals: $\tau_1$ obtained by
concatenating the portion of $\sigma_1$ from $p$ to $x$, with the portion of $\sigma_2$ from $x$ to
$q$, and $\tau_2$ concatenating the portion
of $\sigma_2$ from $p$ to $x$, with the portion of $\sigma_1$ from $x$ to
$q$, would satisfy:
$$length(\tau_1) = length(\tau_2)=length(\sigma_1)=length(\sigma_2).$$
Similarly to the construction in the proof of Corollary \ref{cor3},
see also Figure \ref{Fig4}, we could then approximate $\tau_1$
(and also $\tau_2$) by a piecewise $\mathcal{C}^1$ curve
$\tau:[0,1]\to\bar\Omega\setminus L$ with $length(\tau)$ strictly less
than the four coinciding lengths above. This would contradict
$\sigma_i$-s being geodesics.

\smallskip

Let $a_{i_m}=a_{j_l}$ be the first common vertex of $\sigma_1$ and
$\sigma_2$, beyond $p$. If $i_m=i_1$ and $j_l=j_1$, then we include
$\overline{pa_{i_1}}$ as the starting portion of $\sigma$; otherwise we
choose $\overline{pa_{i_1} a_{i_2}\ldots a_{i_m}}$ in case the concatenation
$\overline{pa_{i_1} \ldots  a_{i_m}a_{j_{l-1}}\ldots a_{j_1}p}$ has
the counterclockwise orientation with respect to the bounded open
region it encloses, or $\overline{pa_{j_1} a_{j_2}\ldots a_{j_l}}$ in
the reverse case. Let $a_{i_{\bar m}}=a_{j_{\bar l}}$ be the second common vertex of $\sigma_1$ and
$\sigma_2$, beyond $a_{i_m}$;  we choose the least of $\overline{a_{i_m} \ldots
a_{i_{\bar m}}}$ and $\overline{a_{j_l} \ldots a_{j_{\bar l}}}$, as
above, to be concatenated with the previous portion of $\sigma$. 
By such inductive procedure, we obtain a required geodesic $\sigma$
that satisfies $\sigma\preceq \sigma_1$ and $\sigma\preceq
\sigma_2$. From minimality, it follows that $\sigma=
\sigma_1=\sigma_2$, and so $\sigma=\sigma_{min}$ is the least element
for $\preceq$.

\smallskip

{\em Step 2.} To prove (ii), we set $\sigma_1 =\sigma_{min}$ and
$\sigma_N =\sigma_{max}$ for some $N\geq 2$. If $\sigma_N$
covers $\sigma_1$, then $\sigma_1\preceq\sigma_N$ is the required
chain. Otherwise, there exists a geodesic $\sigma\not\in\{\sigma_1,
\sigma_N\}$ such that $\sigma_1\preceq \sigma\preceq \sigma_N$. If
$\sigma$ covers $\sigma_1$ then we write $\sigma_2=\sigma$, if it is
covered by $\sigma_N$ then we set $\sigma_{N-1}=\sigma$. If none of
the above holds, there must exist a geodesic $\tau\not\in\{\sigma_1,
\sigma, \sigma_N\}$ such that:
$$\sigma_1\preceq \tau\preceq \sigma\quad \mbox{ or }\quad
\sigma\preceq\tau\preceq\sigma_N.$$
We continue in this fashion until the process is stopped, which will
occur in finitely many steps due to the finite number of geodesics
from $p$ to $q$ in $\Omega\setminus L$. 
\end{proof}

\begin{figure}[htbp]
\centering
\includegraphics[scale=0.6]{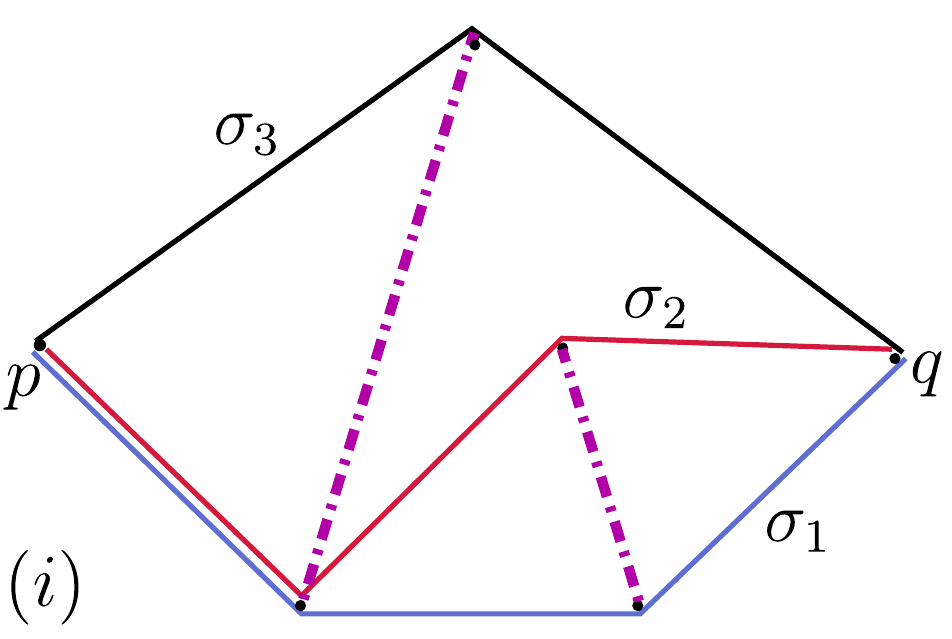}\qquad 
\quad \includegraphics[scale=0.6]{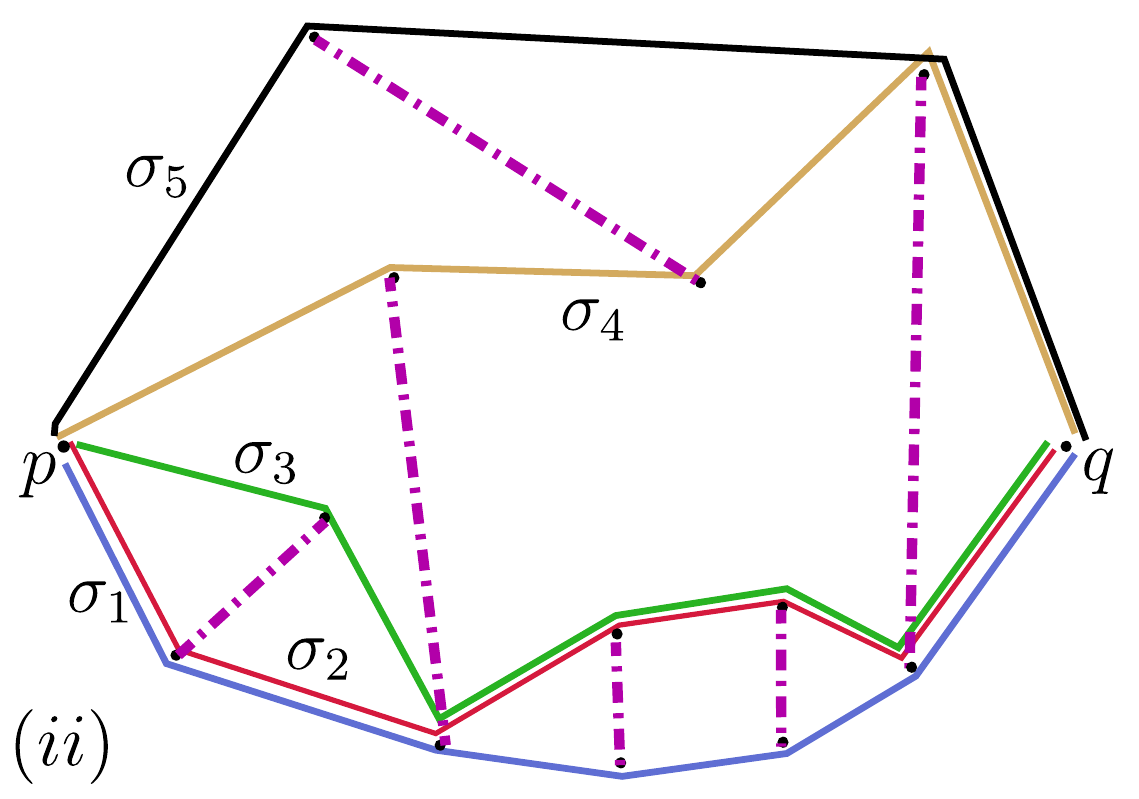}
\caption{{Examples of sequence of geodesics
    produced in Lemma \ref{lem4}: diagram (i) refers to the
    configuration $L, p,q$ in Figure \ref{Fig1} (ii) where the resulting
    sequence consists of the following geodesics: $\sigma_1=\sigma_{min}$ in blue, $\sigma_2$ in red
    and $\sigma_3=\sigma_{max}$ in black; in a more complex diagram (ii) the
    sequence consists of: $\sigma_1=\sigma_{min}$ in blue, $\sigma_2$
    in red, $\sigma_3$ in green, $\sigma_4$ in brown and $\sigma_5=\sigma_{max}$ in black.}}
\label{Fig5}
\end{figure}

\begin{lemma}\label{lem5}
In the above setting, let
$\sigma_1\preceq\sigma_2\ldots\preceq\sigma_N$ be as in Lemma
\ref{lem4} (ii). For each $r=1\ldots N-1$, let $R_r$ be the open, bounded region enclosed by the concatenation
  $\sigma_r\ast (\sigma_{r+1})^{-1}$. We set $R_0=\Omega\setminus
  \bigcup_{r=1}^{N-1} \bar R_r$ to be the exterior region relative to the
  concatenation $\sigma_1\ast (\sigma_{N})^{-1}$. Then, for each tree
  $T$ that is a connected component of $G$, there holds:
\begin{itemize}
\item[(i)] $T$ has nonempty intersection with the interior of exactly one region $\bar R_r$,
\item[(ii)] if $T\subset \bar R_r$ for $r=1\ldots N-1$, then $T$ has
  vertices on both $\sigma_r$ and $\sigma_{r+1}$.
\end{itemize}
Moreover, if $p,q\in \partial\Omega$, then there are no trees in $R_0$.
\end{lemma}
\begin{proof}
{\em Step 1.}  If $T$ violated the condition in (i) then there would
exist a path $\alpha\subset T$ and
two distinct points $A, B\in\alpha$ (which are not necessarily the vertices in
$V$) such that $A\in R_{r-1}$, $B\in R_{r}$ for some
$r=1\ldots N$ (where we set $R_N=R_0$), and such that the portion of
$\alpha$ between $A$ and $B$ crosses the geodesic $\sigma_r$.  This would contradict condition
(ii) in definition (\ref{G}). % for the geodesic $\sigma_r$.

\smallskip

{\em Step 2.} To prove (ii), assume without loss of generality that
all vertices of a maximal tree $T\subset \bar R_r$ belong to
$\sigma_r$. Call $A$ the leaf of $T$ that is closest to $p$ along
$\sigma_r$, and $B$ the leaf that is closest to $q$. By Corollary
\ref{cor3} (ii) and since each tree has at least two leaves, there
must be $A\neq B$. Consider the (unique) path $\alpha\subset T$
connecting $A$ and $B$. If $\alpha\subset \sigma_r$, then there would
be $T=\alpha$. Also, in this case $\alpha$ together with edges of
$\sigma_r$ immediately preceding $A$ and immediately succeeding $B$
would form a straight segment, contradicting the minimality of
$G$. Thus, $\alpha$ passes through $R_r$.

Call $A'$ the first vertex on $\alpha$ whose immediate successor
belongs to $R_{r}$, and $B'$ the last vertex whose immediate
predecessor belongs to $R_{r}$; there may be $A'=A$ or
$B'=B$. Call $\alpha' = \overline{A'a_{i_1}\ldots a_{i_l} B'}\subset
\alpha$ the unique path in $T$ connecting $A'$ with $B'$, and denote
by $D\subset R_r$ the region enclosed by the concatenation of
$\alpha'$ and the portion of $\sigma_r$ between $A'$ and $B'$. 

\smallskip

\begin{figure}[htbp]
\centering
\includegraphics[scale=0.65]{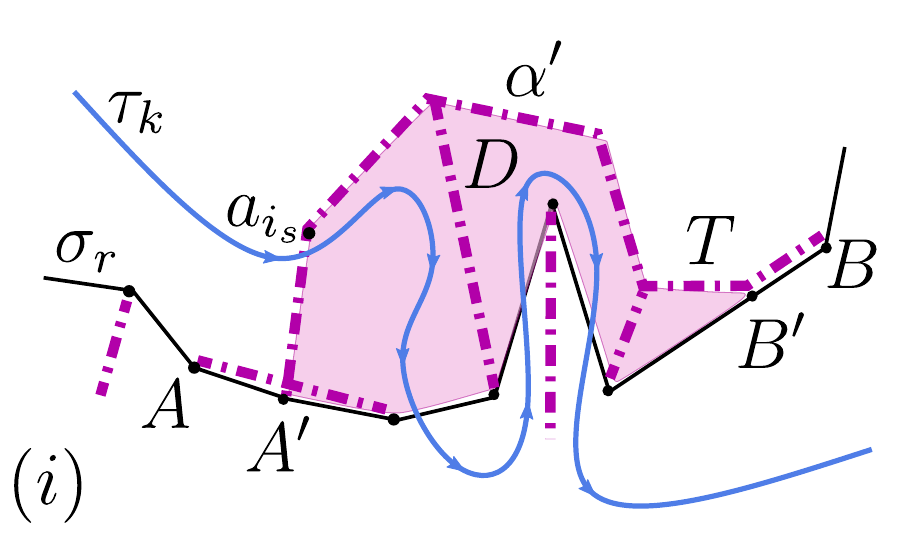}
\qquad \qquad \includegraphics[scale=0.65]{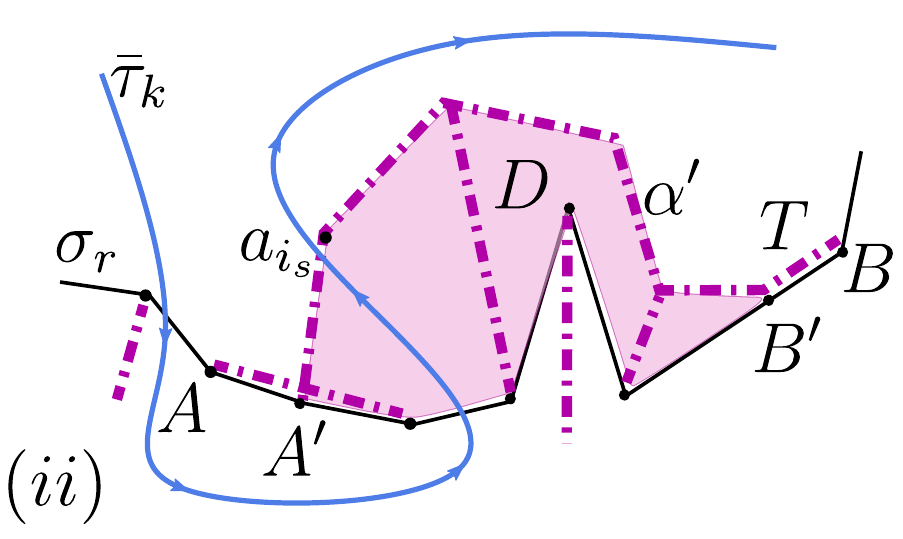}
\caption{{Notation in the proof of Lemma \ref{lem5}, Step 2: in (i)
    the curve $\tau_k$ enters the (shaded) region $D$, hence the indicated
    vertex $a_{i_s}\in\alpha'$ is of type I from left; in (ii), existence of a
    shortening path $\bar\tau_k$ which exits $D$ via the removed edge
    portion preceding $a_{i_s}$ in $\alpha'$, implies that $a_{i_s}$ is of type II from left.}}
\label{Fig5.2}
\end{figure}

We now label each vertex $a_{i_s}\in\alpha'\setminus \sigma_r$ as
{\em type I/II from left}, provided that there exists a
sequence of piecewise $\mathcal{C}^1$ paths
$\{\tau_k:[0,1]\to\bar\Omega\setminus L(1/k)\}_{k=1}^\infty$ with
$\tau_k(0)=p$, $\tau_k(1)=q$ and $length(\tau_k)<\mbox{dist}(p,q)$,
where $L(1/k)$ denotes the modified cut set $L$ in which the edge
  $\overline{a_{i_{s-1}}a_{i_s}}$ in the graph $G$ is replaced by
the shortened segment $\overline{a_{i_{s-1}}a_{i_s}(1/k)}$ with $a_{i_s}(1/k)=
  a_{i_s}-\frac{1}{k}(a_{i_s}-a_{i_{s-1}})$. Further, we request that:
\begin{description}%[leftmargin=10mm]
\item[{$~~$ \normalfont \em Type I from left}] each $\tau_k$ {\em enters} 
$D$, only once, through the removed segment portion
  $\overline{a_{i_s}(1/k) a_{i_s}}$. 
\item[{$~~$ \normalfont \em Type II from left}] each $\tau_k$ {\em exits}
 the region $D$, only once, through $\overline{a_{i_s}(1/k) a_{i_s}}$. 
\end{description}

\smallskip

\noindent Similarly, we label $a_{i_s}\in\alpha'\setminus \sigma_r$ as
{\em type I/II from right}, when there exists a sequence of piecewise $\mathcal{C}^1$ paths
$\{\bar\tau_k:[0,1]\to\bar\Omega\setminus \bar L(1/k)\}_{k=1}^\infty$ with
$\bar\tau_k(0)=p$, $\bar\tau_k(1)=q$, $length(\bar\tau_k)<\mbox{dist}(p,q)$,
and where $\bar L(1/k)$ stands for the modified cut set $L$ in which 
  $\overline{a_{i_{s}}a_{i_{s+1}}}$ is replaced by
$\overline{a_{i_{s}}(1/k) a_{i_{s+1}}}$ with $a_{i_s}(1/k)=
  a_{i_s}+\frac{1}{k}(a_{i_{s+1}}-a_{i_{s}})$. Moreover, we request that:
\begin{description}%[leftmargin=10mm]
\item[{$~~$ \normalfont \em Type I from right}] each $\bar\tau_k$ {\em
    enters} $D$, only once, through the removed segment portion
  $\overline{a_{i_s} a_{i_s}(1/k)}$. 
\item[{$~~$ \normalfont \em Type II from right}] each $\bar\tau_k$ {\em exits}
 $D$, only once, through $\overline{a_{i_s}a_{i_s}(1/k)}$. 
\end{description}
In the definitions above (see diagrams in Figure \ref{Fig5.2}), we set $a_{i_0}=A'$ and
$a_{i_{l+1}}=B'$. By the minimality of $G$, each $a_{i_s}$ must be of
type I or type II from left (it may be both), and it must be of
type I or type II from right (it may be both).

\smallskip

{\em Step 3.} We claim that $a_{i_1}$ has to be of type I from
left. We argue by contradiction and hence assume that $a_{i_1}$
is of type II from left. Note that the length of the portion of the shortening curve
$\bar \tau_k$ between $p$ and the exit point from $D$ is strictly
larger than the distance from $p$ to $a_{i_1}$ in $\Omega\setminus L$,
because all the internal (with respect to $R_r$) angles along $\alpha$ from $A$ to
$A'$ are not greater than $\pi$, whereas the angle at $A'$ is strictly
smaller than $\pi$. Concatenating with the remaining portion of
$\bar\tau_k$ and taking the limit $k\to\infty$, it follows that
there is a geodesic from $p$ to $q$ in $\Omega\setminus L$ passing
through $a_{i_1}$. This contradicts the fact that
$a_{i_1}\not\in\sigma_k$, and proves the claim.

By a similar argument, we can show that if $a_{i_s}$ is of type I from
left, then it is also of type I from right. We argue by contradiction
and hence assume that $a_{i_s}$ is of type II from right. Consider the
curves $\tau_k$ and $\bar\tau_k$ corresponding to the two assumed properties
of $a_{i_s}$; they must intersect at some point $C$ occurring after
$\tau_k$ enters $D$ and before $\bar\tau_k$ exits from $D$. Define the curves:
$\eta_k$ as the concatenation of the portion of $\tau_k$ from $p$ to $C$ with
the portion of $\bar\tau_k$ from $C$ to $q$, and $\bar\eta_k$ as the
concatenation of the portion of $\bar \tau_k$ from $p$ to $C$ with the portion of $\tau_k$
from $C$ to $q$. Since $\bar\eta_k\subset \Omega\setminus L$, it
follows that $length(\bar\eta_k)\geq \mbox{dist}(p,q)$. Consequently:
$$length(\eta_k) = length(\tau_k) + length(\bar\tau_k)-length(\bar\eta_k)<\mbox{dist}(p,q).$$ 
The only possibility for this when taking the limit $k\to\infty$, we obtain the existence of a
geodesic from $p$ to $q$ in $\Omega\setminus L$ passing
through $a_{i_s}$. This contradicts the fact that $a_{i_s}$ is not on any geodesic.

\begin{figure}[htbp]
\centering
\includegraphics[scale=0.65]{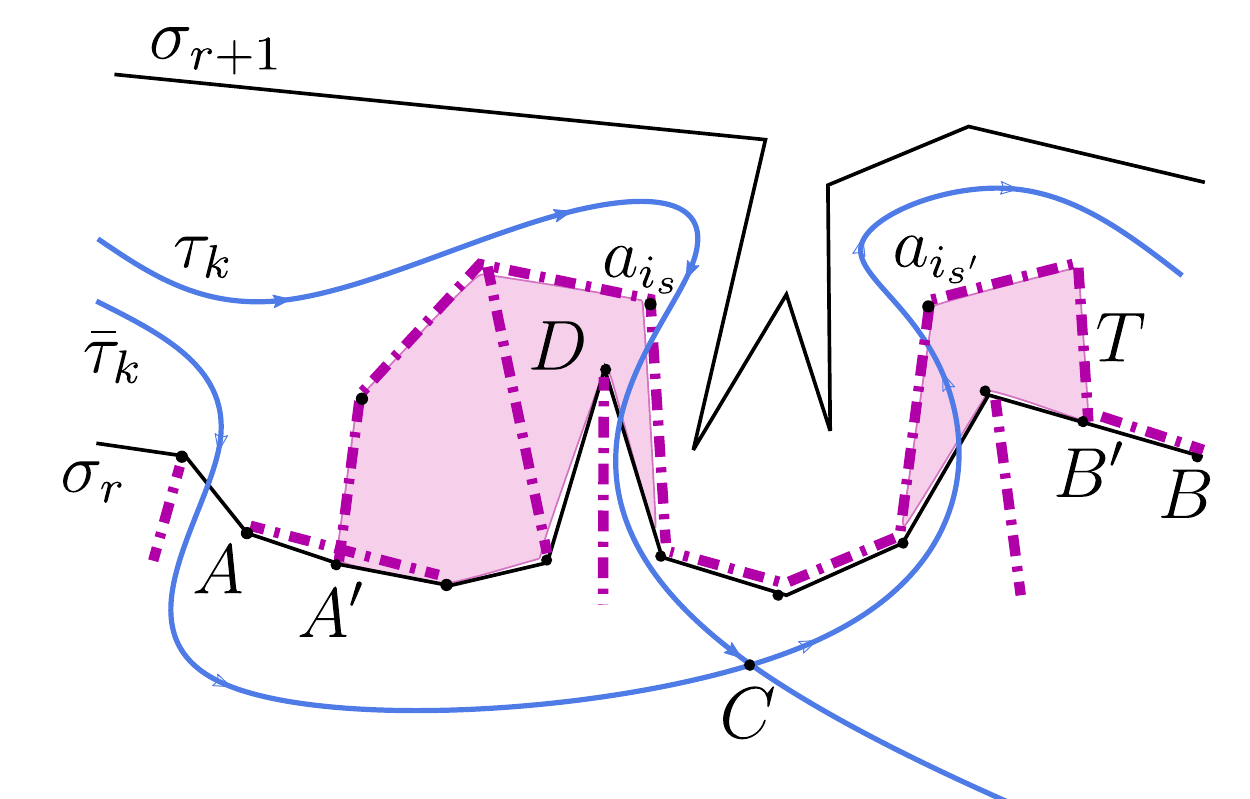}
%\qquad \quad \includegraphics[scale=0.6]{config9.pdf}
\caption{{Concatenating curves $\tau_k$ and $\bar\tau_k$ at the
    intersection $C$ in the proof of Lemma \ref{lem5} Step 3, when the
    region $D$ has multiple connected components.}}
\label{Fig5.5}
\end{figure}

Finally, we argue that if $a_{i_s}$ is of type I from right, then the
next vertex $a_{i_s'}$ on $\alpha'$ that belongs to $R_r$ must be
of type I from left as well. If not, then $a_{i_s'}$ is of type II from
left and we can define the point $C$ and the concatenated curves $\eta_k$
and $\bar\eta_k$ as in the previous reasoning. Again,
$length(\bar\eta_k)\geq \mbox{dist}(p,q)$, so $length(\eta_k)<
\mbox{dist}(p,q)$. However, we may replace the portion of $\eta_k$
between the entry point of $\tau_k$ to $D$, and the exit point of
$\bar\tau_k$ from $D$, by a shorter curve (see Figure \ref{Fig5.5})
which is completely contained in $\Omega\setminus L$. Indeed, when
$i_{s'}=i_{s+1}$, then the said curve may follow the segment
$\overline{a_{i_s}a_{i_{s'}}}\subset \alpha'$. When  $i_{s'}\neq i_{s+1}$, then 
the portion of the polygonal $\alpha'$ between $a_{i_s}$ and
$a_{i_{s'}}$ has all internal angles (with respect to $R_r$) not
greater than $\pi$, so one can simply take the geodesic from $a_{i_s}$ to
$a_{i_{s'}}$ in $\Omega\setminus L$. As a consequence and passing to
the limit with $k\to\infty$, we obtain a 
geodesic from $p$ to $q$ in $\Omega\setminus L$ passing
through $a_{i_s}$ and $a_{i_{s'}}$. This contradicts $a_{i_s}, a_{i_{s'}}$ not being on any geodesic.

\smallskip

{\em Step 4.} Applying the observations from Step 3, it follows that
the vertex $a_{i_l}$ must be of type I from right. However this is
impossible by a symmetric argument to $a_{i_1}$ not being of type II
from left. This ends the proof of (ii). 
In case $p,q\in\partial\Omega$, the region $R_0$ consists of two
connected components, and hence any tree $T\in R_0$ would have
vertices either only on $\sigma_1$ or only on $\sigma_N$. By the same arguments as
above, this is impossible, which implies the final statement of the lemma.
\end{proof}

%\section{Proof of Theorem \ref{thm2}. Step 3: ordering regions between
%consecutive geodesics}\label{sec22.5}

\begin{figure}[htbp]
\centering
\includegraphics[scale=0.6]{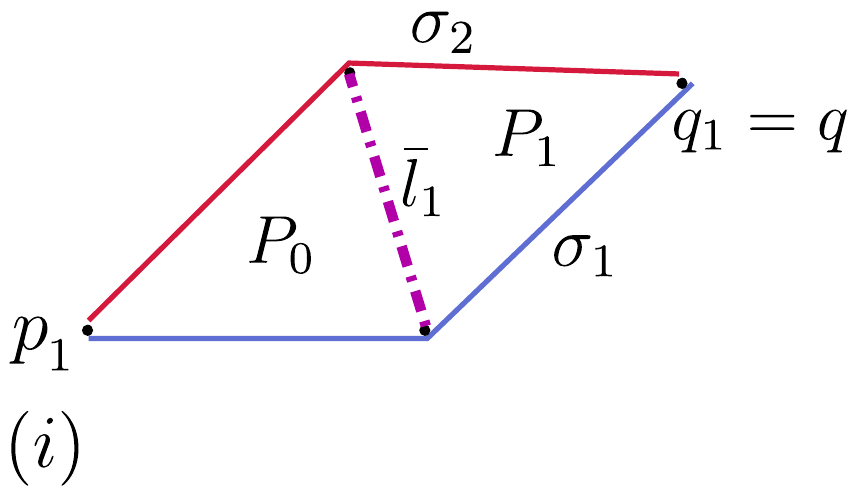}\qquad 
\quad \includegraphics[scale=0.6]{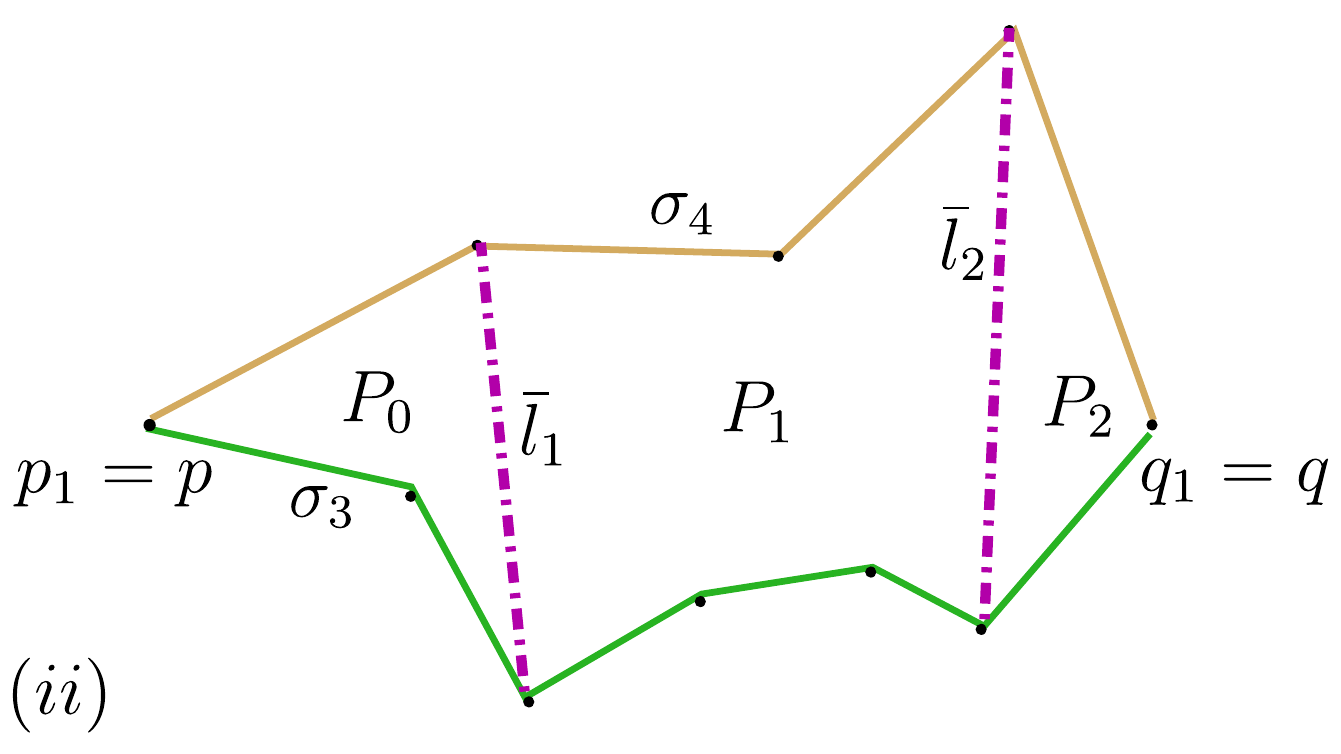}
\caption{{Partitions $\{P_m\}_{m=0}^s$ defined in Lemma \ref{lem5.5}: (i) depicts
    partition of the region $R_1$ corresponding to Figure \ref{Fig5} (i), while (ii)
    depicts the region $R_3$ in Figure \ref{Fig5} (ii). In both figures
    the trees $T_m$ coincide with paths $\alpha_{m-1}^{right}=\alpha_m^{left}$ that are single cuts, and
    consequently all intermediate polygonal collections $Q_m$ are empty.}}
\label{Fig6}
\end{figure}

\begin{figure}[htbp]
\centering
\includegraphics[scale=0.65]{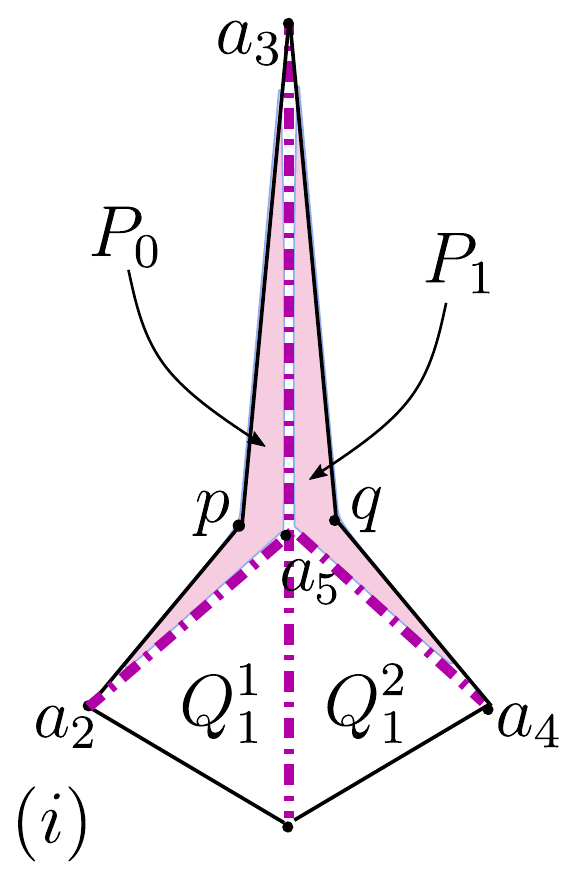}\qquad 
\qquad \includegraphics[scale=0.65]{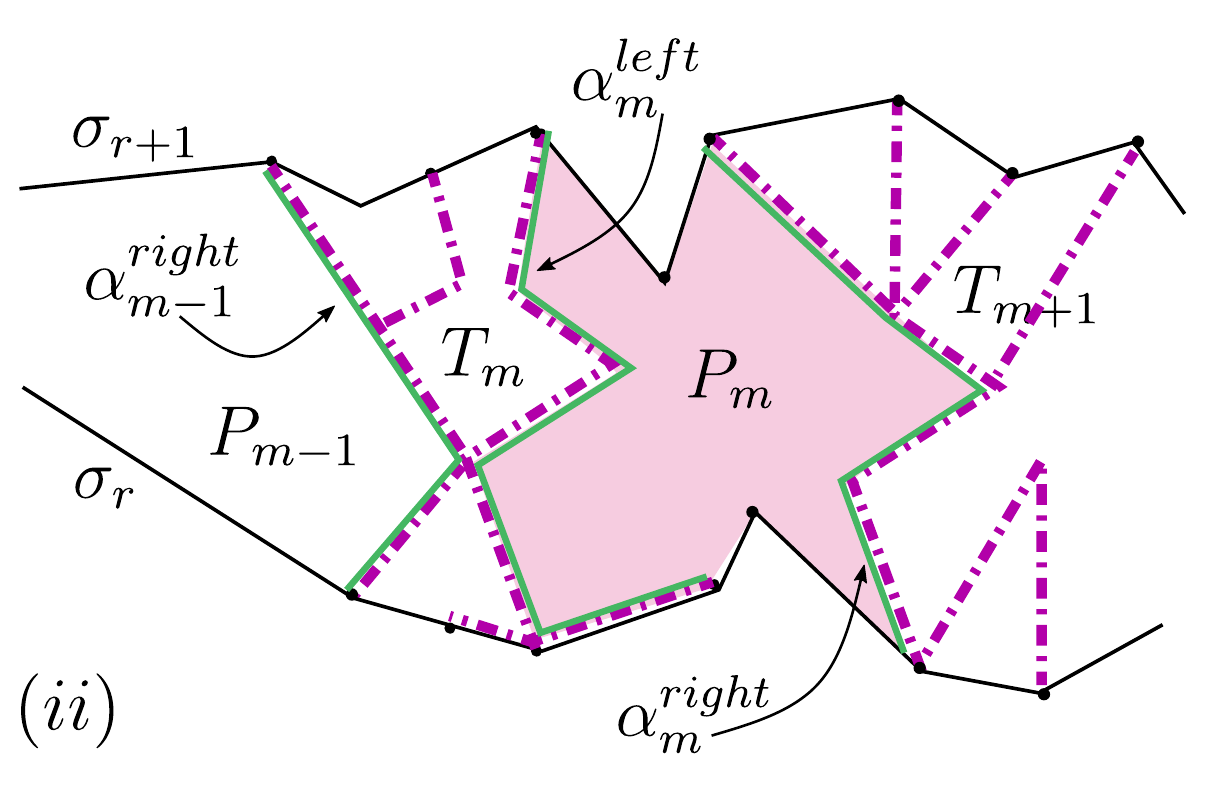}
\caption{{Polygons $\{P_m\}_{m=0}^s$ and polygon families $\{Q_m\}_{m=1}^s$ defined in
    Lemma \ref{lem5.5}: (i) corresponds to the unique region $R_1$ in
    Figure \ref{Fig1.5} (iii) with paths: $\alpha_0^{left}=\overline{pp}$, $\alpha_0^{right}=\overline{a_2a_5a_3}$,
    $\alpha_1^{left}=\overline{a_4a_5a_3}$, $\alpha_1^{right}=\overline{qq}$; (ii) is a general diagram
    depicting the partition of the region $R_r$.}}
\label{Fig6.5}
\end{figure}

We close the above discussion by pointing out that in case
$p,q\not\in \partial\Omega$, there may be a tree (or even multiple
trees) in $R_0$, with vertices both 
on $\sigma_1$ and $\sigma_N$ (see Figure \ref{Fig9.1} in section
\ref{sec24}).  The next main result of this section allows for the
lexicographic ordering of the connected components of
$\bar\Omega\setminus L$. Namely, we have (see example in Figures
\ref{Fig6} and \ref{Fig6.5}):

\begin{lemma}\label{lem5.5}
In the above setting, let $\sigma_1\preceq\sigma_2\ldots\preceq\sigma_N$ be as in Lemma
\ref{lem4} (ii). Fix $r=1\ldots N-1$ and consider the region $R_r$
enclosed between two consecutive geodesics $\sigma_r = \overline{pa_{i_1} a_{i_2}\ldots a_{i_k}q}$ and
$\sigma_{r+1}=\overline{pa_{j_1} a_{j_2}\ldots a_{j_l}q}$, as in Lemma
\ref{lem5}. Consider further the set of maximal trees $\{T_m\}_{m=1}^s$ which
are the connected components of $G$ contained in $\bar R_r$. Then we have:
% which covers $\sigma_r$.  
%Then there exists cuts $\{\bar l_1\ldots \bar
%l_s\}\subset L$ such that there holds:
\begin{itemize}
\item[(i)]  the ordering $T_1,\ldots, T_s$ can be made so that each
  leaf of $T_i$ on $\sigma_r$ (respectively $\sigma_{r+1}$) precedes
  each leaf of $T_j$ on $\sigma_r$  (resp. $\sigma_{r+1}$),
  when $i<j$.
\end{itemize}
The region $R_r\setminus L$ is the union of $s+1$ (open) polygons $\{P_m\}_{m=0}^s$ and
of additional families of polygons $\{Q_m\}_{m=1}^s$, described as follows:
\begin{itemize}
\item[(ii)] we denote $\alpha_0^{left}=\overline{p_1p_1}$ and
  $\alpha_s^{right}=\overline{q_1q_1}$, where $p_1, q_1$ are two common vertices of
  $\sigma_r$ and $\sigma_{r+1}$, such that $\overline{pa_{i_1}\ldots p_1}=
\overline{pa_{j_1}\ldots p_1}$ and $\overline{q_1\ldots a_{i_k}q}=
\overline{q_1\ldots a_{j_l}q}$ (we take the last, along
$\sigma_r$, vertex with the said property to be $p_1$ and the first
vertex to be $q_1$). For each $m=1\ldots s$ we denote
$\alpha_{m-1}^{right}$ (respectively, $\alpha_m^{left}$) the unique
path in $T_m$ joining its first (resp. its last) vertex
on $\sigma_r$ with its first (resp. the last) vertex on $\sigma_{r+1}$, both counting
from $p_1$. Note that there may be $\alpha_m^{left}= \alpha_{m-1}^{right}.$ 
Then, the boundary of each $P_m$ consists of paths $\alpha_m^{left}$,
$\alpha_m^{right}$ and of the intermediate portions of $\sigma_r$ and
$\sigma_{r+1}$ which are concave with respect to $P_m$. Namely,
all interior angles of  $P_m$ which are not on $\alpha_m^{left}\cup
\alpha_m^{right}$ are not less than $\pi$. %
Finally, there are no cuts in $P_m$.
\item[(iii)] each family $Q_m$ consists of finitely many polygons
  $Q_m^f$, that are the connected components of $R_r\setminus L$
  enclosed between $\alpha_m^{left}$, $\alpha_{m-1}^{right}$ and the
 portions of $\sigma_r$ and $\sigma_{r+1}$. The
  boundary of each $Q_m^f$ consists of a single path within $T_m$ plus
  a single portion of the geodesic $\sigma_r$ or $\sigma_{r+1}$. It has
  all interior angles not at vertices belonging of $T_m$ concave with
  respect to $R_m$.
\end{itemize}
\end{lemma}
\begin{proof}
For (i), consider first $\sigma_r$ and recall that each vertex
$a_{i_1},\ldots, a_{i_k}$ is an endpoint of some cut, which belongs to
some maximal tree $T\subset G$. If $T$ extends inside the region $R_r$, then 
it must have vertices on both $\sigma_r$ and $\sigma_{r+1}$, by Lemma \ref{lem5}.
The same reasoning can be applied to cuts emanating from $\sigma_{r+1}$.
We can now order the trees $\{T_m\}_{m=1}^s$, based on how many vertices (along $\sigma_r$ and
$\sigma_{r+1}$) separate their leaves from $p$. This ordering is
well defined, as trees are non-intersecting. Assertions (ii) and (iii)
follow directly by construction and since $\sigma_r$, $\sigma_{r+1}$ are geodesics. 
\end{proof}

Concluding, we see that the assumption (\ref{S}) may be
replaced by the following modified setup:

\begin{equation}\tag{S1}\label{S1}
\left[~\mbox{\begin{minipage}{15.2cm}
The set of cuts $L$ satisfies assertions of Corollary \ref{cor3}.
The chain of geodesics $\{\sigma_r\}_{r=1}^N$, each from $p$ to $q$ in
$\Omega\setminus L$, satisfies condition in Lemma \ref{lem4} (ii)
with respect to the partial order in (\ref{O}). In agreement with Lemma
\ref{lem5.5}, the set $\Omega\setminus L$ is partitioned into $N$
regions $\{R_r\}_{r=0}^N$:
\begin{itemize}
\item[(i)]  for each $r=1\ldots N-1$, the \lq\lq interior" bounded region
  $R_r$ which is enclosed by $\sigma_r\ast (\sigma_{r+1})^{-1}$, and
  partitioned into polygonal sub-regions
  $\{P_m\}_{m=0}^s\cup\{Q_m\}_{m=1}^s$ corresponding to the
  consecutive trees $\{T_m\}_{m=1}^s$ (we suppress the dependence on
  $r$ in this notation), as specified in Lemma \ref{lem5.5},
\item[(ii)] the \lq\lq exterior" region $R_0 =
  \Omega\setminus\bigcup_{r=1}^{N-1}\bar R_r$.
\end{itemize}
We also define the segment $I = \overline{0,length(\sigma_1)e_1}\subset\R^3$.
\end{minipage}}\right.
\end{equation}

\section{Proof of Theorem \ref{thm2}, a simplified case. Step 3: isometric immersion on 
  interior regions between consecutive cuts}\label{sec23}

Assume (\ref{S}) and let $p,q$ be two distinct points in
$\bar \Omega\setminus L$. In view of the results in previous
sections, the goal is to construct an isometry $u$ as in Theorem \ref{thm2},
separately on each $R_r$ identified in (\ref{S1}).
We first concentrate on the interior case $r=1\ldots N-1$, while in section \ref{sec24} we address
the case $r=0$. 

In this section we treat a simplified scenario in which
all trees $T_1,\ldots,T_s$ consist of single cuts; note that this
occurs, in particular, if all cuts in the original graph $G$  are
non-intersecting: % (as in Theorem \ref{thm3}):

\begin{lemma}\label{lem6}
Assume (\ref{S}) and (\ref{S1}). Fix $r=1\ldots N-1$ and further
assume that:
$$T_m = \bar l_m = \overline{a_{m_1}a_{m_2}} \quad \mbox{ for all } \;
m=1\ldots s.$$
Then there exists a continuous, piecewise affine isometric immersion
$u: R_r\setminus \bigcup_{m-1}^s\bar l_m\to \R^3$, with:
$$u(p) = 0, \quad u(q)=length(\sigma_1)e_1, \quad u(\sigma_r)=u(\sigma_{r+1}) = I.$$
\end{lemma}
\begin{proof}
We will inductively find the matching isometric
immersions $u$ of the consecutive polygons $\{P_m\}_{m=0}^s$. Note
that polygons in $\bigcup_{m=1}^sQ_m$ are absent in the presently discussed case.

\smallskip

{\em Step 1.} On $P_0$, we first fold its
``top'' part so that the image of the portion of $\sigma_{r+1}$ from $p_1$
to the endpoint $B_1$ of the cut $\bar l_1=\overline{A_1 B_1}$
coincides with the sub-interval: 
$$\overline{length(\overline{p 
    a_{i_1}\ldots p_1})e_1, length(\overline{p a_{i_1}\ldots  B_1})e_1 }\subset I.$$
This can be achieved because all internal angles of $\sigma_{r+1}$
at vertices between (but not including) $p_1$ and $B_1$ are at least
$\pi$. A symmetric fold construction can be performed on the
``bottom'' part of $P_0$, along the boundary portion contained in $\sigma_r$. 

\smallskip

As a result, the vector $u(B_1) -u(A_1)$ equals $ \big(length(\overline{p 
    a_{i_1}\ldots B_1}) - length(\overline{p a_{j_1}\ldots  A_1})\big)e_1$
and we consecutively have to find an isometric immersion of the
polygon $P_1$ with the property that writing $\bar
l_2=\overline{A_2B_2}$ with $A_2\in\sigma_r, B_2\in\sigma_{r+1}$, the
length of the vector $\big(u(B_2) - u(A_2)\big)- \big(u(B_1)-u(A_1)\big)$ is
prescribed, and that the images of 
portions of: geodesic $\sigma_r$ between $A_1$ and $A_2$, and of
geodesic $\sigma_{r+1}$ between $B_1$ and $B_2$, are contained in $\R e_1$. 

\smallskip

{\em Step 2.} Assume that $u$ has been constructed on $P_1\cup\ldots
P_{m-1}$, for some $m\leq s-1$. Consider the polygon $P_m$ and the
two related closed convex sets $S_A$ and $S_B$. The set $S_B$ is
defined by specifying its boundary to consist of: the
portion of $\sigma_{r+1}$ between the endpoints  $B_m$ and $B_{m+1}$
of the cuts $\bar l_m$, $\bar l_{m+1}$, respectively, and of the segment
$\overline{B_m B_{m+1}}$. The boundary of the set $S_A$ is: the
portion of $\sigma_r$  between the remaining endpoints  $A_m$ and $A_{m+1}$
of the cuts $\bar l_m$, $\bar l_{m+1}$, and of the segment
$\overline{A_m A_{m+1}}$. We note that the interior of the defined
sets may be empty; for example if $\overline{A_m A_{m+1}}\subset \sigma_r$
then $S_A = \overline{A_m A_{m+1}}\subset \sigma_r$.

\begin{figure}[htbp]
\centering
\includegraphics[scale=0.6]{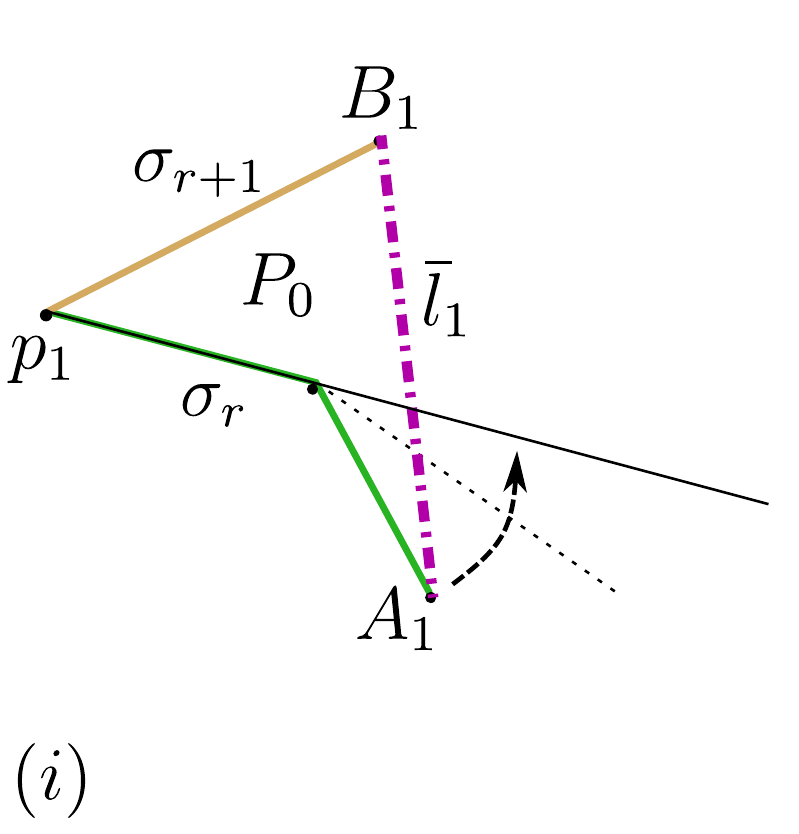}\qquad 
\quad \includegraphics[scale=0.6]{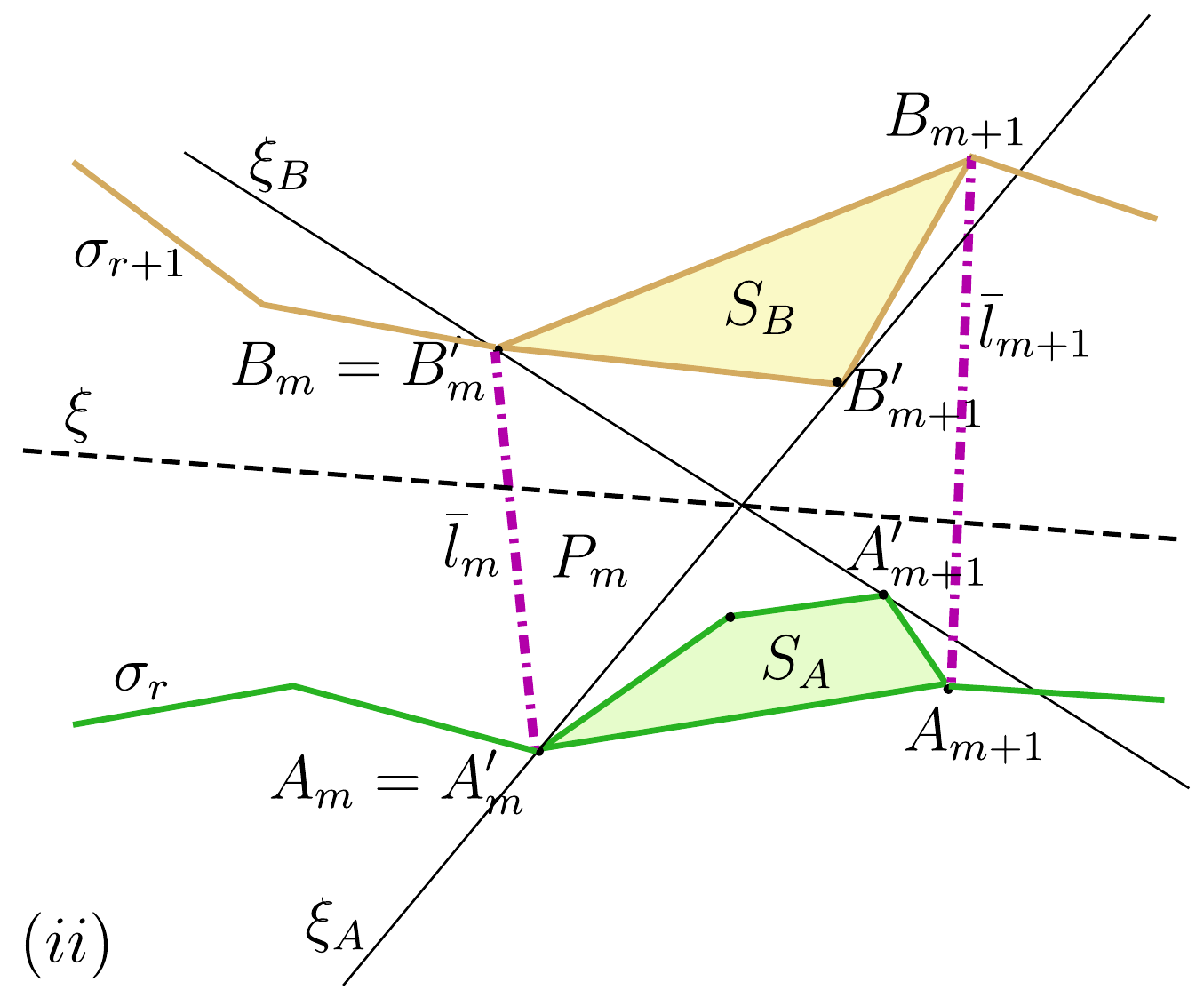}
\caption{{The folding patterns in the proof of Lemma \ref{lem6}: (i)
    corresponds to Step 1 and the interior polygon $P_0$ in Figure
    \ref{Fig6} (ii), the arrow indicates the direction of folding; (ii) corresponds to Step 2 and
    the polygon $P_1$ in Figure \ref{Fig6} (ii).}}
\label{Fig7}
\end{figure}

Since $S_A$ and $S_B$ are closed, convex and disjoint, there exist
precisely two lines $\xi_A, \xi_B$ which are supporting to both sets. Each of these
lines intersects $S_A$ and $S_B$ either at a single vertex (which may
be $A_m$ or $A_{m+1}$ for $S_A$ and $B_m$ or $B_{m+1}$ for $S_B$) or
along the whole segment which is the union of some consecutive edges
in $\sigma_r, \sigma_{r+1}$. Denote
$B_m'\in(\xi_A\cup\xi_B)\cap\sigma_{r+1}$ the vertex which is closest to $B_m$
along $\sigma_{r+1}$, and let
$B_{m+1}'\in(\xi_A\cup\xi_B)\cap\sigma_{r+1}$ be the vertex that is closest to
$B_{m+1}$ (along $\sigma_{r+1}$). Similarly, define $A_{m}', A_{m+1}'\in(\xi_A\cup\xi_B)\cap\sigma_{r}$ as vertices on
$\sigma_r$ which are closest (possibly equal) to, respectively, $A_m$ and $A_{m+1}$
along $\sigma_r$. 

Observe that $\xi_A$ is precisely the line through $A_m'$ and
$B_{m+1}'$ and that it intersects the closures of the cuts $\bar l_m$ and $\bar l_{m+1}$. By
consecutive folding as in Step 1, one constructs an isometric
immersion $v$ of $P_m$ with the property that its both boundary
polygonal sides: from $B_m$ to $B_{m+1}$ and from $A_m$ to $A_{m+1}$
are mapped onto $\xi_A$. By a further rotation, we may ensure that
$\xi_A=\R e_1$. Then:
$$v(B_m) - v(A_m) = \Big(-length(\overline{B_m\ldots B'_{m+1}}) +
length(\overline{A_m' B'_{m+1}}) + length(\overline{A_m\ldots A'_{m}})
\Big)e_1 \doteq \alpha_v e_1.$$

Similarly, by folding on the line $\xi_B$ through $B_m'$ and
$A_{m+1}'$, one obtains an isometric immersion $w$ of $P_m$ with both
polygonal sides (namely, the sides distinct from $\bar l_m$ and $\bar
l_{m+1}$) mapped on $\xi_B$. By a further rotation, we ensure that
$\xi_B=\R e_1$, so that there holds:
$$w(B_m) - w(A_m) = \Big(-length(\overline{B_m\ldots B'_{m}}) -
length(\overline{A_{m+1}' B'_{m}}) + length(\overline{A_m\ldots A'_{m+1}})
\Big)e_1 \doteq \alpha_w e_1.$$

\smallskip

{\em Step 3.} We now estimate the length of the vector $u(B_m)-u(A_m)$
that we need to achieve, and that is determined through previous steps
in the construction. There clearly holds:
\begin{equation}\label{gr}
u(B_m) - u(A_m) = \Big(length(\overline{pa_{j_1}\ldots B_{m}}) -
length(\overline{pa_{i_1\ldots }A_m}) \Big) e_1. 
\end{equation}
Since $P_m$ contains no cuts in its interior, the polygonal:
$\overline{pa_{i_1}\ldots A_{m}A_{m}'B_{m+1}'B_{m+1}\ldots a_{j_s}q}$
(this polygonal follows the portion of $\sigma_r$ up to $A_m'$, then
switches to $\sigma_{r+1}$ along the segment
$\overline{A_m'B_{m+1}'}\subset P_m$, and continues to $q$ along
$\sigma_{r+1}$) cannot be shorter than
$length(\sigma_{r+1})$. Equivalently, we obtain:
\begin{equation*}
\begin{split}
& length(\overline{pa_{i_1}\ldots A_{m}}) +
length(\overline{A_m \ldots A_{m}'}) + length(\overline{A_m'  B_{m+1}'}) 
\\ & \qquad \geq length(\overline{pa_{j_1}\ldots B_{m+1}'}) = 
length(\overline{pa_{j_1}\ldots B_{m}}) + length(\overline{B_m\ldots B'_{m+1}}).
\end{split}
\end{equation*}
By (\ref{gr}) the above yields:
$$\langle u(B_m) - u(A_m), e_1\rangle \leq \alpha_v.$$
By a parallel argument, in which we concatenate $\sigma_{r+1}$ up to $B_m'$ with
the segment $\overline{B_m'A_{m+1}'}$ and then with the portion of
$\sigma_r$ from $A_{m+1}'$ to $q$, there follows the bound:
\begin{equation*}
\begin{split}
& length(\overline{pa_{j_1}\ldots B_{m}}) +
length(\overline{B_m \ldots B_{m}'}) + length(\overline{B_m'  A_{m+1}'}) 
\\ & \qquad \geq length(\overline{pa_{i_1}\ldots A_{m}}) + length(\overline{A_m\ldots a'_{m+1}}),
\end{split}
\end{equation*}
so (\ref{gr}) results in:
$$\langle u(B_m) - u(A_m), e_1\rangle \geq \alpha_w.$$

\smallskip

{\em Step 4.} Let now $\xi$ be any line passing through the
intersection point $\xi_A\cap\xi_B$ and disjoint from the interiors of
$S_A$ and $S_B$ (see Figure \ref{Fig7} (ii)). There exist exactly one
line $\xi(A)$ which is supporting to the convex set $S_A$ and parallel
to $\xi$, and exactly one line $\xi(B)$ supporting to $S_B$ and
parallel to $\xi$. As before, we may fold the top portion of $P_m$ so
that the image of $\overline{B_m\ldots B_{m+1}}$ is a segment within
$\xi(B)$, and fold the bottom portion of $P_m$ so that  the image of
$\overline{A_m\ldots A_{m+1}}$ is a segment in $\xi(A)$. We now
perform two more folds, which map both $\xi(A)$, $\xi(B)$  onto $\xi$,
plus a rigid rotation that maps $\xi$ onto $\R e_1$. Call the
resulting isometric immersion $u_\xi$ and observe that:
$$\mbox{the function }\; \xi\mapsto u_\xi(B_m) - u_\xi(A_m) \;\mbox{
  is continuous}.$$
Since $u_{\xi_A} = v$ and $u_{\xi_B}=w$, the intermediate value
theorem implies that $\langle u_\xi(B_m) - u_\xi(A_m), e_1\rangle$
achieves an arbitrary value within the interval:
$$[\alpha_v,\alpha_w]= \Big[\langle u_{\xi_A}(B_m) - u_{\xi_A}(A_m),
e_1\rangle, \langle u_{\xi_B}(B_m) - u_{\xi_B}(A_m), e_1\rangle\Big].$$
In conclusion, there exists a line $\xi$ such that the corresponding
$u_\xi$ on $P_m$ gives:
$$u_\xi(B_m) - u_\xi(A_m) = u(B_m) - u(A_m).$$ 
We set $u_{\mid P_m}\doteq u_\xi$.

\smallskip

{\em Step 5.} The final step is to construct $u$ on $P_s$. This can be
done by the same folding technique described in Step 1 for $P_0$. We
then note that $\langle u(B_s) - u(A_s), e_1\rangle$ automatically
equals: $length(\overline{A_s\ldots a_{j_l}q_1}) -
length(\overline{B_s\ldots a_{i_k}q_1})$, because $length(\sigma_r) =
length(\sigma_{r+1})$. The proof is done.
\end{proof}

\section{Proof of Theorem \ref{thm2}, the general case. Step 3: isometric immersion on 
  interior regions between and within consecutive trees}\label{sec23.2}

In this section we exhibit a procedure of constructing
an isometric immersion on $R_r$, in the general setting
(\ref{S1}). Namely, we prove the following version of Lemma \ref{lem6}:

\begin{lemma}\label{lem6.2}
Assume (\ref{S}), (\ref{S1}) and fix $r=1\ldots N-1$. 
Then, there exists a continuous, piecewise affine isometric immersion $u$ of $R_r\setminus
\bigcup_{m=1}^s T_m$ into $\R^3$, which satisfies:
$$u(p) = 0, \quad u(q)=length(\sigma_1)e_1, \quad u(\sigma_r)=u(\sigma_{r+1}) = I.$$
\end{lemma}
\begin{proof}
We will inductively find the matching isometric
immersions (always denoted by $u$) of the consecutive polygons in
$P_0, \{Q_m\cup P_m\}_{m=1}^s$. Recall that we have defined two
families of cut paths within each tree $T_m$: the path
$\alpha_{m-1}^{right}$ joining vertices $A_{m-1}^{right}\in\sigma_r$
with $B_{m-1}^{right}\in \sigma_{r+1}$, and the path
$\alpha_{m}^{left}$ joining vertices $A_{m}^{left}\in\sigma_r$
with $B_{m}^{right}\in \sigma_{r+1}$.

\smallskip

{\em Step 1.}  Thus, the polygon $P_0$ is bounded by the the
concatenation of: the portion of $\sigma_{r+1}$ from $p_1$ to
$B_0^{right}$, with $\alpha_0^{right}$, with the portion of $\sigma_r$
from $A_0^{right}$ to $p_1$. We first fold the indicated portion of
$\sigma_{r+1}$ so that it coincides with the sub-interval:
$$\overline{length(\overline{p 
    a_{i_1}\ldots p_1})e_1, length(\overline{p a_{i_1}\ldots  B_0^{right}})e_1 }\subset I.$$
This can be achieved as all internal angles of $\sigma_{r+1}$
at vertices between $p_1$ and $B_0^{right}$ are at least
$\pi$. A symmetric folding can be done along  the portion of
$\sigma_r$ within the boundary of $P_0$, see Figure \ref{Fig7.2} (i).
This construction is similar to Step 1 in the proof of Lemma \ref{lem6} (i).

\smallskip

As a result, the vector $u(B_0^{right}) -u(A_0^{right})$ equals $ \big(length(\overline{p 
    a_{i_1}\ldots B_0^{right}}) - length(\overline{p a_{j_1}\ldots  A_0^{right}})\big)e_1$
and we consecutively have to find an isometric immersion of each
region in the family $Q_1$, with the property that the vector $u(B_0^{right}) - u(A_0^{right})$ is
prescribed, and that the images of the
portion of $\sigma_r$ between $A_0^{right}$ and $A_1^{left}$, and of
$\sigma_{r+1}$ between $B_0^{right}$ and $B_1^{left}$, are contained in $\R e_1$. 

\begin{figure}[htbp]
\centering
\includegraphics[scale=0.6]{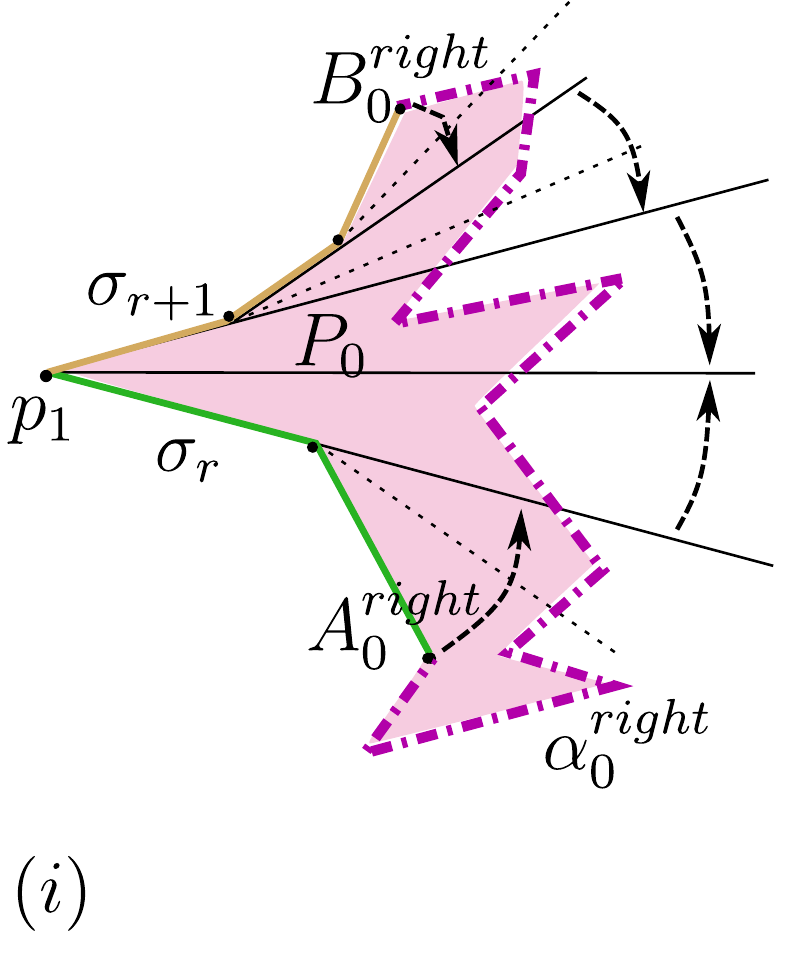}\qquad 
\quad \includegraphics[scale=0.6]{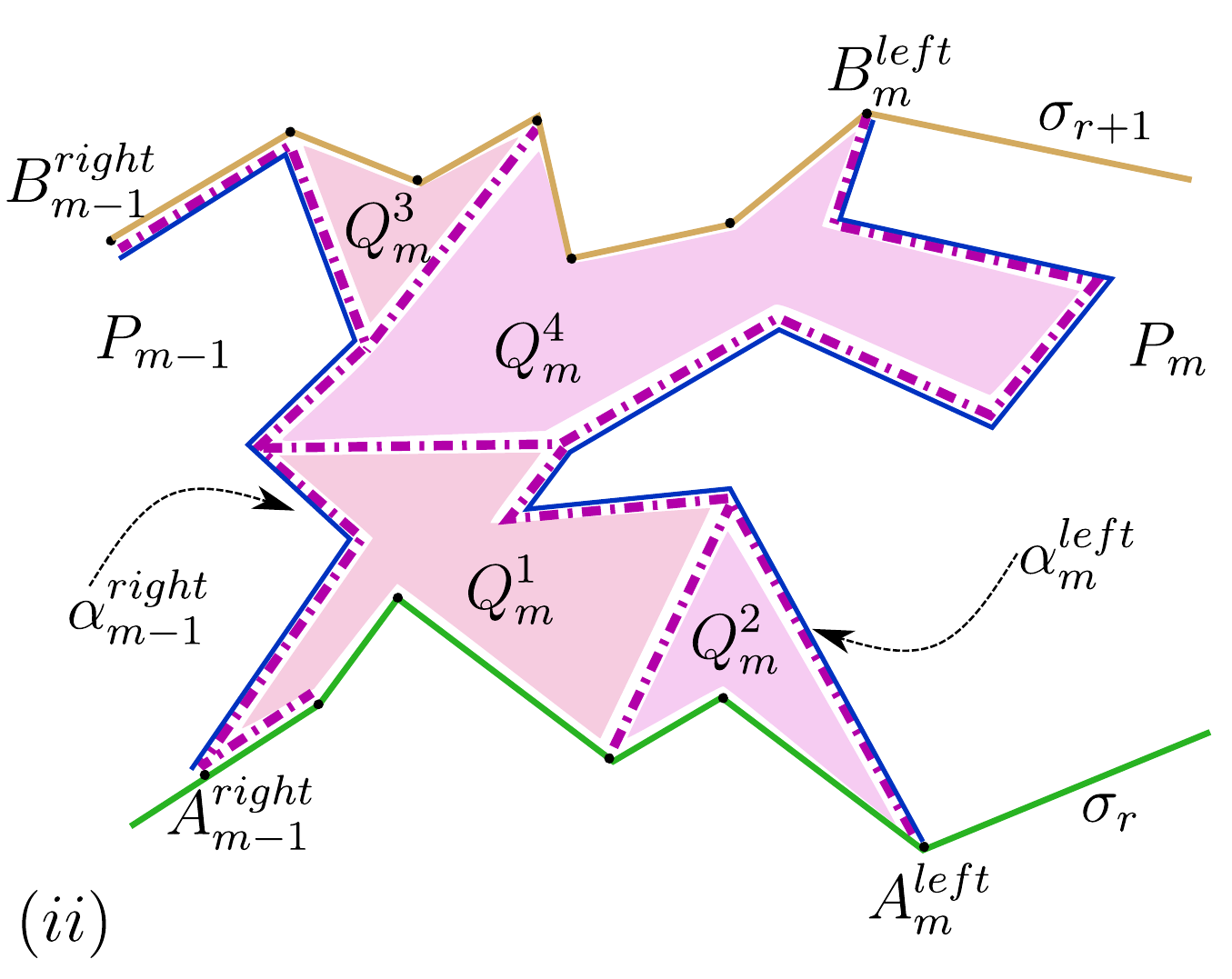}
\caption{{The folding patterns in the proof of Lemma \ref{lem6.2}: (i)
    corresponds to Step 1 and the polygon $P_0$, the arrows indicate
    the directions of folding; (ii) corresponds to Step 2 and the
    collection of polygons $Q_m$ within the tree $T_m$.}}
\label{Fig7.2}
\end{figure}

\smallskip

{\em Step 2.} Assume that $u$ has been constructed on $P_0\cup
Q_1\cup\ldots P_{m-1}$ for some $m\leq s$. Consider the tree $T_m$ and
the corresponding family of polygons $Q_m$ enclosed between the paths
of cuts $\alpha_{m-1}^{right}$, $\alpha_m^{left}$ (both contained in
$T_m$) and the portions of $\sigma_r$ (respectively $\sigma_{r+1}$)
between the vertices $A_{m-1}^{right}$ and $A_m^{left}$ 
(resp. $B_{m-1}^{right}$ and $B_m^{left}$), see Figure \ref {Fig7.2} (ii).
By Lemma \ref{lem5.5} (iii), each polygon $Q_m^f\subset Q_m$ has an
isometric immersion $u$ in which the image of its boundary portion
included in $\sigma_r\setminus T_m$ (or in $\sigma_{r+1}\setminus
T_m$) is  a segment on $\mathbb{R}e_1$. This construction consists of
a collection of simple foldings as in Step 1 and Figure \ref{Fig7.2} (i), that can
be implemented because all the internal (with respect to $R_r$)
angles of $Q_m^f$ at the vertices $\sigma_r\setminus T_m$ and $\sigma_{r+1}\setminus
T_m$, are at least $\pi$. 

\smallskip

{\em Step 3.} Assume that $u$ has been constructed on $P_0\cup Q_1\cup \ldots
P_{m-1}\cup Q_m$ for some $m\leq s-1$. We now aim at describing $u$ on
the polygon $P_m$; note that the vector $u(B_m^{left})-u(A_m^{left})$
is a prescribed, by the previous steps of the proof, scalar
multiple of $e_1$. The construction below is based on the ideas of
Steps 2-4 in the proof of Lemma \ref{lem6}, however the present
setting of trees $T$ replacing the single cuts $\bar l$ requires
taking care of the additional details below.

\smallskip

Call $\xi_A$ (respectively $\xi_B$) the shortest path in $P_m$ that
joins $A_m^{left}$ with $B_m^{right}$ (resp. $A_m^{right}$ with
$B_m^{left}$). We now estimate the length of the vector
$u(B_m^{left})-u(A_m^{left})$, namely:
$$\big\langle u(B_m^{left})-u(A_m^{left}), e_1\big\rangle =
length(\overline{p a_{j_1}\ldots B_m^{left}}) - length(\overline{p a_{i_1}\ldots A_m^{left}}).$$
Since there are no cuts in the interior of $P_m$, it follows that the
concatenation of the portion of $\sigma_r$ from $p$ to $A_m^{left}$
with $\xi_A$ and then with $\sigma_{r+1}$ from $B_m^{right}$ to $q$, cannot be
shorter than $\sigma_{r+1}$. Equivalently: 
$$length(\overline{p a_{i_1} \ldots A_m^{left}}) + length(\xi_A) \geq length(\overline{p a_{j_1} \ldots B_m^{right}}).$$ 
Also, concatenating $\sigma_{r+1}$ up to $B_m^{left}$ with
$\xi_B$ and then with the portion of $\sigma_r$ from $A_m^{right}$ to $q$, yields:
$$length(\overline{p a_{j_1} \ldots B_m^{left}}) + length(\xi_B) \geq length(\overline{p a_{i_1} \ldots A_m^{right}}).$$ 

The three last displayed bounds imply that:
\begin{equation}\label{eq6.1}
\begin{split}
& \big\langle u(B_m^{left})-u(A_m^{left}), e_1\big\rangle \in
[\alpha_w,\alpha_v], \\ & \mbox{where }\; \alpha_v \doteq
length(\xi_A) - length\big(B_m^{left}\ldots (\sigma_{r+1})\ldots B_m^{right}\big), \\ 
&\qquad \quad \, \alpha_w \doteq
length(\xi_B) - length\big(A_m^{left}\ldots (\sigma_{r})\ldots A_m^{right}\big)
\end{split}
\end{equation}

\smallskip

{\em Step 4.} Since $\xi_A$ has no
self-intersections, it divides $P_m$ into  two connected
components, and the endpoints $B_m^{left}, A_m^{right}$ of $\xi_B$
belong to the closures of the distinct components. Hence
$\xi_A\cap\xi_B\neq\emptyset$.

\begin{figure}[htbp]
\centering
\includegraphics[scale=0.6]{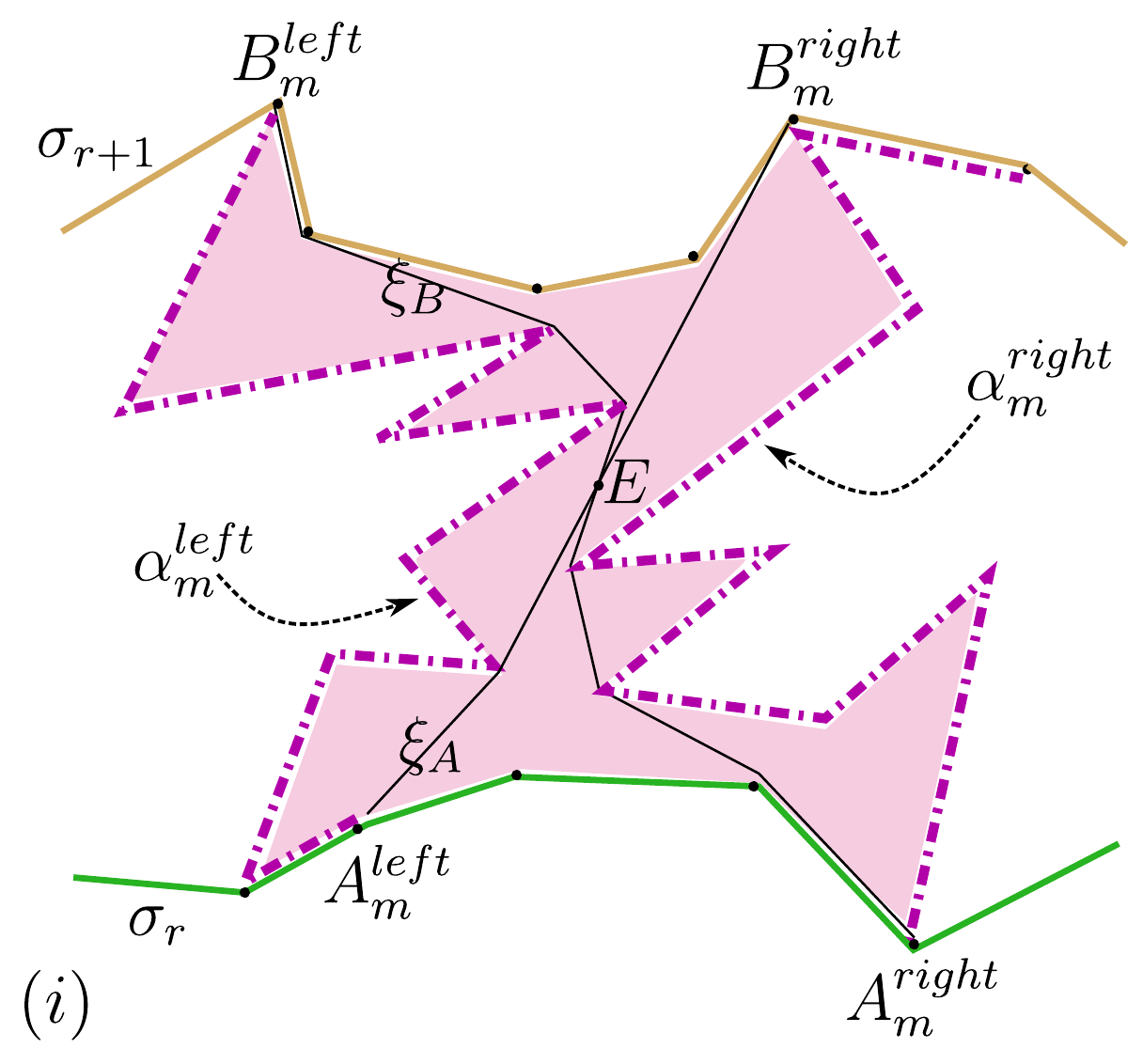}\qquad 
\quad \includegraphics[scale=0.6]{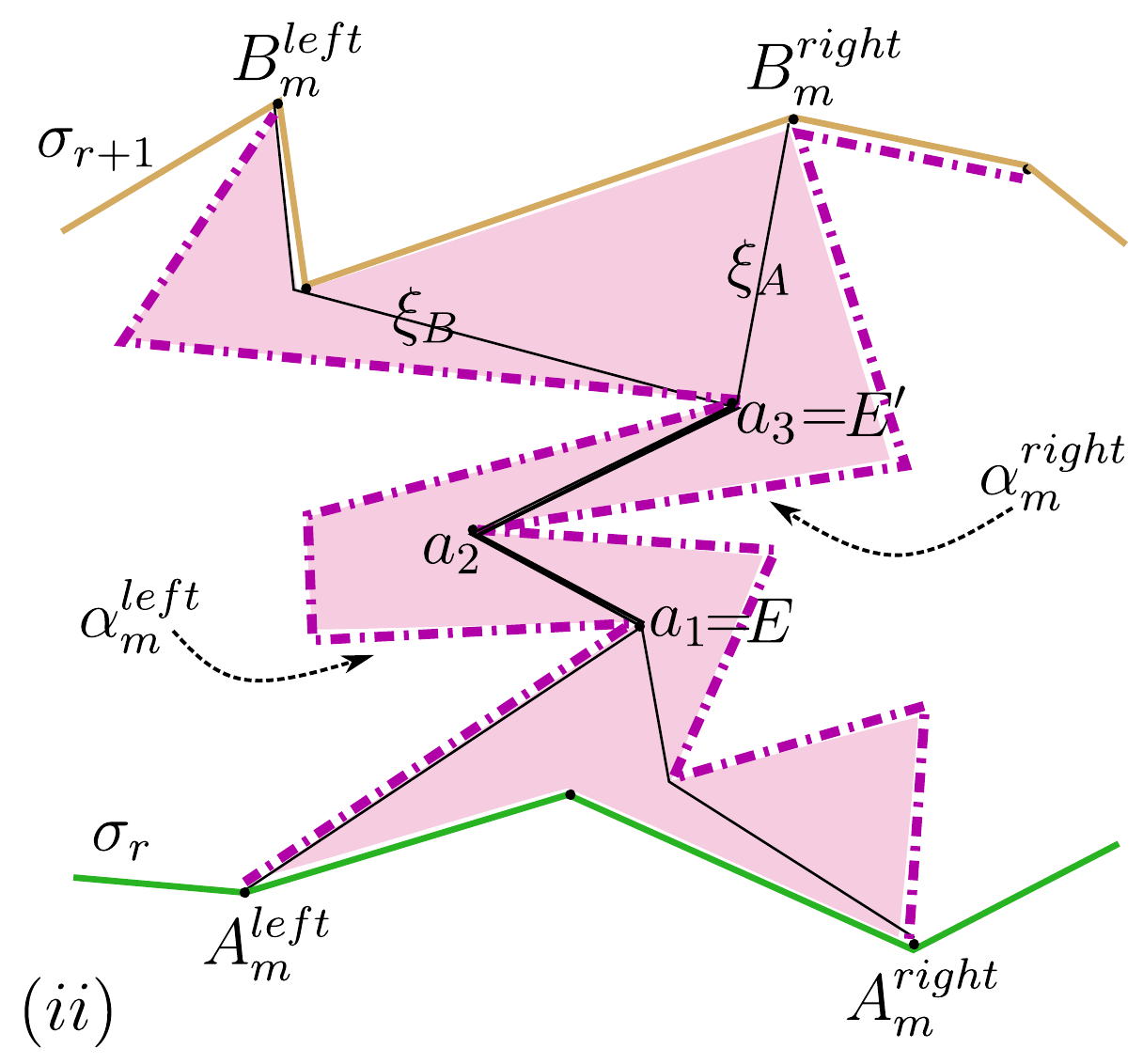}
\caption{{ Two types of regions $P_m$ (shaded) and the supporting polygonals
$\xi_A$, $\xi_B$ in the proof of Lemma \ref{lem6.2} Step 4: in (i) the
intersection $\xi_A\cap\xi_B$ consists of a single point $E$; in (ii) $\xi_A$
and $\xi_B$ intersect along the polygonal $\overline{a_1a_2a_3}$.}}
\label{Fig7.3}
\end{figure}

Define $C^{left}$ to be the vertex at which 
$\xi_A$ detaches itself from $\sigma_r$ ($C^{left}$ may be equal to $A_m^{left}$) and
$C^{right}$ to be the vertex where $\xi_B$ detaches from
$\sigma_r$. Similarly, we have the detachment vertices
$D^{left}\in\sigma_{r+1}\cap\xi_B$ and $D^{right}\in\sigma_{r+1}\cap\xi_A$. We observe in passing, that 
$\xi_A$ and $\xi_B$ used in the proof of Lemma \ref{lem6} are
precisely the portions of $\xi_A, \xi_B$ that we use presently,
between $C^{left}, D^{right}$ and $C^{right}, D^{left}$. We now argue
that $C^{left}$ precedes or equals to $C^{right}$ along $\sigma_r$
(with the usual order from $p$ to $q$). Assume by contradiction that
$C^{right}$ strictly precedes $C^{left}$ and call $\gamma$ the line spanned
by the edge interval of $\sigma_r$ that precedes $C^{left}$. By the minimizing
length of $\xi_A$, its portion right after the detachment from
$\sigma_r$ stays in the half-plane $S$ that is on the same side of
$\gamma$ as $\sigma_r\cap \partial P_m$. 

Assume first that $\alpha_m^{left}$ has no
intersection with $\xi_A$. In this case, $\xi_A$ is the straight
line from $C^{left}$ up to $D^{right}$, and $\alpha_m^{right}\subset
S$. Thus, $\alpha_m^{right}$ has no intersection with $\xi_B$
and both $\xi_B$ and $\alpha_m^{left}$ are contained in $S$, with
$\xi_B$ being a straight line from $C^{right}$ to $D^{left}$. 
The fact that both $D^{left}\neq D^{right}$ are in $\sigma_{r+1}\cap
S$ contradicts the concavity of $\sigma_{r+1}\cap \partial P_m$ and its
disjointness from $\sigma_r$. 

It hence follows that the polygonal $\xi_A$ must have a common vertex
with $\alpha_m^{left}$. Let $a$ be the first such vertex (in order
from $C^{left}$ to $D^{right}$; it necessarily belongs to $S$. Then
$D^{left}$ is contained in the region bounded by the concatenation of
the portion of $\alpha_m^{left}$ from $a$ to $A_m^{left}$, with
$\sigma_r$ from $A_m^{left}$ to $C^{left}$, with the straight segment
of $\xi_A$ from $C^{left}$ to $a$. Further, $B_m^{right}$ cannot
belong to the said region, unless $B_m^{right}=a$. This again
contradicts the convexity of $\sigma_{r+1}$ between $D^{left}$ and
$B_m^{right}$.

So indeed $C^{left}$ (respectively $D^{left}$) precedes or equals
$C^{right}$ (resp. $D^{right}$) along $\sigma_r$ (resp. $\sigma_{r+1}$). 

\smallskip

{\em Step 5.} We now make further observation about the supporting
polygonals $\xi_A, \xi_B$. Firstly, $\xi_A$ (respectively $\xi_B$)
have common vertices only with $\alpha_m^{left}$
(resp. $\alpha_m^{right}$) from its endpoint on $\sigma_r$ up to
its first intersection point $E$ with $\xi_B$ (resp. $\xi_A$). Indeed, the first time that $\xi_A$
encounters $\alpha_m^{right}$ it must also encounter $\xi_B$, and the
first time $\xi_B$ encounters $\alpha_m^{left}$ it must encounter
$\xi_A$ as well.

Secondly, the angles formed by $\xi_A$ or $\xi_B$ at these common vertices (up to $E$), interior
with respect to the polygon with vertices $C^{left}$, $E$, $C^{right}$
and with boundary along the appropriate portions of $\sigma_r$,
$\xi_A$ and $\xi_B$, are at least $\pi$. This fact is an easy
consequence of $\xi_A$, $\xi_B$ being shortest paths. 

\smallskip

The two symmetric statements to those made above,  are likewise valid for the portions of $\xi_A$,
$\xi_B$ from their endpoints on $\sigma_{r+1}$ up to their respective
first intersection point $E'$ (there may be $E'=E$). 
Our last statement is that between $E$ and $E'$, the polygonals
$\xi_A, \xi_B$ coincide. For otherwise, $\xi_A$ and $\xi_B$ would be
non-intersecting between some common vertices $\bar E\neq\bar E'$,
hence bounding a polygon whose at least one angle other than $\bar E,\bar E'$ would
be less than $\pi$ (with respect to the polygon's interior). This
would result in $\alpha_m^{left}$ or $\alpha_m^{right}$ having a
vertex at the indicated angle, and hence necessarily intersecting the
opposite boundary portion (of the polygon), which is a contradiction.

The two scenarios of $\xi_A\cap\xi_B$ consisting of a single
intersection point, and of a polygonal with vertices in $V$ (we recall
that $V$ is the set of vertices of the graph $G$) are
presented in Figure \ref{Fig7.3}.

\smallskip

{\em Step 6.} In this and the next Step we assume that: 
\begin{equation}\label{cas1}
E=E'\not\in V.
\end{equation}
We claim that there exists a folding pattern on $P_m$ such that:
\begin{itemize}
\item[(i)] the polygonal $\overline{C^{right}a_{h_1}\ldots
  a_{h_z}D^{left}}\subset \xi_B$ has its image contained in some straight line
$\bar\xi_B$, 
\item[(ii)] the same polygonal has its image length unchanged from the original length,
\item[(iii)]  the images of the portions of geodesics
$\sigma_r\cap \partial P_m$,  $\sigma_{r+1}\cap \partial P_m$ are their rigid
motions.
\end{itemize}

\smallskip

For this construction, we write
$E\in\overline{a_{h_t}a_{h_{t+1}}}$ and first consider the polygonal
$\overline{Ea_{h_{t+1}}\ldots a_{h_z}D^{left}}$  (see Figure
\ref{Fig7.6} (i) with $t=2$). We start at the vertex $a_{h_{t+1}}$ and perform the simple fold
of the half-line that extends $\overline{a_{h_{t+1}}a_{h_{t+2}}}$
beyond $a_{h_{t+1}}$, intersected with $P_m$, onto the line spanned by
the edge $\overline{a_{h_t}a_{h_{t+1}}}$. In doing so, we
take advantage of the fact that both indicated lines intersect
$\alpha_m^{right}$ before they may intersect $\sigma_{r+1}$. Next, we
fold by bisecting the angle at the vertex $a_{h_{t+2}}$: the half-line
extending $\overline{a_{h_{t+3}}a_{h_{t+2}}}$ beyond $a_{h_{t+2}}$
becomes the part of half-line that extends the (previously modified)
segment $\overline{a_{h_{t+1}}a_{h_{t+2}}}$. We continue in this
fashion, until we align $\overline{D^{left}a_{h_z}}$ with (previously
modified) $\overline{a_{h_{z-1}}a_{h_z}}$. 
Similarly, we consider the polygonal $\overline{C^{right}a_{h_1}\ldots
a_{h_t}E}$, and starting from $a_{h_t}$ we align $\overline{a_{h_t}a_{h_{t-1}}}$ with
$\overline{a_{h_t}a_{h_{t+1}}}$. The final fold in this construction
is that of half-line extending $\overline{C^{right}a_{h_1}}$ beyond
$a_{h_1}$, intersected with $P_m$, onto the line spanned by
$\overline{a_{h_{1}}a_{h_{2}}}$. As a result, the portion of $\xi_B$
between $C^{right}$ and $D^{left}$ has been straightened onto the line
$\bar\xi_B$ spanned by the segment $\overline{a_{h_t}a_{h_{t+1}}}$.

We now fold both portions of geodesics $\sigma_r\cap \partial P_m$ and
$\sigma_{r+1}\cap \partial P_m$ onto $\bar \xi_B$, using their
concavity, in the same manner as was done in the proof of Lemma \ref{lem6} in the simplified context of section
\ref{sec23}. By a further rotation we may exchange $\bar\xi_B$ into
$\mathbb{R}e_1$ and denote the resulting isometric immersion of $P_m$ by
$w$. Recalling the notation in (\ref{eq6.1}), it directly follows that:
$$ w(B_m^{left})-w(A_m^{left})= \alpha_w e_1. $$
Similarly, by folding $(\sigma_r\cup\sigma_{r+1})\cap \partial P_m$ onto the line $\bar\xi_A$ obtained as the
straightening of the portion of $\xi_A$ from $C^{left}$ to $D^{rigth}$, we get
an isometric immersion $v$ of $P_m$ with the property that:
$$ v(B_m^{left})-v(A_m^{left})= \alpha_v  e_1. $$

\smallskip

\begin{figure}[htbp]
\centering
\includegraphics[scale=0.6]{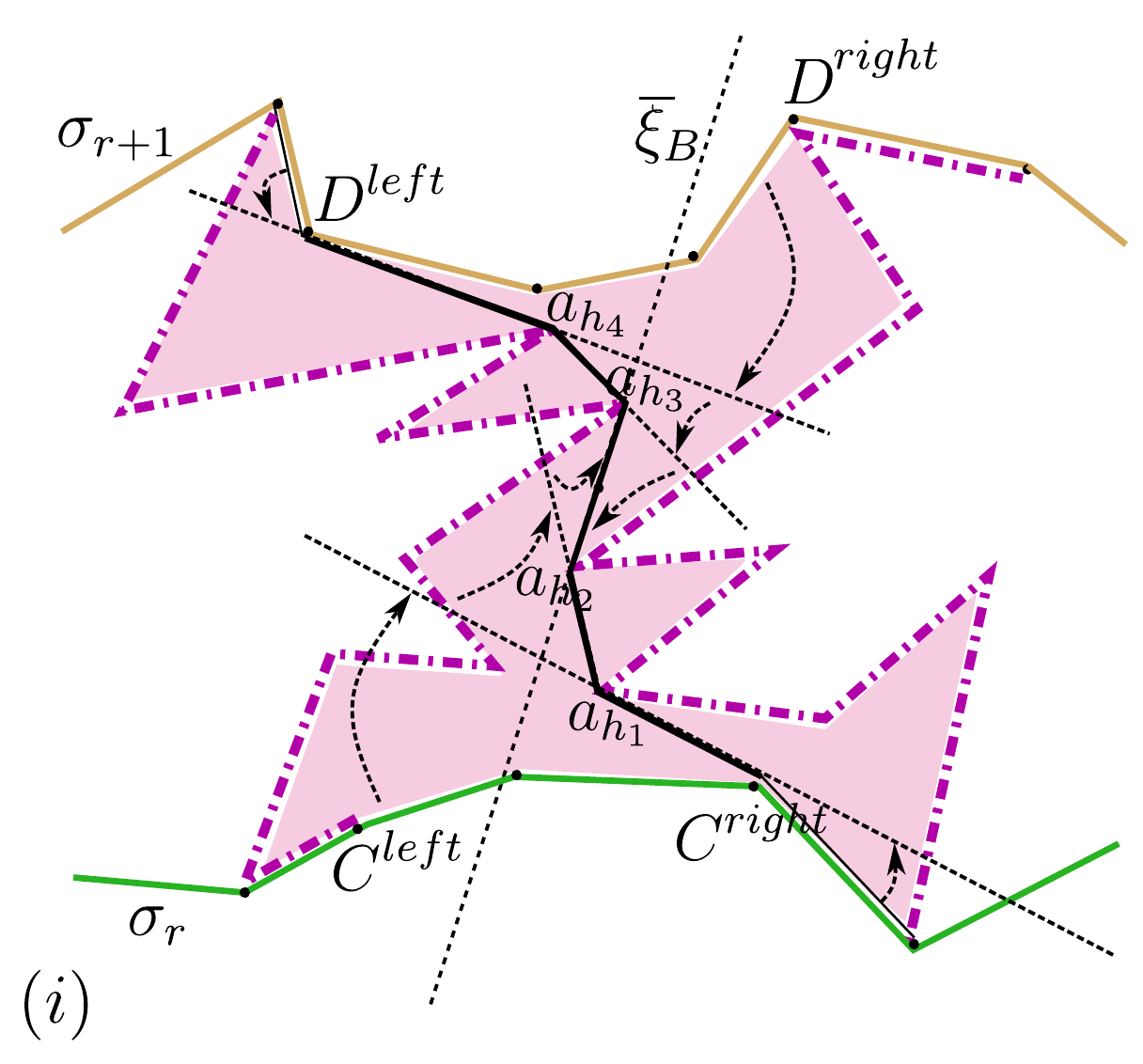}\qquad 
\quad \includegraphics[scale=0.6]{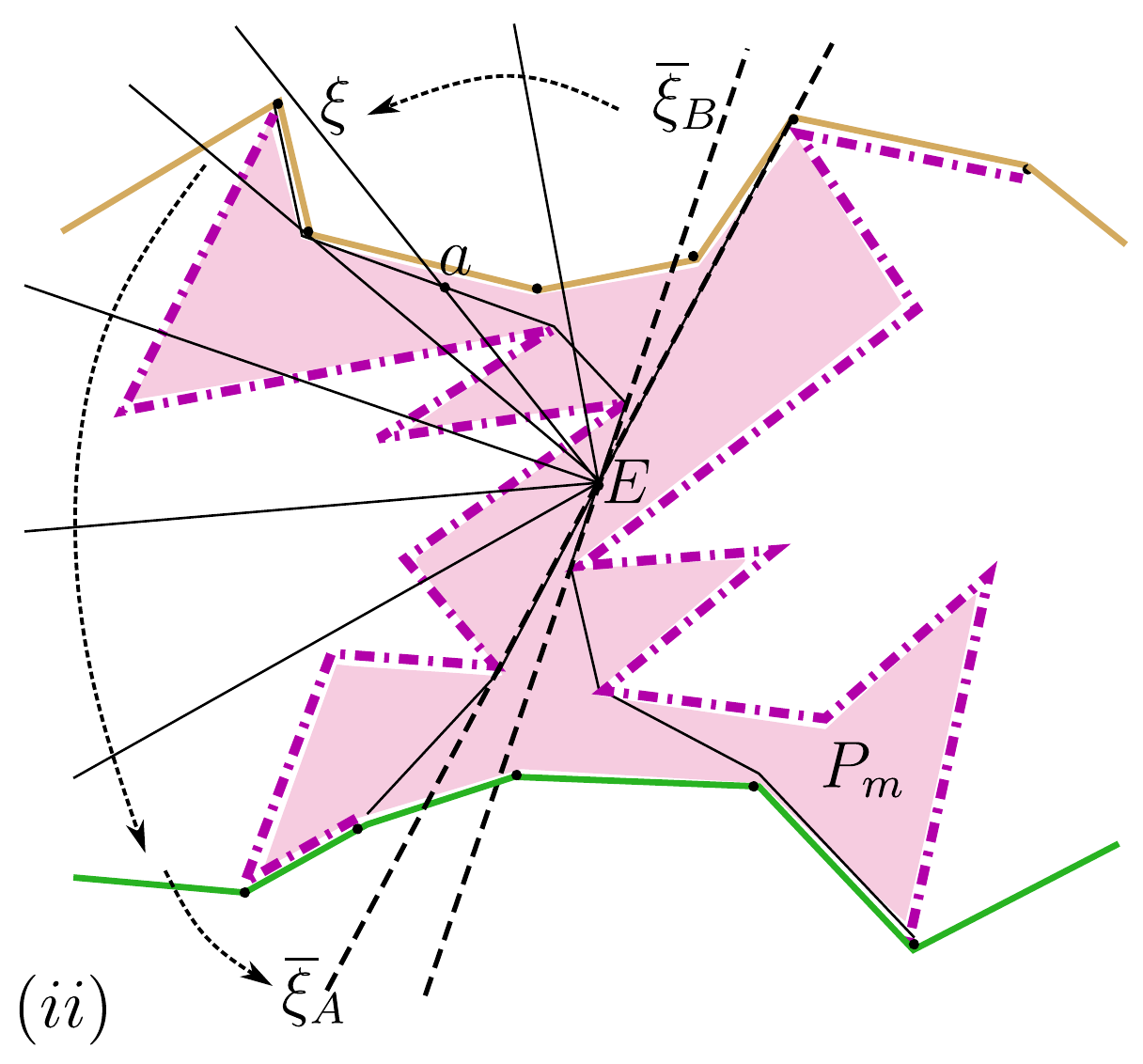}
\caption{{Elements of the proof of Lemma \ref{lem6.2}: diagram (i)
    depicts construction of the isometric immersion $w$ on the region
    $R_r$ as in Figure \ref{Fig7.3} (i). The arrows indicate the consecutive folds and
    the resulting straightenings of the intermediate polygonal
    $\overline{C^{right}a_{h_1}a_{h_2}a_{h_3}a_{h_4}D^{left}}\subset
    \xi_B$ to the line $\bar\xi_B$, and the projections of the boundary
    portions of $\sigma_r$, $\sigma_{r+1}$ onto $\bar\xi_B$
    in Step 6; diagram (ii) depicts the direction of rotating from $\bar\xi_B$
    to $\bar\xi_A$, through the intermediate direction lines $\xi$ in
    Step 7. The intersection point of the given half-line $\xi^1$ with the polygonal
    $\xi_B$ is called $a$ (provided it exists).}}
\label{Fig7.6}
\end{figure}

\smallskip

{\em Step 7.} Consider the family of lines $\{\xi\}$ obtained by
rotating $\bar\xi_B$ around $E$ onto $\bar\xi_A$. The direction of
rotation (see Figure \ref{Fig7.6} (ii)) is so that the half-line from
$E$ through $a_{h_{t+1}}$ gets rotated onto the half-line from $E$ to
the vertex on $\xi_A$ that is closest to $E$ between $E$ and
$C^{left}$, without passing through $D^{right}$ and $C^{right}$ along
the way. For each such line $\xi$ we will describe the folding
and the resulting isometry $u_\xi$ on $P_m$, with the property that the function:
$$\xi\mapsto u_\xi(B_m^{left}) - u_\xi(A_m^{left})$$
is continuous and that $u_{\bar\xi_B} = w$,  $u_{\bar\xi_A} = v$. By a
further rotation, we may map $\xi$ onto $\mathbb{R}e_1$ and hence
conclude that the scalar function
$\xi\mapsto \langle u_\xi(B_m^{left}) - u_\xi(A_m^{left}), e_1\rangle$
attains all values in the interval $[\alpha_w, \alpha_v]$. In virtue
of (\ref{eq6.1}), this will end the proof of Lemma \ref{lem6.2} under
the assumption (\ref{cas1}).

\smallskip

Fix $\xi$ as above and denote by $\xi^1$ the half-line emanating from
$E$ which is the rotated image of the half-line obtained by extending
$\overline{E a_{h_{t+1}}}\subset \bar\xi_B$ beyond $a_{h_{t+1}}$. We
also denote $\xi^2=\xi\setminus \xi^1$. Now, if $\xi^1$ intersects the
portion of the polygonal $\xi_B$ between $E$ and $D^{left}$ (we will
refer to this intersection point by calling it $a$) , we utilize the same
folding construction as in Step 6, but we replace the portion of
$\xi_B$ between $E$ and $a$ by the segment $\overline{Ea}\subset
\xi$. Observe that $\overline{Ea}$ must intersect the boundary
$\alpha_m^{left}$ of $P_m$ before it reaches $a$. Thus there exists a
simple fold which results in rotating (around $a$) of the edge of
$\xi_B$ containing $a$, onto $\xi$, and in such a way that the
position of $a$ remains unchanged, and that $\sigma_{r+1}\cap \partial P_m$ is
also only transformed via a rigid motion. We then continue the
straightening procedure of $\xi_B$ beyond $a$ as before. 

On the other hand, if the first intersection point $a$ of $\xi^1$ with $\xi_B$ occurs
between $D^{left}$ and $B_m^{left}$, we again take advantage of the
fact that the open segment $\overline{Ea}$ must intersect
$\alpha_m^{left}$; this allows for a single simple fold which rotates
the segment edge of $\xi_B$, to which $a$ belongs (this edge must be
contained in $\xi_{r+1}$) around $a$ and onto $\xi$. 
In both so far described cases, the geodesic portion $\sigma_{r+1}\cap
\partial P_m$ may be subsequently folded onto $\xi$, due to its concavity and
the fact that its one point ($D^{left}$ in the former case, $a$ in the
latter) already belongs to $\xi^1$. 

In the third case when $\xi^1$ has no intersection with $\xi_B$ beyond
$E$ (hence $\xi^1\cap \sigma_{r+1}\cap \partial P_m=\emptyset$), we utilize the
construction from proof of Lemma \ref{lem6}. This entails identifying
the line $\gamma$ that is parallel to $\xi$ and supporting to $\sigma_{r+1}\cap \partial P_m$.
We then first fold $\sigma_{r+1}\cap \partial P_m$ onto $s$, and then fold $\gamma$
onto $\xi$. This, again, can be done without altering $\sigma_{r}\cap
\partial P_m$ beyond possibly applying a rigid motion to it, because both lines $\xi^2$ and
$\gamma$ intersect $\alpha_m^{right}$ before they possibly intersect
$\sigma_r\cap \partial P_m$. 

\smallskip

Rotating $\xi^1$ further, we have it eventually pass through
$A_m^{left}$, then $C^{left}$, then intersect the polygonal
$\xi_A\setminus \sigma_{r}$, and finally coincide with the appropriate
half-line in $\bar\xi_A$. In each of these listed scenarios, we perform the corresponding (in the reverse
order of appearance)  folding construction relative to the polygonal
$\xi_A$ rather than $\xi_B$. In the same fashion, we define the
folding patterns relative to the half-line $\xi^2$. This concludes the definition of
each $u_\xi$ in case (\ref{cas1}).

\smallskip

\begin{figure}[htbp]
\centering
\includegraphics[scale=0.6]{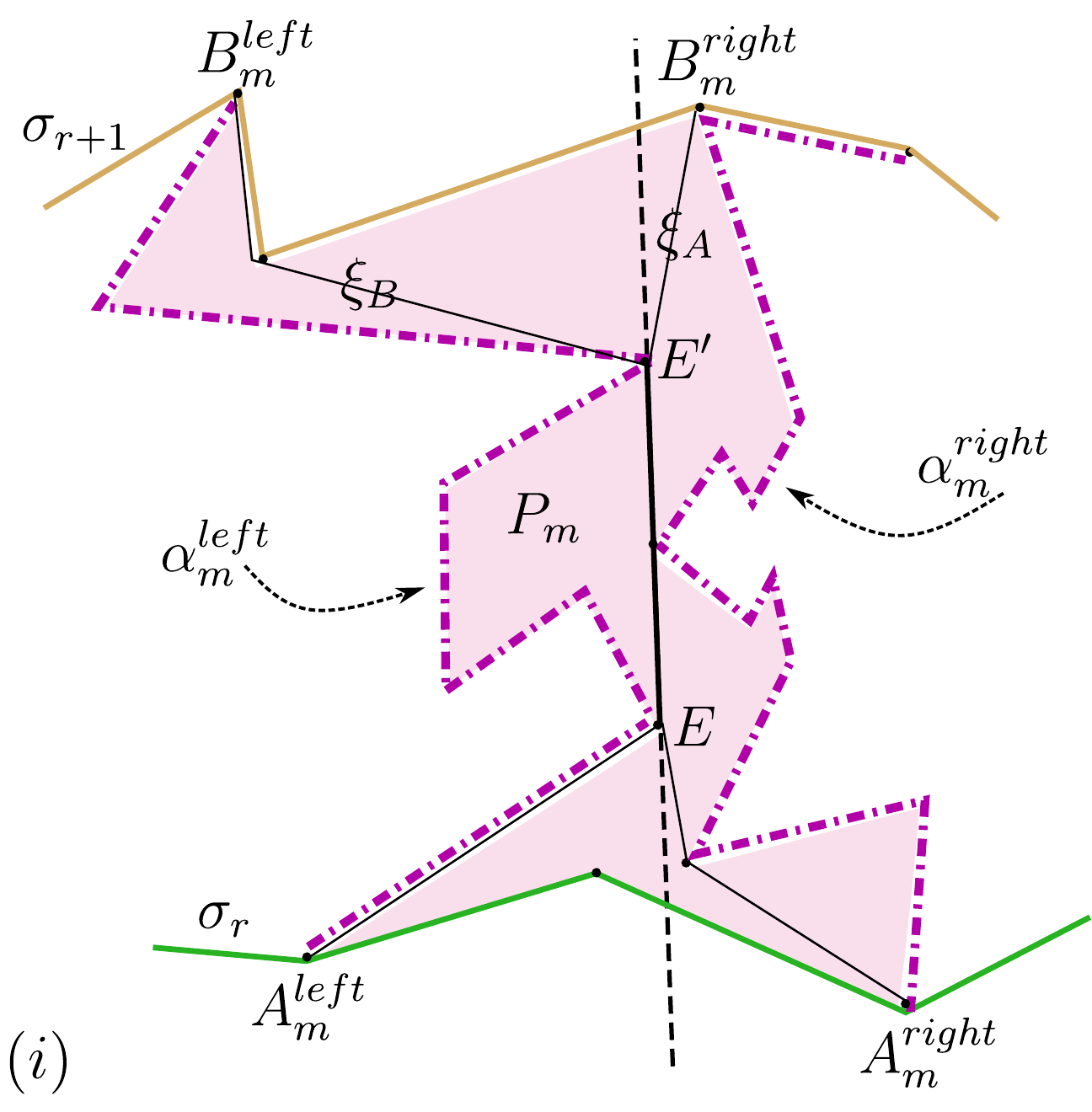}\qquad 
\includegraphics[scale=0.6]{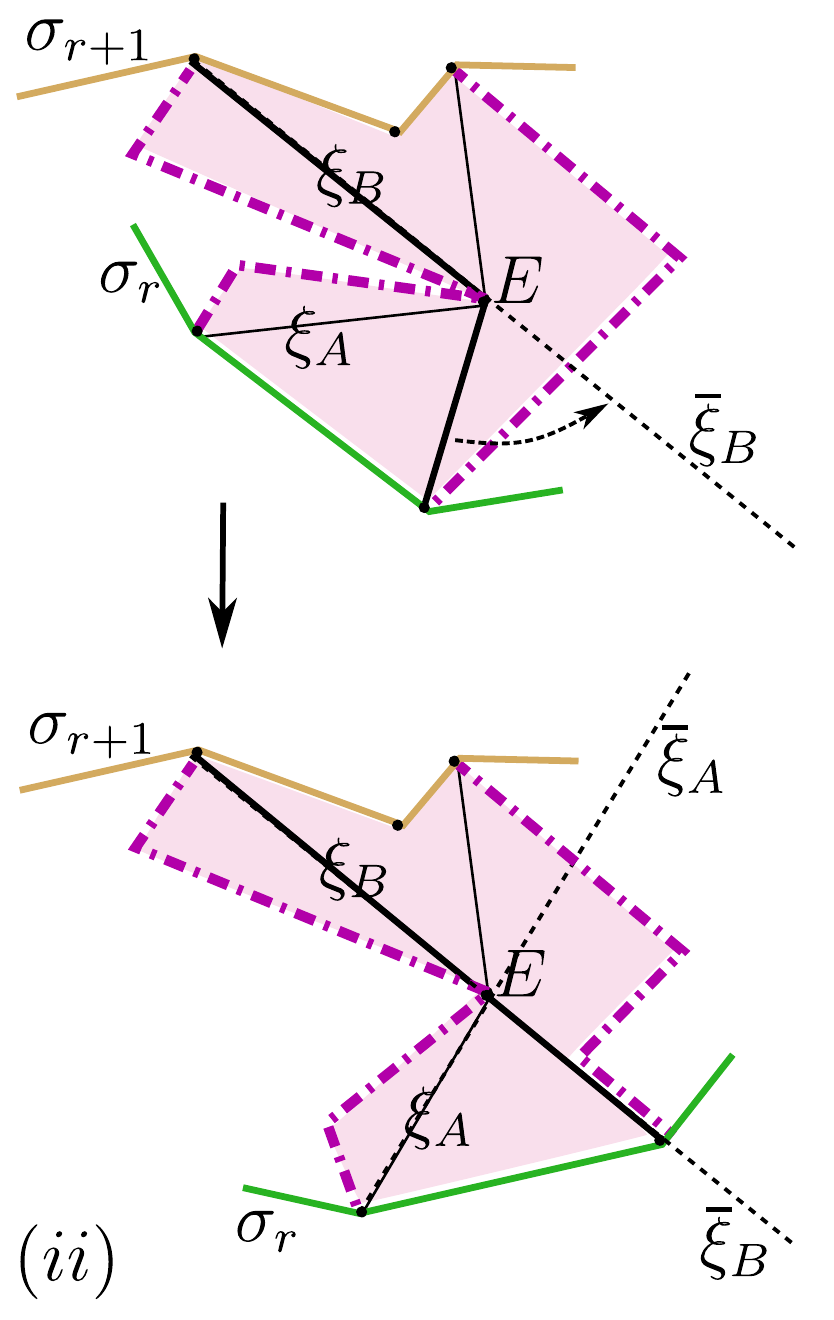}
\caption{{Construction in Step 8 of the proof of Lemma
 \ref{lem6.2}: (i) depicts the result of the initial folding,
 straightening $\xi_A\cap\xi_B$ into segment $\overline{EE'}$, applied
 to $P_m$ in Figure \ref{Fig7.3} (ii); diagram (ii) indicates the
 initial folding in case $E=E'$.}}
\label{Fig7.7}
\end{figure}

{\em Step 8.} Note that $E\neq E'$ implies $E, E'\in
V$ (see Figure \ref{Fig7.3} (ii)). In this Step we assume that:
\begin{equation}\label{cas2}
E=E'\in V \quad \mbox{or} \quad E\neq E'.
\end{equation}
When $E\neq E'$, we first perform several simple folds which straighten the polygonal $\xi_A\cap\xi_B$ 
into a segment with endpoints that we continue to denote $E$ and $E'$,
and such that the line spanned by the new $\overline{EE'}$ enters each of the two angles between
the pairs of distinct edges in $\xi_A$ and $\xi_B$, emanating from $E$ and
$E'$ (see Figure \ref{Fig7.7} (i)). As a result, the new polygonals
$\xi_A$ and $\xi_B$ have the same convexity properties as in case
(\ref{cas1}). This allows for applying the folding construction in
Step 7, relative to the center point $(E+E')/2$ % instead of $E$, 
and the projection lines $\bar\xi_B=\bar\xi_A$ spanned by $\overline{EE'}$.

\smallskip

When $E=E'$, we first perform one simple fold at the vertex $E$, which
either: (i) makes the two edges of $\xi_B$ with common vertex $E$
collinear, and keeps the angle (internal to $P_m$) between the two edges of $\xi_A$
adjacent to $E$ not smaller than $\pi$; or (ii) makes the two edges of $\xi_A$ with common vertex $E$
collinear, and keeps the angle between the two edges of $\xi_B$
adjacent to $E$ not smaller than $\pi$. In what follows we will assume,
without loss of generality, the former scenario as in Figure \ref{Fig7.7} (ii).

We let $\bar\xi_B$ to be the line spanned by the segment of $\xi_B$ passing
through $E$ and $\bar\xi_A$ to be spanned by the segment with vertex
$E$ and the successive vertex along $\xi_A$ towards $C^{left}$. We also
call $\xi_A^2$ the half-line from $E$ through the successive vertex
along $\xi_A$ towards $D^{right}$. We now apply the construction from
Steps 6 and 7, where we rotate the line $\bar\xi_B$ onto $\bar\xi_A$
around $E$ and perform a family of foldings onto each intermediate
line $\xi$. In case $\xi_A^2\not\subset \bar\xi_A$, the same
construction is applied to each half-line emanating from $E$ and
intermediate to $\bar\xi_A$ and $\xi_A^2$, completed by an extra
simple fold at $E$ that aligns the said lines. This ends the
definition of each $u_\xi$ in case (\ref{cas2}).

\smallskip

{\em Step 9.} The final step is to construct $u$ on $P_s$. This can be
done by the same folding technique as in Step 1. We also get: $\langle
u(B_s^{left}) - u(A_s^{left}), e_1\rangle =
length(\overline{A_s^{left}\ldots a_{j_l}q_1}) - 
length(\overline{B_s^{left}\ldots a_{i_k}q_1})$, because $length(\sigma_r) =
length(\sigma_{r+1})$.  This ends the proof of Lemma \ref{lem6.2}.
\end{proof}

\section{Proof of Theorem \ref{thm2}. Step 4: isometric immersion on the
  exterior region. A counterexample when $p,q\not\in \partial\Omega$}\label{sec24}

In this section, we first construct an isometric immersion $u$ on the remaining region $R_0$. 

%we complete the proof of Theorem \ref{thm2} under a
%simplifying assumption of the exterior region $R_0$ being free of
%cuts. This condition holds automatically when all trees in $G$ are segments $l\in E$, as in
%section \ref{sec23}. We take $p\neq q\in \Omega\setminus L$ and
%construct an isometric immersion $u$ on $R_0$. 

\begin{lemma}\label{lem7}
Assume (\ref{S}) and (\ref{S1}). If $p,q\in\partial\Omega$ then there exists a
continuous, piecewise affine isometric immersion $u$ of $R_0$ into $\R^3$, satisfying:
$$u(p) = 0, \quad u(q)=length(\sigma_1)e_1, \quad u(\sigma_1)=u(\sigma_{N}) = I.$$
\end{lemma}
\begin{proof}
By Lemma \ref{lem5} we have:  $R_0\cap L=\emptyset$ 
and both the least and the greatest geodesics
$\sigma_1, \sigma_N$ are convex, i.e. the region $\Omega\setminus R_0$
has all the (internal) angles at the vertices distinct from $p$, $q$, not greater than $\pi$.
Indeed, consider an intermediate vertex $A\not\in\{p,q\}$ of
$\sigma_{min}$. If the internal angle at $A$ was strictly larger than
$\pi$, then the cut $l =\overline{AB}\in E$ emanating from $A$ would have to point
inside $R_0$, as otherwise $\sigma_{min}$ could be shortened, contradicting the fact that it is a
geodesic. The argument for $\sigma_{max}$ is similar.
One can now apply the usual sequence of simple folds to obtain $u$ on $R_0$.
\end{proof}

\medskip

The proof of Theorem \ref{thm2} is now complete. Note that the
constructed isometric immersion $u$ consists exclusively of planar
folds and returns the image that is a subset of $\R^2$.

\begin{figure}[htbp]
\centering
\includegraphics[scale=0.6]{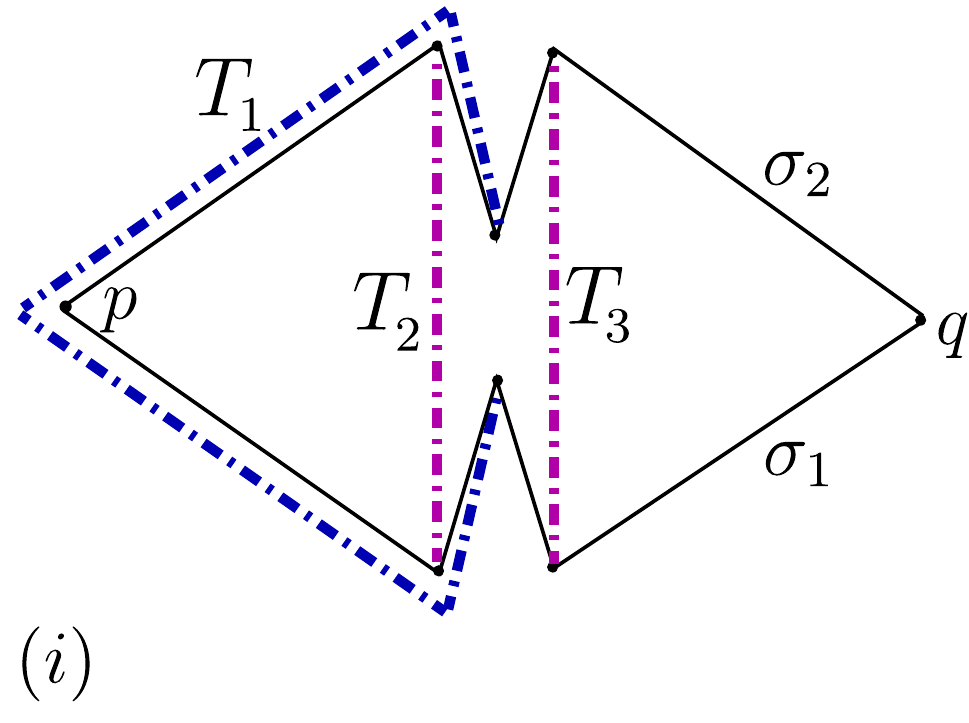} \qquad \quad
\includegraphics[scale=0.6]{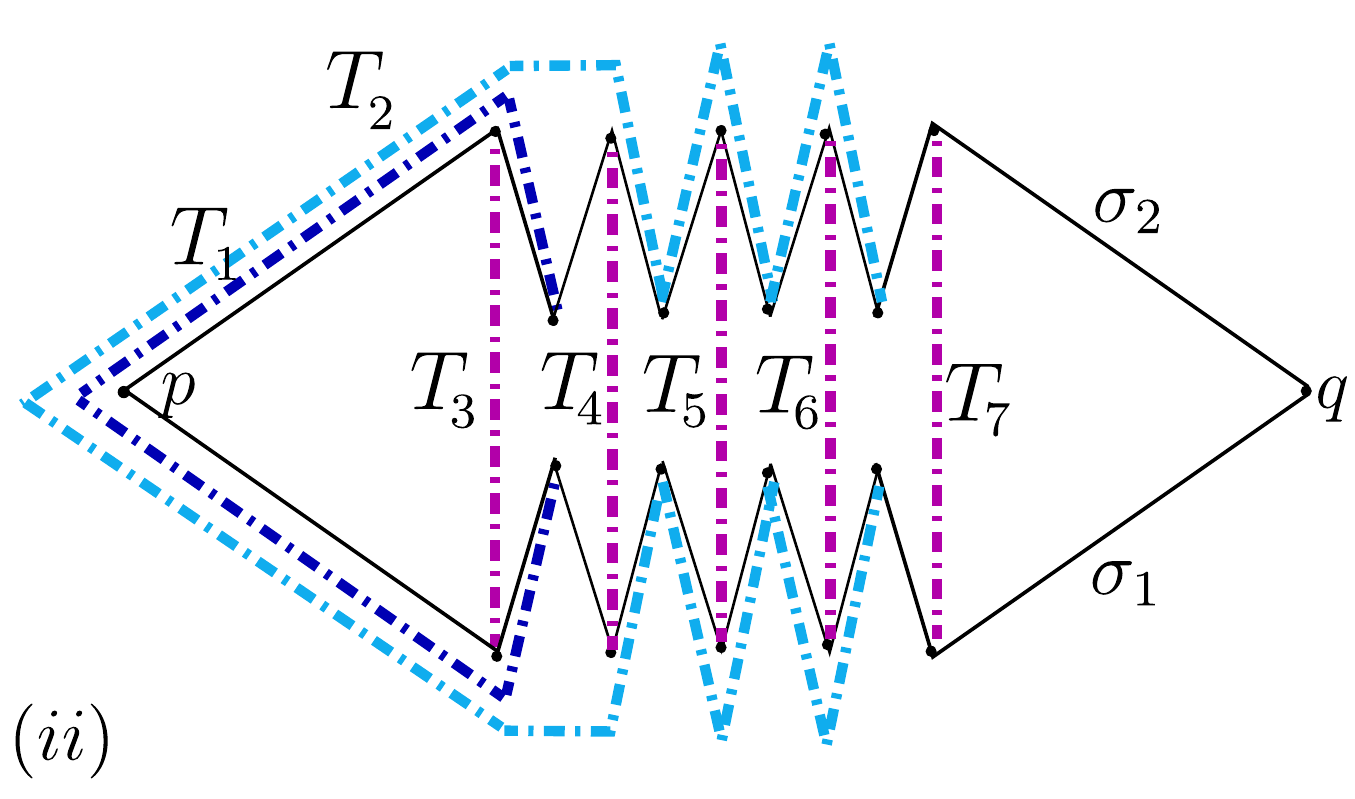} \\
\includegraphics[scale=0.6]{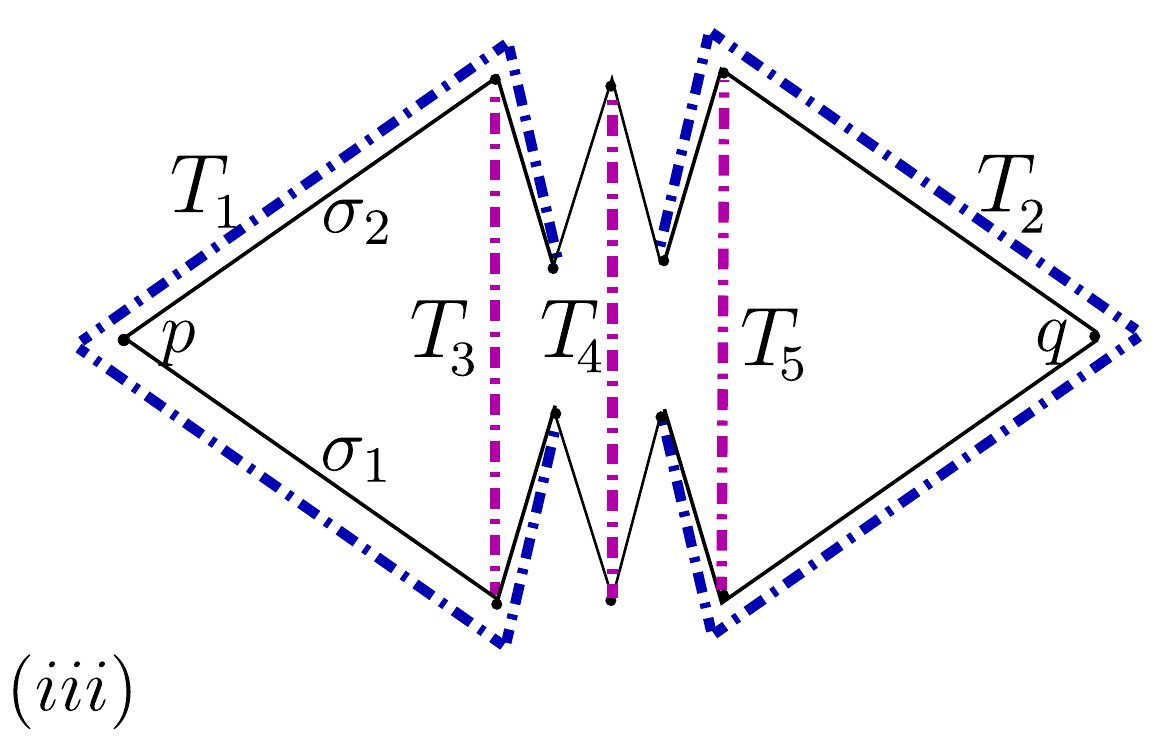} \qquad \quad
\includegraphics[scale=0.6]{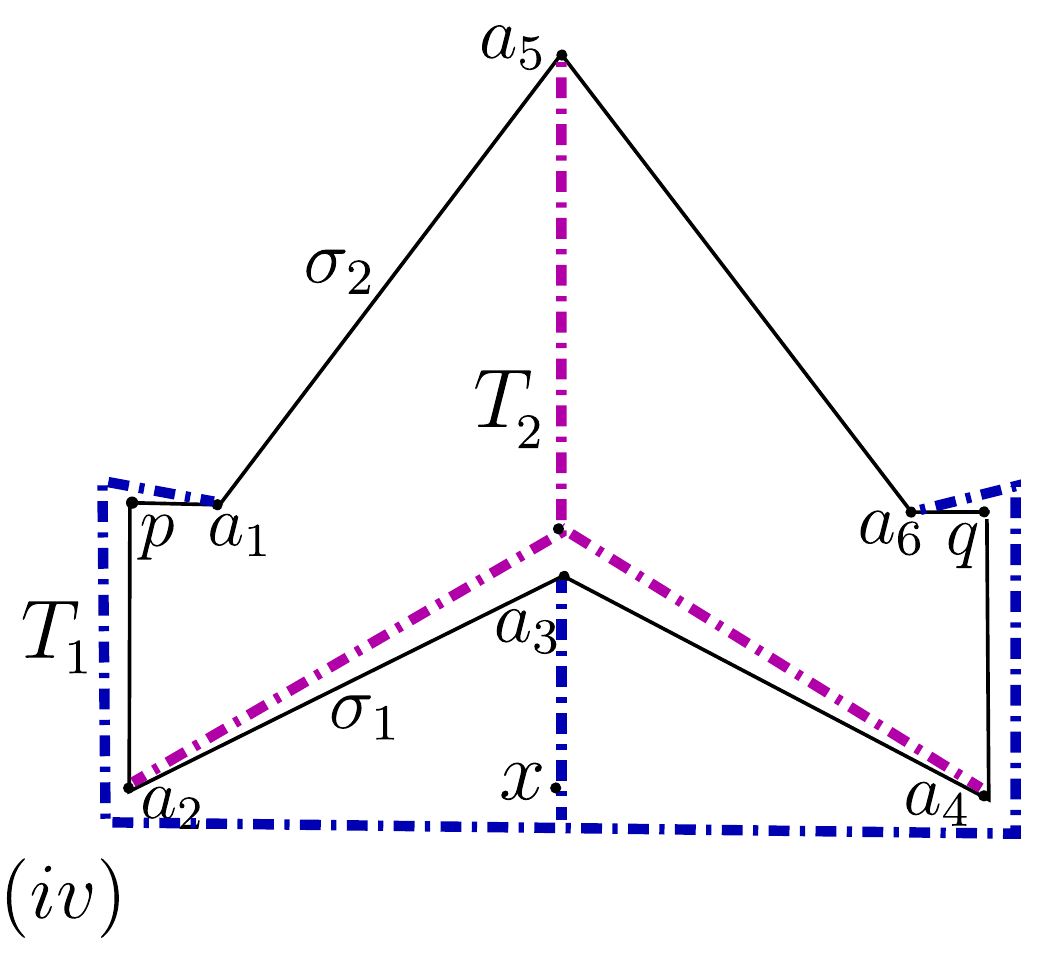} 
\caption{Examples of minimal configurations with the external region $R_0$
containing non-sealable cuts: in (i) the graph $G$ consists of three trees
    $\{T_i\}_{i=1}^3$. There are two geodesics $\sigma_1,\sigma_2$ from $p$ to $q$ in
    $\Omega\setminus L$, a single internal region $R_1$ (which
    contains $T_2$, $T_3$) and the external region $R_0$ which contains
  $T_1$. Note that $T_1$ has leaves on both $\sigma_1$ and $\sigma_2$;
in (ii) $R_0$ contains two ``nested'' trees $T_1, T_2$, while there are five
more trees in $R_1$; in (iii) $R_0$ contains two trees $T_1, T_2$
(not ``nested''); in (iv) we assume that $|\overline{pa_1}| < |\overline{a_2a_3}| -
|\overline{a_2x}|$, for minimality. There
is a single tree $T_1$ in $R_0$ and a single tree $T_2$ in
$R_1$. The two geodesics are: $\sigma_1=\overline{pa_2a_3a_4q}$ and $\sigma_2=\overline{pa_1a_5a_6q}$.}
\label{Fig9.1}
\end{figure}

\bigskip 

The fact that $R_0\cap L=\emptyset$ is directly related to the
assumption that $p,q\in\partial\Omega$. Indeed, examples in Figure \ref{Fig9.1}
show there may be (even multiple) external trees, necessarily with
vertices on both $\sigma_1$ and $\sigma_N$, when the said assumption is removed.
This type of configuration may also be used to show
that an isometric immersion $u$ of $\Omega\setminus L$ with the property
that the Euclidean distance between $u(p)$ and $u(q)$ equals the
geodesic distance from  $p$ to $q$ in $\Omega\setminus L$, may in
general not exist. 
%when at least one of points $p,q$ lies in the interior of $\Omega$.

\begin{figure}[htbp]
\centering
\includegraphics[scale=0.6]{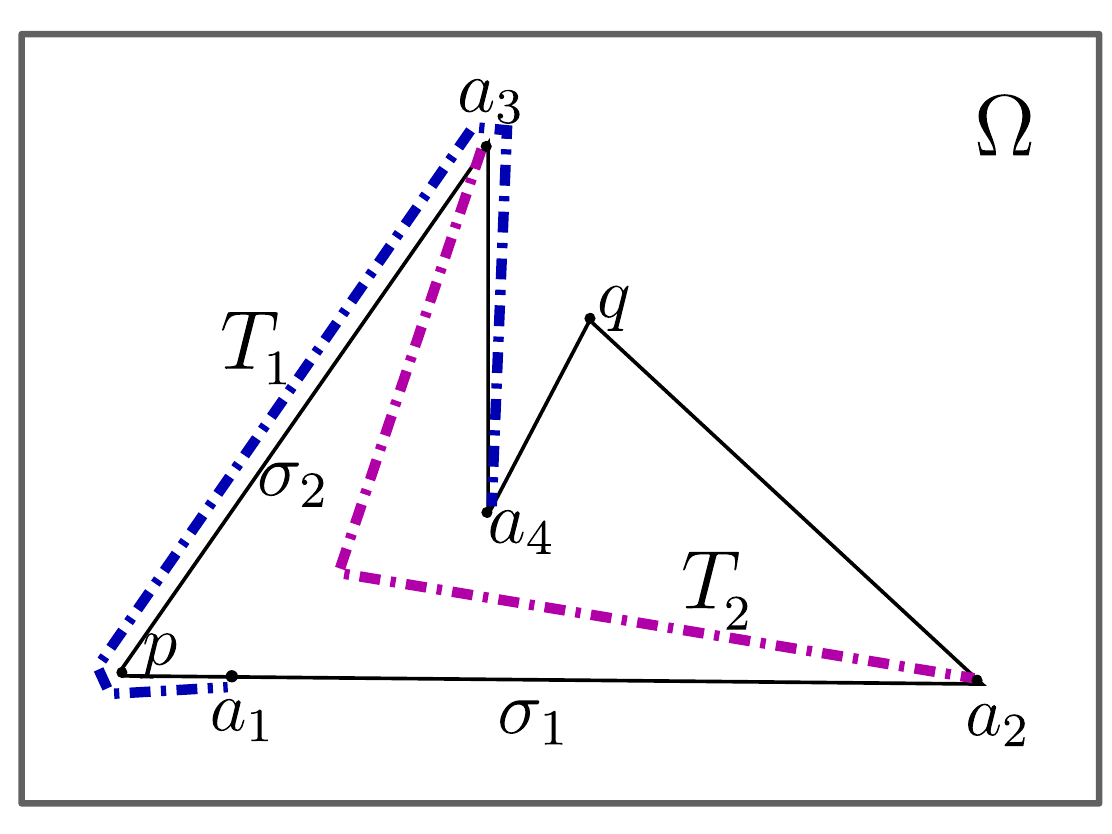} 
\caption{A configuration of $G, \Omega$ and $p,q\in\Omega$ for which
  the conclusion of Theorem \ref{thm1} fails. The set of vertices of
  $G$ is: $V=\big\{a_1=(0,0), a_2=(1,0), a_3=(\frac{1}{\sqrt{2}}-c, \frac{1}{\sqrt{2}}),
  a_4=(\frac{1}{\sqrt{2}}-c, \alpha)\big\}$ and $p=(-c,0),
  q=(\frac{1}{2},\frac{1}{2})$. The set of cuts is:
  $L=T_1\cup T_2$, where $T_1$ is the single exterior tree and $T_2$ is 
the single interior tree. For every $\frac{1}{\sqrt{2}} - \frac{1}{2}<c<\frac{1}{\sqrt{2}}$ there exists
  $0<\alpha<\frac{1}{2}$ so that there are two geodesics $\sigma_1 = \overline{pa_1
a_2q}$, $\sigma_2 = \overline{pa_3a_4q}$ satisfying:
  $length(\sigma_1)=length(\sigma_2)=1+\frac{1}{\sqrt{2}}+c$. When
  additionally $c<\sqrt{2}-1$, then the above configuration is minimal. }
\label{Fig9.2}
\end{figure}

Consider the example in Figure \ref{Fig9.2}.
It is easy to check that $\mbox{dist}_{\Omega\setminus L}(p,q) =
c+1+\frac{1}{\sqrt{2}}=length(\sigma_1)=length(\sigma_2)$, when:
$$\alpha=\alpha(c) = \frac{1-\frac{1}{\sqrt{2}}+(1-\sqrt{2})c}{2c+1}.$$
Also, the constraints $\frac{1}{\sqrt{2}}-c<\frac{1}{2}$ and
$0<\alpha<\frac{1}{2}$ (implying that the interior region $R_1$ has exactly
the shape indicated in Figure \ref{Fig9.2}) hold, in particular, when taking:
\begin{equation}\label{mm}
\frac{1}{\sqrt{2}}-\frac{1}{2}<c<\sqrt{2}-1.
\end{equation}
The minimality of the configuration $L=T_1\cup T_2$ is guaranteed by
requesting that: $length(\overline{pa_1pa_3q})<\mbox{dist}(p,q)$, which
upon a simple calculation reduces to: $c<\sqrt{2}-1$, guaranteed in (\ref{mm}).

We now claim that there is no isometry $u$ of $\Omega\setminus
(T_1\cup T_2)$, which straightens the polygonal
$\overline{a_1a_2q}$. This is because otherwise there would be:
$$length(\overline{a_1a_2q})=|u(q)-u(a_1)|\leq\mbox{dist}_{\Omega\setminus R_1}(a_1,q).$$
However, the inequality above is violated when the tree $T_1$ approximates closely the
polygonal path $\overline{a_1pa_3q}$, in view of the bound $length(\overline{a_1pa_3q})<length(\overline{a_1a_2q})$
which again follows from (\ref{mm}).

\section{Discussion}\label{final}

Our two geometrical theorems are inspired by simple observations of
the mechanical response of a sheet of paper that has cuts in it,
valid only in the limit when the sheet is mapped to itself via a
piecewise (non-unique) affine map that is isometric to the plane. To
remove this non-uniqueness, we must account for the energetic penalty
of deforming a sheet of small but finite thickness, by bending it out
of the plane. When this physical fact is accounted for, a kirigamized
sheet will deform into a complex shape constituted of conical and
cylindrical structures glued together. 

Understanding the mechanics and
mathematics of these objects, while also solving the inverse problem
of how to design the number, size, orientation and location of the
cuts, remain open problems. And while we have limited ourselves to the
study of Euclidean case, our study naturally raises questions
about the nature and form of geodesics in non-Euclidean surfaces with
co-dimension one obstructions, and higher-dimensional generalizations
that might be relevant for traffic, fluid flow and stress transmission
in continuous and discrete geometries.

\end{document}